\documentclass[reprint,english,twocolumn,amsfonts,amssymb,amsmath]{revtex4-1}

\usepackage[T1]{fontenc}
\usepackage[latin9]{inputenc}
\usepackage[dvipsnames]{xcolor}
\usepackage{float}
\usepackage{fancybox}
\usepackage{calc}
\usepackage{mathtools}
\usepackage{url}
\usepackage{dsfont}
\usepackage{graphicx}
\usepackage{enumitem}
\usepackage[normalem]{ulem}
\usepackage{xr}
\usepackage{blkarray}
\usepackage{tikz,tikz-cd}
\usetikzlibrary{cd,decorations.pathreplacing,fit}

\usepackage[unicode=true,pdfusetitle,
 bookmarks=true,bookmarksnumbered=false,bookmarksopen=false,
 breaklinks=false,backref=false,colorlinks=true]
 {hyperref}
\hypersetup{linkcolor=blue, urlcolor=blue, citecolor=blue}

\makeatletter

\renewcommand{\d}{\mathrm{dist}}
\newcommand{\pa}{\mathrm{pa}}
\newcommand{\ddt}{\frac{\mathrm{d}}{\mathrm{d}t}}

\usepackage{stmaryrd}

\usepackage{trimclip}
\DeclareRobustCommand{\shortto}{\mathrel{\mathpalette\short@to\relax}}
\newcommand{\short@to}[2]{\mkern2mu\clipbox{{.5\width} 0 0 0}{$\m@th#1\vphantom{+}{\shortrightarrow}$}}

\newcommand\indep{\protect\mathpalette{\protect\indepT}{\perp}}
\def\indepT#1#2{\mathrel{\rlap{$#1#2$}\mkern2mu{#1#2}}}

\newcommand{\jsq}[1]{\textbf{\color{red}[JS: #1]}}
\newcommand{\jsr}[1]{}
\newcommand{\js}[1]{{\color{blue}#1}}

\newcommand{\de}{\mathrm{d}}
\makeatother

\newcommand{\bwq}[1]{\textbf{\color{RawSienna}[BW: #1]}}

\newcommand{\bwr}[1]{}

\usepackage{svg}
\graphicspath{ {./figures/motifsvgs/} }

\newcommand{\mota}{\vcenter{\hbox{\includesvg[inkscapelatex=false,height=.76em]{a}}}}
\newcommand{\motb}{\vcenter{\hbox{\includesvg[inkscapelatex=false,height=.76em]{b}}}}
\newcommand{\motc}{\vcenter{\hbox{\includesvg[inkscapelatex=false,height=.76em]{c}}}}

\newcommand{\motaNL}{\vcenter{\hbox{\includesvg[inkscapelatex=false,height=.57em]{a_nolab}}}}
\newcommand{\motab}{\vcenter{\hbox{\includesvg[inkscapelatex=false,height=.76em]{ab2}}}}
\newcommand{\motabNL}{\vcenter{\hbox{\includesvg[inkscapelatex=false,height=.57em]{ab2_nolab}}}}
\newcommand{\motbc}{\vcenter{\hbox{\includesvg[inkscapelatex=false,height=.76em]{bc}}}}
\newcommand{\motac}{\vcenter{\hbox{\includesvg[inkscapelatex=false,height=.76em]{ac}}}}
\newcommand{\motca}{\vcenter{\hbox{\includesvg[inkscapelatex=false,height=.76em]{ca}}}}
\newcommand{\motbd}{\vcenter{\hbox{\includesvg[inkscapelatex=false,height=.76em]{bd}}}}
\newcommand{\motabc}{\vcenter{\hbox{\includesvg[inkscapelatex=false,height=.76em]{abc}}}}
\newcommand{\motabd}{\vcenter{\hbox{\includesvg[inkscapelatex=false,height=.76em]{abd}}}}
\newcommand{\motcbd}{\vcenter{\hbox{\includesvg[inkscapelatex=false,height=.76em]{cbd}}}}
\newcommand{\motabctr}{\vcenter{\hbox{\includesvg[inkscapelatex=false,height=1.615em]{abctr}}}}
\newcommand{\motabcdtre}{\vcenter{\hbox{\includesvg[inkscapelatex=false,height=1.71em]{abcdtre}}}}

\newcommand{\motlatsix}{\vcenter{\hbox{\includesvg[inkscapelatex=false,height=1.71em]{lat_abcdef}}}}
\newcommand{\motlatfour}{\vcenter{\hbox{\includesvg[inkscapelatex=false,height=1.71em]{lat_abcd}}}}

\newcommand{\motabcNL}{\vcenter{\hbox{\includesvg[inkscapelatex=false,height=.57em]{abc_nolab}}}}
\newcommand{\mottr}{\vcenter{\hbox{\includesvg[inkscapelatex=false,height=1.33em]{abctr_nolab}}}}

\newcommand{\motsqi}{\vcenter{\hbox{\includesvg[inkscapelatex=false,height=1.425em]{sqi}}}}
\newcommand{\motsqii}{\vcenter{\hbox{\includesvg[inkscapelatex=false,height=2.1375em]{sqii}}}}
\newcommand{\motstr}{\vcenter{\hbox{\includesvg[inkscapelatex=false,height=1.425em]{str}}}}

\newcommand{\motS}{\vcenter{\hbox{\includesvg[inkscapelatex=false,height=.63em]{S}}}}
\newcommand{\motI}{\vcenter{\hbox{\includesvg[inkscapelatex=false,height=.63em]{I}}}}
\newcommand{\motIS}{\vcenter{\hbox{\includesvg[inkscapelatex=false,height=.63em]{IS}}}}
\newcommand{\motSS}{\vcenter{\hbox{\includesvg[inkscapelatex=false,height=.63em]{SS}}}}
\newcommand{\motII}{\vcenter{\hbox{\includesvg[inkscapelatex=false,height=.63em]{II}}}}

\newcommand{\motSSI}{\vcenter{\hbox{\includesvg[inkscapelatex=false,height=.5985em]{SSI}}}}
\newcommand{\motSII}{\vcenter{\hbox{\includesvg[inkscapelatex=false,height=.5985em]{SII}}}}
\newcommand{\motISI}{\vcenter{\hbox{\includesvg[inkscapelatex=false,height=.5985em]{ISI}}}}
\newcommand{\motSSS}{\vcenter{\hbox{\includesvg[inkscapelatex=false,height=.5985em]{SSS}}}}
\newcommand{\motSIS}{\vcenter{\hbox{\includesvg[inkscapelatex=false,height=.5985em]{SIS}}}}
\newcommand{\motIII}{\vcenter{\hbox{\includesvg[inkscapelatex=false,height=.5985em]{III}}}}

\newcommand{\motSSStr}{\vcenter{\hbox{\includesvg[inkscapelatex=false,height=1.539em]{SSStr}}}}
\newcommand{\motSSItr}{\vcenter{\hbox{\includesvg[inkscapelatex=false,height=1.539em]{SSItr}}}}
\newcommand{\motSIItr}{\vcenter{\hbox{\includesvg[inkscapelatex=false,height=1.539em]{SIItr}}}}
\newcommand{\motIIItr}{\vcenter{\hbox{\includesvg[inkscapelatex=false,height=1.539em]{IIItr}}}}

\newcommand{\motSSSIsqo}{\vcenter{\hbox{\includesvg[inkscapelatex=false,height=1.92375em]{SSSIsqo}}}}
\newcommand{\motSISIsqo}{\vcenter{\hbox{\includesvg[inkscapelatex=false,height=1.92375em]{SISIsqo}}}}
\newcommand{\motSSIIsqo}{\vcenter{\hbox{\includesvg[inkscapelatex=false,height=1.92375em]{SSIIsqo}}}}
\newcommand{\motSIIIsqo}{\vcenter{\hbox{\includesvg[inkscapelatex=false,height=1.92375em]{SIIIsqo}}}}

\newcommand{\motSISSsqi}{\vcenter{\hbox{\includesvg[inkscapelatex=false,height=1.2825em]{SISSsqi}}}}
\newcommand{\motSISIsqi}{\vcenter{\hbox{\includesvg[inkscapelatex=false,height=1.2825em]{SISIsqi}}}}
\newcommand{\motISSSsqi}{\vcenter{\hbox{\includesvg[inkscapelatex=false,height=1.2825em]{ISSSsqi}}}}
\newcommand{\motIISSsqi}{\vcenter{\hbox{\includesvg[inkscapelatex=false,height=1.2825em]{IISSsqi}}}}
\newcommand{\motIISIsqi}{\vcenter{\hbox{\includesvg[inkscapelatex=false,height=1.2825em]{IISIsqi}}}}
\newcommand{\motISISsqi}{\vcenter{\hbox{\includesvg[inkscapelatex=false,height=1.2825em]{ISISsqi}}}}
\newcommand{\motIIISsqi}{\vcenter{\hbox{\includesvg[inkscapelatex=false,height=1.2825em]{IIISsqi}}}}

\newcommand{\motSSSIsqii}{\vcenter{\hbox{\includesvg[inkscapelatex=false,height=1.92375em]{SSSIsqii}}}}
\newcommand{\motSSIIsqii}{\vcenter{\hbox{\includesvg[inkscapelatex=false,height=1.92375em]{SSIIsqii}}}}
\newcommand{\motSIIIsqii}{\vcenter{\hbox{\includesvg[inkscapelatex=false,height=1.92375em]{SIIIsqii}}}}

\newcommand{\motSSSIstr}{\vcenter{\hbox{\includesvg[inkscapelatex=false,height=1.2825em]{SSSIstr}}}}
\newcommand{\motSISIstr}{\vcenter{\hbox{\includesvg[inkscapelatex=false,height=1.2825em]{SISIstr}}}}
\newcommand{\motSSIIstr}{\vcenter{\hbox{\includesvg[inkscapelatex=false,height=1.2825em]{SSIIstr}}}}

\newcommand{\motISSSstr}{\vcenter{\hbox{\includesvg[inkscapelatex=false,height=1.2825em]{ISSSstr}}}}

\newcommand{\motISSIstr}{\vcenter{\hbox{\includesvg[inkscapelatex=false,height=1.2825em]{ISSIstr}}}}
\newcommand{\motIISIstr}{\vcenter{\hbox{\includesvg[inkscapelatex=false,height=1.2825em]{IISIstr}}}}
\newcommand{\motISIIstr}{\vcenter{\hbox{\includesvg[inkscapelatex=false,height=1.2825em]{ISIIstr}}}}

\newcommand{\motSSSI}{\vcenter{\hbox{\includesvg[inkscapelatex=false,height=.55575em]{SSSI}}}}
\newcommand{\motSISI}{\vcenter{\hbox{\includesvg[inkscapelatex=false,height=.55575em]{SISI}}}}
\newcommand{\motSSII}{\vcenter{\hbox{\includesvg[inkscapelatex=false,height=.55575em]{SSII}}}}

\newcommand{\motISSI}{\vcenter{\hbox{\includesvg[inkscapelatex=false,height=.55575em]{ISSI}}}}
\newcommand{\motIISI}{\vcenter{\hbox{\includesvg[inkscapelatex=false,height=.55575em]{IISI}}}}

\newcommand{\motSSSItre}{\vcenter{\hbox{\includesvg[inkscapelatex=false,height=1.2825em]{SSSItre}}}}

\newcommand{\motSSIItre}{\vcenter{\hbox{\includesvg[inkscapelatex=false,height=1.2825em]{SSIItre}}}}

\newcommand{\motISIItre}{\vcenter{\hbox{\includesvg[inkscapelatex=false,height=1.2825em]{ISIItre}}}}

\renewcommand{\thesection}{\Roman{section}}
\renewcommand{\thesubsection}{\Alph{subsection}}
\renewcommand{\theparagraph}{\alph{paragraph}}

\makeatletter
\def\p@paragraph{\thesection\,\thesubsection\,}
\makeatother

\begin{document}

\title{Mean-field models of dynamics on networks via moment closure: \\an automated procedure}

\author{Bert Wuyts}
\affiliation{College of Engineering, Mathematics and Physical Sciences, 
University of Exeter, EX4 4QF, UK}
\author{Jan Sieber}
\affiliation{College of Engineering, Mathematics and Physical Sciences, 
University of Exeter, EX4 4QF, UK}
\email{b.wuyts@ex.ac.uk}


\begin{abstract}
In the study of dynamics on networks, moment closure is a commonly used method to obtain low-dimensional evolution equations amenable to analysis. The variables in the evolution 
equations are mean counts of subgraph states and are referred to as moments. Due to interaction between neighbours, each moment equation is a function of higher-order moments, such that an infinite hierarchy of equations arises. Hence, the derivation requires truncation at a given
order, and, an approximation of the highest-order moments in terms of lower-order ones,
known as a closure formula. Recent systematic approximations 
have either restricted focus to closed moment equations for SIR epidemic spreading or to unclosed moment equations for arbitrary dynamics.
In this paper, we develop a general procedure that automates both derivation and closure of arbitrary order moment equations for dynamics with 
nearest-neighbour interactions on undirected networks. 
Automation of the closure step was made possible by our generalised closure scheme,
which systematically decomposes the largest subgraphs into their smaller components.
We show that this decomposition is exact
if these components form a tree, there is independence at distances
beyond their graph diameter, and there is spatial homogeneity. Testing our method for SIS epidemic spreading on lattices and random networks confirms that biases are larger
for networks with many short cycles in regimes with long-range dependence. A \texttt{Mathematica} package that automates the moment closure is available for download. 

\end{abstract}

\maketitle


\section{Introduction} 
The dynamics of complex systems are usually most accurately represented by high-dimensional
stochastic simulation models. However, their large
state space makes exact mathematical analysis prohibitive. Therefore, one often looks for
low-dimensional approximations that permit analysis. 
In \emph{moment closure},
one achieves this by studying the time evolution of a finite set of ``moments''
rather than that of the full probability distribution of the considered 
stochastic dynamical system \cite{kuehn2016}. The complete set of moment
equations forms an infinite hierarchy of ordinary differential equations (ODEs), with lower-order moments depending on higher order ones.
The approximation
consists then of truncating the hierarchy at a chosen order and replacing the highest order
moments by functions of the lower-order moments. Such functions are known as
closure formulas, and they can be obtained in various ways \cite{kuehn2016}, such as via an assumption of
statistical independence \cite[e.g.][]{Sharkey2015,Sharkey2015a,Sharkey2011}, physical principles (e.g. maximum entropy \cite{Rogers2011}),
time scale separation \cite{GrossKevrekidis2008}, and assumptions on the type of probability distribution \cite[e.g.][]{Isham1991}.
In this work,
we only consider closures derived from an assumption of statistical independence,
which is equivalent to a \emph{mean-field approximation}, a method originating
from the statistical physics of phase transitions in materials \cite{weiss1907,bragg1934,bethe1935,kikuchi1951}.
First-order moment closure assumes pairwise independence of species counts and corresponds to the `mean field'
or `simple mean field' \cite[e.g][]{Marro1999,Henkel2008,Tome2015}, resulting in 
equations only for total counts of each species. Likewise, second order moment closure assumes independence of pair counts in larger units, and corresponds to the `pair approximation' \cite[e.g.][]{Matsuda1992,Keeling1997,Rand1999,Dieckmann2000,Kefi2007a,Gross2006},
which also includes equations for pair counts. 
We aim to exploit this connection between moment closure and mean-field approximations
to generalise and automate the derivation of arbitrary-order mean-field 
models for arbitrary dynamics with 
nearest-neighbour interactions on undirected networks. In the rest of our introduction, we
introduce our approach with more precision,
discuss the relevant literature, state our aims, and provide an overview
of the paper contents.

In general, moments in the moment closure for dynamics on networks represent the expected frequencies of small subgraph states known as \emph{network motifs} \cite{House2009}.
Derivation of the moment equations proceeds from smaller to larger
sized-motifs, with dynamics of mean motif counts of size $m$ depending
only on mean motif counts of size $m$ and $m{+}1$ if the dynamics has only nearest-neighbour interactions. Hence,
a system of ODEs obtained in such a manner for motif counts 
up to a maximum considered size $k$ (also referred to as the order 
of the moment closure) is always underdetermined,
because it depends on motifs of size $k{+}1$ but does not contain equations for them.
Therefore, as the second step of moment closure, a closure 
approximation is applied by 
expressing counts of $(k{+}1)$-size motifs as functions
of counts of $\{1,...,k\}$-size motifs, closing the system of ODEs. In this substitution,
larger-sized motifs factorise in terms of smaller-sized
ones, which we will justify, as mentioned above, by an assumption of statistical independence.  
For homogeneous networks, closures that are valid at the individual level
(i.e. concerning states of given nodes) are also valid at the population level 
(i.e. concerning total counts or averages of states in the whole network)
\cite{Sharkey2008},
which then permits a compact description in terms of population-level 
quantities. Yet, the number of motif types, and hence equations, increases
combinatorially with $k$. It is then hoped that
the derivation can be stopped at an order low enough for the resulting
system of ODEs to be sufficiently amenable to analytical or numerical
methods and high enough to satisfy the independence assumptions underlying the closure approximately. We note that various other types of approximations exist that
focus on specific types of subgraphs, such as active motifs \cite{Bohme2011}, star graphs \cite{gleeson2013,fennell2019}
or hypergraphs of cliques \cite{marceau2010,stonge2021} or of general motifs \cite{Cui2022}. These have also been referred to as moment closure \cite[e.g. in][]{demirel2014}, but their truncation order and closure formulas are implicit in the method. At present, the equations obtained by the approach of references \cite{gleeson2013,fennell2019,marceau2010,stonge2021,Cui2022} are referred to as approximate master equations.

In context of population dynamics on networks or lattices, moment closure methods
have been used to study applications such as spatial ecology \cite[e.g][]{Matsuda1992,Dieckmann2000,Kefi2007a}, 
epidemics \cite[e.g.][]{Keeling1997,Keeling1999,Rand1999,Kiss2017,House2009,Gross2006},
opinion formation \cite[e.g.][]{demirel2014}, evolution of cooperation \cite[e.g.][]{Szabo2002}, among others. 
Despite the wide use of moment closure in applications, derivation of the moment equations is customarily done separately for each considered process and approximation order, while justifying the used closure formulas only heuristically. Comparatively few studies have shown how moment equations derive in general from the master equation or how
closure formulas arise from precisely defined independence assumptions.
Regarding the former, an automated algorithm to derive (unclosed) moment equations for arbitrary adaptive dynamics on directed networks was recently developed by Danos et al. \cite{Danos2020}. Regarding the latter, 
attention has centred on the specific case of SIR epidemic spreading  
\cite{Sharkey2015,Sharkey2015a,Sharkey2011}
because proving validity of low-order closures is least challenging here. 
In particular, Sharkey et al.
\cite{Sharkey2015a} proved that an exact 
individual-level closure approximation exists for motifs 
that have an all-susceptible set of nodes which cuts all possible chains of infection
between the remaining parts when removed. In tree networks, this is already possible
with three nodes 
such that the largest required motif in the moment equations is of size 2. In case of non-tree networks, larger
motifs need to be taken into account, resulting in a larger number of equations.
Hence, while Danos et al. \cite{Danos2020} have shown that it is feasible to derive moment equations in a generic form, the work on SIR spreading 
\cite{Sharkey2011,Sharkey2015a,Sharkey2015} 
indicates which type of independence assumptions are required to obtain valid closures.

In this paper, we 
provide a first \emph{fully automated} procedure for both derivation and closure of population-level moment equations up to any order and
for arbitrary dynamics with at most nearest-neighbour interactions on undirected networks. 
Automated derivation is made possible by our generic moment equation \eqref{eq:dinvfull1},
which we derived from the master equation. As mentioned above, more general derivations
than ours exist \cite{Danos2020}. Hence, our main contribution is in automating also
the closure, 
which we show and justify in detail (Section \ref{sec:Closure-scheme}). 
Our closure scheme (equations \eqref{eq:dec:mot:chord} or 
\eqref{eq:dec:mot:nchord}) generalises
previous insights \cite{Sharkey2011,Sharkey2015a,Sharkey2015} and relies on the
theory of Markov networks \cite{Pearl1988} to make it applicable to motifs of any type and size, such that it can close any set of moment equations at arbitrary order.
We show that our closure scheme is exact if the motifs that are decomposed by the closure form a tree and if there is independence beyond their graph diameter. 

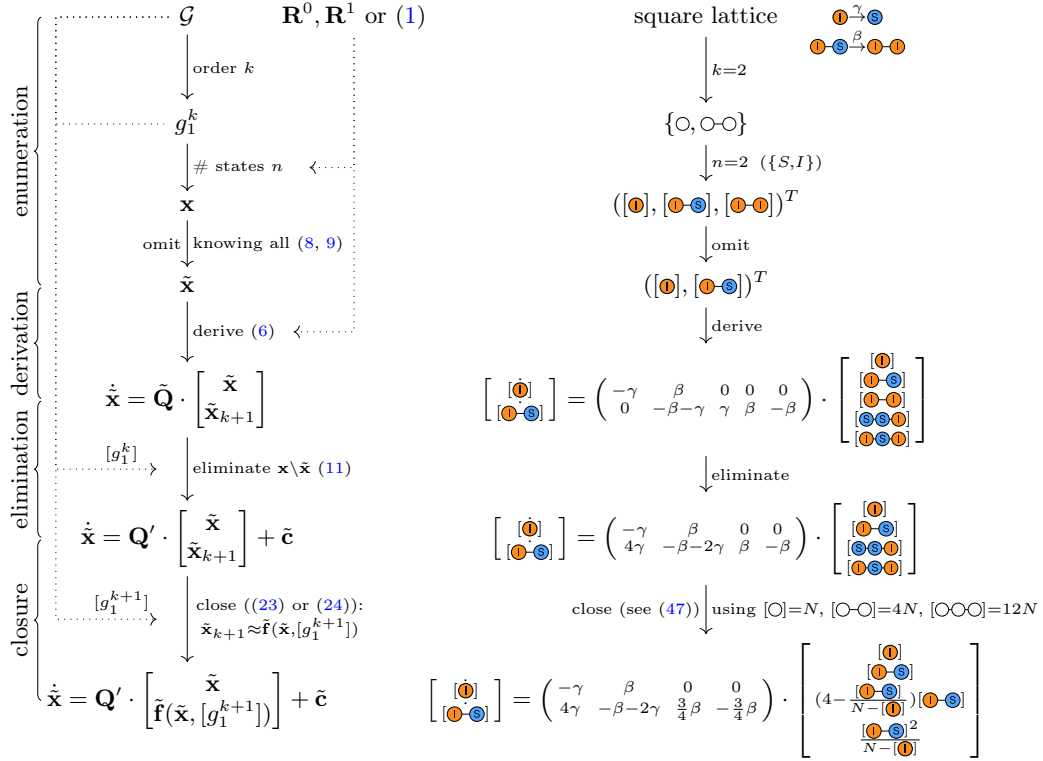
\begin{figure*}
	\begin{tikzcd}[column sep=4pt,row sep=14pt]
		{}&{\cal G}\arrow[d,"\text{order }k"] &[-29pt]\mathbf{R}^0,\mathbf{R}^1\text{ or \eqref{eq:reactionrules}}&[-12] \text{square lattice}\arrow{d}{k=2}&[-80pt]	\begin{smallmatrix}\\[7.3pt]\motI\overset{\gamma}{\to}\motS\\\motIS\overset{\beta}{\rightarrow}\motII\end{smallmatrix}
		\\ 
		{}&g_1^k\arrow[d,"\text{\# states }n",""{name=S}]&&\{\motaNL,\motabNL\}\arrow[d,"n=2\;\; {(\{S,I\})}",""{name=S2}]&
		\\
		{}&\mathbf{x}\arrow[d,"\text{omit}"',"{\text{knowing all (\ref{eq:cons:norm},~\ref{eq:cons2})}}"]&&([\motI],[\motIS],[\motII])^T\arrow{d}{\text{omit}}&
		\\	 
		{}&\tilde{\mathbf{x}}\arrow[d,"\text{derive  }\eqref{eq:dinvfull1}",""{name=D}]&&([\motI],[\motIS])^T\arrow[d,"\text{derive }",""{name=D2}]&
		\\
		{}&\dot{\tilde{\mathbf{x}}}=\tilde{\mathbf{Q}}
		\cdot
		\begin{bmatrix}
			\tilde{\mathbf{x}}\\
			\tilde{\mathbf{x}}_{k+1}\\
		\end{bmatrix}
		\arrow[d,""{name=E},"\text{eliminate }\mathbf{x}\setminus\tilde{\mathbf{x}}\text{ }\eqref{eq:subcons}"]
		&&
		\left[\begin{smallmatrix}
			\dot{[\motI]}\\\dot{[\motIS]}\\
		\end{smallmatrix}\right]
		=
		\left(\begin{smallmatrix}
			-\gamma & \beta & 0 & 0 & 0\\
			0 &-\beta-\gamma & \gamma & \beta & -\beta\\
		\end{smallmatrix}\right)
		\cdot
		\left[\begin{smallmatrix}
			[\motI]\\ [\motIS]\\ [\motII]\\ [\motSSI]\\ [\motISI]\\
		\end{smallmatrix}\right]\arrow{d}{\text{eliminate}}&
		\\
		{}&\dot{\tilde{\mathbf{x}}}=\mathbf{Q}'
	 	\cdot
	 	\begin{bmatrix}
	 		\tilde{\mathbf{x}}\\
	 		\tilde{\mathbf{x}}_{k+1}\\
	 	\end{bmatrix}
	 	+\tilde{\mathbf{c}}
		\arrow[d,"\begin{smallmatrix}\text{close (\eqref{eq:dec:mot:chord} or \eqref{eq:dec:mot:nchord}):}\\{\tilde{\mathbf{x}}_{k+1}\approx\tilde{\mathbf{f}}(\tilde{\mathbf{x}},[g_1^{k+1}])}\end{smallmatrix}",""{name=C}]&&
	 	\left[\begin{smallmatrix}
			\dot{[\motI]}\\ \dot{[\motIS]}\\
		\end{smallmatrix}\right]
		=
		\left(\begin{smallmatrix}
			-\gamma & \beta & 0 & 0\\
			4\gamma &-\beta-2\gamma & \beta & -\beta\\
		\end{smallmatrix}\right)
		\cdot
		\left[\begin{smallmatrix}
			[\motI]\\ [\motIS]\\ [\motSSI]\\ [\motISI]\\
		\end{smallmatrix}\right]\arrow[d,"\text{close (see \eqref{eq:mf2:hom:closure})}"',"\text{using }{[\motaNL]=N,\;[\motabNL]=4N,\;[\motabcNL]=12N}"]&
		\\[5pt]
		{}&\dot{\tilde{\mathbf{x}}}=\mathbf{Q}'
		\cdot
		\begin{bmatrix}
			\tilde{\mathbf{x}}\\
			\tilde{\mathbf{f}}(\tilde{\mathbf{x}},[g_1^{k+1}])\\
		\end{bmatrix}
		+\tilde{\mathbf{c}}
		&&
		\left[\begin{smallmatrix}
			\dot{[\motI]}\\ \dot{[\motIS]}\\
		\end{smallmatrix}\right]
		=
		\left(\begin{smallmatrix}
			-\gamma & \beta & 0 & 0\\
			4\gamma &-\beta-2\gamma & \frac{3}{4}\beta & -\frac{3}{4}\beta\\
		\end{smallmatrix}\right)
		\cdot
		\left[\begin{smallmatrix}
			[\motI]\\ [\motIS]\\ (4-\frac{[\motIS]}{N-[\motI]})[\motIS]\\ \frac{[\motIS]^2}{N-[\motI]}\\
		\end{smallmatrix}\right]
	\arrow[dotted,from=1-3,to=D,to path={ -- (\tikztostart.south) |- (\tikztotarget)},shorten >= 36pt]
	\arrow[dotted,from=1-3,to=S,to path={ -- (\tikztostart.south) |- (\tikztotarget)},shorten >= 42pt]
	\arrow["{[g_1^{k+1}]}",near end,dotted,from=1-2,to=C,to path={ -- ([xshift=-11ex]\tikztostart.west) |- (\tikztotarget)\tikztonodes},shorten >= 3.ex]
	\arrow["{[g_1^{k}]}",near end,dotted,from=1-2,to=E,to path={ -- ([xshift=-11ex]\tikztostart.west) |- (\tikztotarget)\tikztonodes},shorten >= 3.ex]
	\arrow[dotted,dash,from=1-2,to=2-2,to path={ -- ([xshift=-11ex]\tikztostart.west) |- (\tikztotarget)}]
	\arrow[from=1-1,to=4-1, start anchor=west, end anchor={[yshift=.5pt]west}, no head, decorate, decoration={brace,mirror}, "\rotatebox{90}{\text{enumeration}}" left=3pt,xshift=3.5ex]
	\arrow[from=4-1,to=5-1, start anchor={[yshift=-.5pt]west}, end anchor={[yshift=.5pt]west}, no head, decorate, decoration={brace,mirror}, "\rotatebox{90}{\text{derivation}}" left=3pt,xshift=3.5ex]
	\arrow[from=5-1,to=6-1, start anchor={[yshift=-.5pt]west}, end anchor={[yshift=.5pt]west}, no head, decorate, decoration={brace,mirror}, "\rotatebox{90}{\text{elimination}}" left=3pt,xshift=3.5ex]
	\arrow[from=6-1,to=7-1, start anchor={[yshift=-.5pt]west}, end anchor=west, no head, decorate, decoration={brace,mirror}, "\rotatebox{90}{\text{closure}}" left=3pt,xshift=3.5ex]
	\end{tikzcd}
	\caption{Mean-field models for multistate dynamics on networks via moment closure. Left: general form. Right: example (MF2 for SIS spreading). Input (top): network ${\cal G}$, the considered process defined by the reaction rate matrices $\mathbf{R}^0,\mathbf{R}^1$. Output (bottom): closed system of differential equations representing mean frequencies
	of motif counts in the network up to a maximum size $k$. The four main steps are: (\emph{i}) enumerate all required motifs, (\emph{ii}) derive the moment equations, (\emph{iii}) substitute to be eliminated terms via conservation relations, (\emph{iv}) close the system via a closure scheme.} \label{fig:overview}
\end{figure*}

As shown in Figure \ref{fig:overview}, the whole procedure consists of four main steps: (\emph{i}) enumeration of all required motifs
up to size $k$, 
(\emph{ii}) derivation of
the unclosed ODEs for these motifs, 
(\emph{iii}) elimination via conservation relations, 
and (\emph{iv}) closure of the system
of ODEs. We will refer to the final closed system of ODEs as the $k$-th order 
mean field, or MF$k$ in short. Elimination is not strictly necessary but
it increases efficiency, particularly for higher-order approximations.
We developed a \texttt{Mathematica} \cite{Mathematica} package that
derives MF$k$ by performing steps (\emph{i})--(\emph{iv}). 
The required inputs for this
algorithm are: the counts of all induced subgraphs
in the underlying network up to size $k$, and, the matrices
$\mathbf{R}^0$, $\mathbf{R}^1$ with conversion and interaction rates. 
A code example with output is discussed in 
Appendix \ref{sec:MFhoeqs} and the package is available for download \cite{Wuyts2022mfmcl}.

In Sections~\ref{sec:ds-ct-mc} and \ref{sec:motifs}, the underlying Markov chain for the
network and its motifs are introduced. Section~\ref{sec:diff:inv} explains the general
formula \eqref{eq:dinvfull1} for step (\emph{ii}), which follows from
the master equation (for its derivation, see Appendix
\ref{sec:Differential-and}). 
The conservation relations used for elimination of variables from the ODEs in
step (\emph{iii}) are
shown in equations (\ref{eq:cons:norm},~\ref{eq:cons2}) of Section
\ref{sec:alg:iv}. The form of the system of moment equations up to
a truncation order $k$ \eqref{eq:xdotblock} and
the variable elimination are shown in Section \ref{sec:QE}, resulting in the
unclosed system \eqref{eq:xdotsubcons2}. Finally, general expressions to close 
the system of moment equations in step (\emph{iv}) are derived in
Section~\ref{sec:Closure-scheme} (equations
\eqref{eq:dec:mot:chord} or \eqref{eq:dec:mot:nchord}), resulting in form
\eqref{eq:xdotsubconscl}.
In
Section \ref{sec:SIS} we set up MF1-5 models of SIS epidemic spreading
and compare their steady states to those of simulations 
on a selection of networks.  We focus in particular
on the square lattice, for which low-order moment closures 
fail, due to its large number of
cycles of any size, and we compare against random networks and
higher-dimensional lattices, for which they work well.

\section{The underlying discrete-state continuous-time Markov chain}
\label{sec:ds-ct-mc}
We consider a dynamical system on a fixed undirected graph $\cal{G}$
with $N$ nodes and adjacency matrix
$\mathsf{A}\in\{0,1\}^{N\times N}$, where each node may have one out
of $n$ discrete states. We may see the nodes as locations and the
states as species, such that the space at node
$i\in\{1,...,N\}$ is occupied by exactly one species in
$\{1,...,n\}$. We denote the state vector at time $t$ by $X(t)$ such that $[X_i(t)]_{i=1}^N\in\{1,...,n\}^N$. 
The type of dynamics we
consider is a continuous-time Markov chain with two Poisson process transition types, with rates specified by a $n\times n$ matrix $\mathbf{R}^0$ and a
$n\times n\times n$ tensor $\mathbf{R}^1$:
\begin{itemize}[noitemsep,nosep]
\item[(\emph{i})] $\mathbf{R}^0$
specifies \emph{spontaneous conversion} rates. Any node with state $a$ may
change spontaneously into state $b$, with rate $R^0_{ab}$;
\item[(\emph{ii})] $\mathbf{R}^1$ specifies \emph{nearest-neighbour-induced conversion} rates. Any node with state $a$ may
change into state $b$ for each link to a node with state $c$, with
rate $R^1_{abc}$.
\end{itemize}
This corresponds to the reaction rules
\begin{equation}\label{eq:reactionrules}
\mota\overset{R^0_{ab}}{\to}\motb,\qquad\motab\overset{R^1_{abc}}{\rightarrow}\motac.
\end{equation}

A simple example is \emph{susceptible-infected-susceptible} (SIS) \emph{epidemic spreading} on a square two-dimensional lattice
of $\sqrt{N}\times\sqrt{N}$ nodes (say, periodic in both directions).
For SIS spreading, each node may have one of two states, \emph{susceptible} or
\emph{infected}, $\{\motS,\motI\}=\{1,2\}$, for $t\in[0,\infty)$ and
$i\in\{1,...,N\}$. $S$ nodes can become infected at rate $\beta$
per infected neighbour and $I$ nodes recover spontaneously at rate $\gamma$. 
Hence, for SIS spreading, the matrix $\mathbf{R}^0$ has a single non-zero
entry $R^0_{2,1}=\gamma>0$ (for spontaneous \emph{recovery}) and the
tensor $\mathbf{R}^1$ has a single non-zero entry
$R^1_{1,2,2}=\beta>0$ (for \emph{infection} along $IS$ links, denoted
by the symbol $\motIS$ below), corresponding to the reaction rules 
\begin{equation}
	\motI\overset{\gamma}{\to}\motS,\qquad\motIS\overset{\beta}{\rightarrow}\motII.\label{eq:reactionSIS}
\end{equation}

\section{Network motifs and their counts}
\label{sec:motifs}
We define network motifs as (typically small) graphs
with given state labels.
The order of a motif is the number
of nodes it has. E.g., in our example, 
$I$ nodes and $IS$ links are examples of first and second-order motifs.
We use square brackets to denote the
count of occurrences of motifs in $({\cal G},X)$, i.e. the number of occurrences
of $I$ nodes is $[\motI]$. Hence, we can write e.g. 
respectively for $I$ nodes, $IS$ links and $ISI$ chains
\begin{align*}
[\motI]=&\sum_{\mathclap{i\in\{1,...,N\}}}\delta_2(X_i)\mbox{,}
[\motIS]=\sum_{\mathclap{i,j\neq i\in\{1,...,N\}}}\delta_1(\mathsf{A}_{ij})\delta_{2}(X_{i})\delta_{1}(X_{j})\mbox{,}\\[2pt]
[\motISI]=&\sum_{\mathclap{i,j\neq i,k\neq j\in\{1,...,N\}}}\delta_1(\mathsf{A}_{ij})\delta_1(\mathsf{A}_{jk})\delta_0(\mathsf{A}_{ik})\delta_{2}(X_{i})\delta_{1}(X_{j})\delta_{2}(X_{k})\mbox{,}
\end{align*}
where we omitted the dependence on $t$. The Kronecker delta,
$\delta_y(x)$, is $1$ if $x$ equals $y$ and $0$ otherwise. 
By construction,
the motif counts on the left-hand side are random, since $X(t)$ is random. 

Generalizing the above examples, a \emph{network motif} of \emph{order} $m$, 
is a network with $m$
nodes, each of which 
are labelled with a state. It his hence fully characterised by its connectivity pattern between nodes and its state labels on the nodes. The connectivity between motif nodes,
i.e. the motif without labels, will be 
indicated by $\mathsf{a}$, which,
depending on the context, denotes the adjacency matrix of the
motif, or a set
$\mathsf{a}\in(\{1,...,m\}\times\{1,...,m\})^\mu$ of links
between the $m$ nodes (the indices of the non-zero entries in the
adjacency matrix of the motif such that $\mu\leq m(m-1)/2$), or a graphical
representation of the connectivity. For instance, two linked nodes
are displayed according to these representations as ${\big(\begin{smallmatrix}0 & 1\\ 1 & 0 \end{smallmatrix}\big)}$, $\{(1,2)\}$, or $\motabNL$ respectively. As we focus only on undirected networks, each pair in the pair representation is bidirectional, i.e. we write
$\{(1,2)\}$ instead of $\{(1,2),(2,1)\}$. Motif labels will be
indicated via a vector
$\boldsymbol{x}=(x_{1},...,x_{m})\in\{1,...,n\}^m$ with state
labels $x_{p}$ at positions $p=1,...,m$. Hence, the pair of
$\boldsymbol{x}$ and $\mathsf{a}$ describes the motif, which
we write as $\boldsymbol{x}^{\mathsf{a}}$. For example
$(2,1)^{(1,2)}=(\motI,\motS)^{(1,2)}$ is the $\motIS$ link.

For a general network motif $\boldsymbol{x}^\mathsf{a}$, the
total count $[\boldsymbol{x}^{\mathsf{a}}]=[\boldsymbol{x}^{\mathsf{a}}](t)$ in the 
network ${\cal G}$ with fixed adjacency $\mathsf{A}$ and node labels $X=X(t)$ is 
\begin{align}\label{eq:motifcount:all}
\left[\boldsymbol{x}^{\mathsf{a}}\right]=\sum_{\boldsymbol{i}\in S(m,N)}\delta_{\mathsf{a}}(\mathsf{A}_{\boldsymbol{i}})\delta_{\boldsymbol{x}}(X_{\boldsymbol{i}})\mbox{,}
\end{align}
where $S(m,N)$ is the set of all $m$-tuples from $\{1,...,N\}$
without repetition, which has size $|S(m,N)|=N!/(N-m)!$. For instance,
$\sum_{\boldsymbol{i}\in
  S(3,N)}=\sum_{i\in\{1,...,N\}}\sum_{j\in\{1,...,N\}\setminus
  i}\sum_{k\in\{1,...,N\}\setminus \{i,j\}}$. For an index set
  $\boldsymbol{i}\subset\{1,...,N\}^m$ of length $m$ we use the convention that
  $\mathsf{A}_{\boldsymbol{i}}\in\{0,1\}^{m\times m}$ and
  $X_{\boldsymbol{i}}\in\{1,...,n\}^m$ are the restrictions of
matrix $\mathsf{A}$ and vector $X$ to the index set
$\boldsymbol{i}$. We count exact matches between the motif and the subgraph in $({\cal G},X)$.
This means that the counted motif needs to have $\mathsf{a}$ as induced subgraph of ${\cal G}$,
with matching state labels. Counting via \eqref{eq:motifcount:all} leads to
multiple counting of motifs with symmetries (more precisely: automorphisms -- see Appendix \ref{sec:Differential-and}), with multiplicity equal
to the number of symmetries. For instance, $\motIS$ is counted once;
but $\motSS$ or $\motII$ twice, and $\motSSStr$ six times.

\section{Differential equations for motif counts}
\label{sec:diff:inv} In our SIS spreading example, the
expected rate of change for the count of infected nodes is
well known to satisfy
\begin{align}
\ddt\langle[\motI]\rangle & =  -\gamma\langle[\motI]\rangle+\beta\langle[\motIS]\rangle\mbox{,}\label{eq:ov:mf1}
\end{align}
where $\langle\cdot\rangle$ brackets denote expectations over
many independent realisations of the underlying Markov
chain. Relation \eqref{eq:ov:mf1} is exact for finite network sizes
$N$ and can be derived from the Kolmogorov-forward (or master)
equation for the Markov chain, see Appendix
\ref{sec:Differential-and}. Note that the structure in \eqref{eq:ov:mf1}
is such that the expected rate for the frequency of a motif of size
$m=1$ depends on the counts of motifs of size $m=1$
(here $\langle[\motI]\rangle$) from spontaneous conversions (here recovery) and
size $m=2$ (here $\langle[\motIS]\rangle$) from nearest-neighbour-induced conversions (here infection). 

This is true in general such that the count of a general motif
$[\boldsymbol{x}^\mathsf{a}](t)$ of size $m$ satisfies an ordinary differential equation of form
\begin{equation}\label{eq:dinv}
\ddt\langle\left[\boldsymbol{x}^{\mathsf{a}}\right]\rangle=F_{\boldsymbol{x}^\mathsf{a}}\left(\langle[\boldsymbol{y}_1^{\mathsf{b}_1}]\rangle,\langle[\boldsymbol{y}_2^{\mathsf{b}_2}]\rangle,...\right)\mbox{,}    
\end{equation}
where $|\boldsymbol{y}_i|\in\{m,m{+}1\}$. On the right-hand
side 
the $[\boldsymbol{y}_i^{\mathsf{b}_i}](t)$ stand for counts of
motifs of size $m$ or $m{+}1$ on which the dynamics of
$\bigl\langle\left[\boldsymbol{x}^{\mathsf{a}}\right]\bigr\rangle$
depend. Our package expresses the right-hand side of the differential equation
\eqref{eq:dinv} for arbitrary motifs $\boldsymbol{x}^\mathsf{a}$ of size $m$ in the general form 
\begin{widetext}
	\begin{equation}\label{eq:dinvfull1}
\begin{aligned}
\ddt\langle\left[\boldsymbol{x}^{\mathsf{a}}\right]\rangle =&
\sum_{p=1}^m\sum_{k=1}^n\sum_{c=1}^{n}\Biggl\{\left(\frac{R^0_{kx_{p}}}{n}+\kappa_{p,c}^{\boldsymbol{x}^{\mathsf{a}}}R^1_{kx_{p}c}\right)\bigl\langle\left[\boldsymbol{x}^{\mathsf{a}}_{p\shortto k}\right]\bigr\rangle-
\left(\frac{R^0_{x_{p}k}}{n}+\kappa_{p,c}^{\boldsymbol{x}^{\mathsf{a}}}R^1_{x_pkc}\right)\langle\left[\boldsymbol{x}^{\mathsf{a}}\right]\rangle
+\\
&\sum_{\boldsymbol{y}^\mathsf{b}\in{\cal N}^c_p\left(\boldsymbol{x}^{\mathsf{a}}_{p\shortto k}\right)}
R^1_{kx_{p}c}\langle[\boldsymbol{y}^\mathsf{b}]\rangle-
\sum_{\boldsymbol{y}^\mathsf{b}\in{\cal N}^c_p(\boldsymbol{x}^\mathsf{a})}
R^1_{x_{p}kc}\langle[\boldsymbol{y}^\mathsf{b}]\rangle\Biggr\}.
\end{aligned}
\end{equation}
\end{widetext}
Here, 
$\boldsymbol{x}_{p\shortto k}$ is the state label vector obtained by
setting the state label of the $p$th element of $\boldsymbol{x}$ to
$k$. $\kappa_{p,c}^{\boldsymbol{x}^{\mathsf{a}}}$ is the $c$-degree at
position $p$ in motif $\boldsymbol{x}^{\mathsf{a}}$, i.e. it is the
number of connections node $p$ in motif $\boldsymbol{x}^{\mathsf{a}}$
has to nodes with state label $c$. The set
${\cal N}^c_p(\boldsymbol{x}^{\mathsf{a}})$ is defined as
\begin{equation*}\label{overview:nbhdef}
    {\cal N}_{p}^c(\boldsymbol{x}^\mathsf{a}){:=}\bigcup_{\ell=1}^{m+1}
    \left\{\boldsymbol{y}^\mathsf{b}{:}|\boldsymbol{y}|{=}m{+}1, y_\ell{=}c, \boldsymbol{y}^\mathsf{b}_{\ell\shortto\emptyset}{=}\boldsymbol{x}^\mathsf{a},(\ell,p){\in}\mathsf{b} \right\},
\end{equation*}
and contains all 
motifs $\boldsymbol{y}^\mathsf{b}$ of order $m{+}1$
that extend
the state label vector $\boldsymbol{x}$ by one new node with state
label $c$ and extend adjacency $\mathsf{a}$ by links to the new node, 
where $\boldsymbol{y}^\mathsf{b}_{\ell\shortto\emptyset}$ denotes the $m$th order connected motif obtained by deleting the $\ell$th node of $\boldsymbol{y}^\mathsf{b}$ and its links.
The differential equation \eqref{eq:dinvfull1} shows how
the expected count of $\boldsymbol{x}^{\mathsf{a}}$ is increased by
transitions \emph{into} $\boldsymbol{x}^{\mathsf{a}}$ -- the positive
terms -- and decreased by transitions \emph{out of}
$\boldsymbol{x}^{\mathsf{a}}$ -- the negative terms. This happens
through spontaneous conversions (terms with $R^0$), through
nearest-neighbour interaction between nodes within the motif (first
two terms with $R^1$), or through nearest-neighbour interaction with
nodes outside the motif (last two terms with $R^1$). Note that 
equivalent motifs up to permutation (isomorphic motifs -- 
see Appendix \ref{sec:Differential-and}) result in the same equation, such that
we choose the same representative node indexing for each equivalence class.

\section{Conservation relations}
\label{sec:alg:iv}
Conservation relations are linear algebraic relations between motif counts. E.g., 
for SIS spreading on a square lattice with periodic boundary conditions,
we have for first and second-order motifs:
	\begin{align*}
		[\motS](t)+[\motI](t)&=N\mbox{,}&
		[\motSS](t)+[\motIS](t)+[\motII](t)&=4N\mbox{.}
	\end{align*}
Such conservation of node and link counts occurs because our graph $\cal{G}$ is fixed. 
We can use the total counts on the right-hand side as normalising factors such that we
may write 
\begin{align*}
	\llbracket\motS\rrbracket(t)+\llbracket\motI\rrbracket(t)&=1\mbox{,}& \llbracket\motSS\rrbracket(t)+\llbracket\motIS\rrbracket(t)+\llbracket\motII\rrbracket(t)&=1\mbox{.}
\end{align*}
where $\llbracket\cdot\rrbracket$ is our notation for normalised motif counts. For networks
with homogeneous degree, one can also write conservation equations of the type
$\llbracket\motS\rrbracket=\llbracket\motIS\rrbracket+\llbracket\motSS\rrbracket$.

In general, for
each adjacency matrix $\mathsf{a}$ of possible motifs of size $m$ the conservation
relation 
\begin{equation}\label{eq:cons}
	\sum_{\boldsymbol{x}: |\boldsymbol{x}|=m}\left[\boldsymbol{x}^{\mathsf{a}}\right](t)=\left[\mathsf{a}\right]
	\mbox{,}
	\;\;\mbox{where}\;\; [\mathsf{a}]=\sum_{\boldsymbol{i}\in S(m,N)}\delta_{\mathsf{a}}(\mathsf{A}_{\boldsymbol{i}})\mbox{,}
\end{equation}
holds. The overall count $[\mathsf{a}]$ of induced subgraphs
$\mathsf{a}$ in graph $\cal{G}$ is constant in time, and split up between
all possible labellings.
We can therefore use 
$\left[\mathsf{a}\right]$ as a normalisation factor for the (variable) counts
of motifs such that we may consider normalised motif counts
\begin{equation}\label{eq:cons:norm}
	\llbracket\boldsymbol{x}^{\mathsf{a}}\rrbracket:=\left[\boldsymbol{x}^{\mathsf{a}}\right]/\left[\mathsf{a}\right]\mbox{, such that\ }\sum_{\boldsymbol{x}: |\boldsymbol{x}|=m}\llbracket\boldsymbol{x}^{\mathsf{a}}\rrbracket(t)=1
\end{equation}
is the conservation relation for the normalized quantities.
For networks with homogeneous degree there is the additional type of conservation relation
\begin{equation}\label{eq:cons2}
	{\sum_{k=1}^n}\llbracket\boldsymbol{x}^{\mathsf{a}}_{p\shortto k}\rrbracket=\llbracket\boldsymbol{x}^\mathsf{a}_{p\shortto\emptyset}\rrbracket,
\end{equation}
for each stub $p$ of the motif $\boldsymbol{x}^\mathsf{a}$ (a stub is
a node with degree $1$ in $\mathsf{a}$).
The conservation relations \eqref{eq:cons:norm} and \eqref{eq:cons2} can be used to reduce the
number of variables in the moment equations via substitution. This can result
in a substantial reduction in the number of moment equations (see Table \ref{tab:Number-of-equations}). 

\section{Truncation and substitution}\label{sec:QE}

When we use
\eqref{eq:dinvfull1} to express expected rates of change for the set
of motifs up to a chosen maximum size $k$, 
we obtain a truncated hierarchy of moment equations. This linear system of differential equations has the form
\begin{widetext}
	\begin{equation}\label{eq:xdotblock}
	\begin{bmatrix}
		\dot{\mathbf{x}}_1\\
		\dot{\mathbf{x}}_2\\
		\vdots\\
		\dot{\mathbf{x}}_k\\
	\end{bmatrix}
	=
	\begin{bmatrix}
		\mathbf{Q}_1 & \mathbf{Q}_{12} & \mathbf{0} & \cdots & \cdots & \mathbf{0}\\
		\mathbf{0} & \mathbf{Q}_2 & \mathbf{Q}_{23} & \mathbf{0} & \cdots & \vdots\\
		\vdots & \mathbf{0} & \ddots & \ddots & \mathbf{0} & \vdots\\
		\vdots & \cdots & \mathbf{0} & \mathbf{Q}_{k-1} & \mathbf{Q}_{k-1,k} & \mathbf{0}\\
		\mathbf{0} & \cdots & \cdots & \mathbf{0} & \mathbf{Q}_k & \mathbf{Q}_{k,k+1}
	\end{bmatrix}
	.
	\begin{bmatrix}
		\mathbf{x}_1\\
		\mathbf{x}_2\\
		\vdots\\
		\mathbf{x}_k\\
		\mathbf{x}_{k+1}\\
	\end{bmatrix}
\;=:\;
\begin{blockarray}{cccccccc}
	&& \mathbf{x}_1 & \mathbf{x}_2 & \cdots & \mathbf{x}_{k-1} & \mathbf{x}_k & \mathbf{x}_{k+1}\\
	\begin{block}{cc[cccccc]}
		\dot{\mathbf{x}}_1 &\quad & \mathbf{Q}_1 & \mathbf{Q}_{12} & \mathbf{0} & \cdots & \cdots & \mathbf{0}\\
		\dot{\mathbf{x}}_2 &\quad & \mathbf{0} & \mathbf{Q}_2 & \mathbf{Q}_{23} & \mathbf{0} & \cdots & \vdots\\
		\vdots &\quad & \vdots & \mathbf{0} & \ddots & \ddots & \mathbf{0} & \vdots\\
		\dot{\mathbf{x}}_{k-1} &\quad & \vdots & \cdots & \mathbf{0} & \mathbf{Q}_{k-1} & \mathbf{Q}_{k-1,k} & \mathbf{0}\\
		\dot{\mathbf{x}}_k &\quad & \mathbf{0} & \cdots & \cdots & \mathbf{0} & \mathbf{Q}_k & \mathbf{Q}_{k,k+1}\\
	\end{block}
\end{blockarray}\;,
\end{equation}
\end{widetext}
where we defined a more compact notation on the right.
In \eqref{eq:xdotblock}, $\mathbf{x}_m$ are vectors with all dynamically relevant motif counts of size $m$, $\mathbf{Q}_m$ the 
coefficients
for $\dot{\mathbf{x}}_m:=d\mathbf{x}_m/dt$ with motifs of the same size, and $\mathbf{Q}_{m,m+1}$ the 
coefficients for $\dot{\mathbf{x}}_m$ with motifs of the size $m{+}1$. The block diagonal
form arises because the change of size-$m$ motif counts depends only on motif counts of size $m$ and of size $m{+}1$. In Appendix \ref{sec:MFSIS3}, 
\eqref{eq:xdotblock} is shown for SIS spreading up to a maximum motif size of $k=3$.

The substitution via conservation relations can be written as
\begin{equation}\label{eq:subcons}
	\mathbf{x}=\mathbf{E}\cdot\tilde{\mathbf{x}}+\mathbf{c},
\end{equation}
where $\mathbf{x}=\mathbf{x}_1,...,\mathbf{x}_k$, $\tilde{\mathbf{x}}$
is $\mathbf{x}$ with the to be substituted elements omitted. $\mathbf{E}$, $\mathbf{c}$
contain the coefficients of linear dependence from (\ref{eq:cons},~\ref{eq:cons2}), with
their $i$th row/element corresponding to the identity transformation for motifs that
are not substituted. Substituting this
into \eqref{eq:xdotblock} results in the system of equations for the remaining
motifs $\tilde{\mathbf{x}}$:
\begin{equation}\label{eq:xdotsubcons1}
	\dot{\tilde{\mathbf{x}}}=\tilde{\mathbf{Q}}_{1\cdots k,1\cdots k}\cdot(\mathbf{E}\cdot\tilde{\mathbf{x}}+\mathbf{c})+\tilde{\mathbf{Q}}_{1\cdots k,k+1}
	\cdot
	\tilde{\mathbf{x}}_{k+1}
\end{equation}
(where the tilde omits the to be substituted rows/elements), such that we can write
the system in the same form as \eqref{eq:xdotblock}, but now with an added 
constant vector:
\begin{equation}\label{eq:xdotsubcons2}
	\dot{\tilde{\mathbf{x}}}=\mathbf{Q}'
	\cdot
	\begin{bmatrix}
		\tilde{\mathbf{x}}\\
		\tilde{\mathbf{x}}_{k+1}\\
	\end{bmatrix}
	+\tilde{\mathbf{c}}.
\end{equation}
For an example, see Section \ref{sec:Second-order}.

\section{Closure scheme}\label{sec:Closure-scheme}
Because in \eqref{eq:xdotsubcons2}, counts for the largest motifs $\mathbf{x}_{k+1}$ appear
on the right-hand side but not on the left-hand side, \eqref{eq:xdotsubcons2} is underdetermined.
A closure scheme provides a way of expressing the
undetermined parts $\tilde{\mathbf{x}}_{k+1}$ in \eqref{eq:xdotblock}
through a nonlinear function $\tilde{\mathbf{x}}_{k+1}\approx \tilde{\mathbf{f}}(\tilde{\mathbf{x}})$, creating
a closed system of ODEs:
\begin{equation}\label{eq:xdotsubconscl}
	\dot{\tilde{\mathbf{x}}}=\mathbf{Q}'
	\cdot
	\begin{bmatrix}
		\tilde{\mathbf{x}}\\
		\tilde{\mathbf{f}}(\tilde{\mathbf{x}})\\
	\end{bmatrix}
	+\tilde{\mathbf{c}},
\end{equation}
where $\tilde{\mathbf{f}}$ also depends on counts of induced subgraphs of order up to $k{+}1$
($[g_1^{k+1}]$ in Figure \ref{fig:overview}) if the motifs were not normalised in advance.
In this section, we develop a closure scheme that decomposes $\tilde{\mathbf{x}}_{k+1}$
into its smaller-sized components.
Our final formula generalises closures hitherto most commonly used, as shown in e.g. \citet{House2009}. 
We will show that the decomposition
is valid when: (\emph{i}) counts of components are conditionally independent given the node states in their intersection and the adjacency
structure between them is a tree,
(\emph{ii}) the network is spatially homogeneous, and (\emph{iii}) the network is sufficiently large, such that the law of large numbers applies. We start with some introducing examples in 
Section \ref{sec:introex} and defer detailed explanation to Sections \ref{sec:closuredefbg},\ref{sec:motifdecomp}. Examples are given in Sections \ref{sec:clex} and \ref{sec:clexs}.

\subsection{Introduction}\label{sec:introex}

When truncating the moment hierarchy at a chosen order $k$, we approximate
the order $k{+}1$ motifs appearing in the equations for order $k$ motifs in terms of
lower-order motifs. E.g., looking at the moment equations for SIS spreading in 
Appendix \ref{sec:MFSIS3}, when truncating at $k{=}1$, we would need
an expression of $[\motIS]$ in terms of $[\motI]$ and $[\motS]$. When assuming
statistical independence of neighbouring node states, the resulting expression
is of the form $[\motIS]\propto [\motI][\motS]$ (ignoring proportionality constants for now). 
Similarly, for
truncation at $k{=}2$, we would need an expression for all the 3-chains and triangles
on the right-hand side of (\ref{eq:mf4},~\ref{eq:mf5}). For instance,
the 3-chain $[\motSSI]$ is typically decomposed as $[\motSSI]\propto [\motSS][\motIS]/[\motS]$.
Using the shorthand $x_i$ for the event $X_i{=}x$, this is justified when there is a conditional independence relation of the form $P(S_i,S_j,I_k){=}P(S_i,S_j)P(S_j,I_k{\mid} S_j)$ and if
the component probabilities are the same everywhere in the network, such that node indices
$i,j,k$ do not matter. This can be generalised to larger chains, such as e.g.
$[\motSSII]\propto[\motSSI][\motSII]/[\motIS]$, which similarly follows from
the assumed conditional independence relation
$P(S_i,S_j,I_k,I_l){=}P(S_i,S_j,I_k)P(S_j,I_k,I_l{\mid} S_j,I_k)$ and homogeneity
in the network. 

It is possible to generalise the 
examples above to larger subgraphs by starting from the chain rule of probability,
\begin{equation*}
	P(\boldsymbol{x})=P(\boldsymbol{x}_{\boldsymbol{i}_n}\mid \boldsymbol{x}_{\boldsymbol{i}_{n-1}},..., \boldsymbol{x}_{\boldsymbol{i}_2},\boldsymbol{x}_{\boldsymbol{i}_1})... P(\boldsymbol{x}_{\boldsymbol{i}_2}\mid \boldsymbol{x}_{\boldsymbol{i}_1})P(\boldsymbol{x}_{\boldsymbol{i}_1}),
\end{equation*}
where $\boldsymbol{x}$ is a vector of state labels on a given network 
motif, and subsequently simplify with assumed conditional independence
relations. For instance, for our second example above, $[\motSSI]$, we have $P(S_i,S_j,I_k){=}P(I_k {\mid} S_i,S_j)P(S_j{\mid} S_i)P(S_i)$. If now node $k$ is conditionally independent of node $i$,
we substitute $P(I_k {\mid} S_i,S_j){=}P(I_k {\mid} S_j)$, such that we obtain
the expression found above. In general, a simplification of the chain
rule in terms of subgraphs is possible if we can order the chosen 
sets of subgraphs (with node indices $\boldsymbol{i}_1,...,\boldsymbol{i}_n$) 
without creating cycles
and if the states of adjacent subgraphs are conditionally independent given 
their shared nodes. To make this precise, we need the concept of
independence map. The
independence map is a graph in which link absence
between two nodes means that there is no direct dependence between them. 
For instance, if in the four-node 
graph shown in row  4 column 1 of Table \ref{tab:clex} there is statistical independence
between node states that are further than two steps removed, its independence map
is the graph shown in row 4 column 3 of Table \ref{tab:clex}. The condition mentioned
above for simplification of the chain rule in terms of subgraphs reduces
to the requirement that the independence map be chordal, because this allows
a tree composition in terms of maximal cliques, as shown in column 4 of 
Table~\ref{tab:clex}. A graph is chordal if every 
cycle greater than three is cut short by a link between two non-consecutive
nodes of the cycle.

\begin{table*}
	\includesvg[inkscapelatex=false,width=2\columnwidth]{./figures/table_closures}
	\caption{Examples of our method to obtain closure formulas. Because
		the closures can be written independent of the state labels, the 
		decomposition is shown for node-indexed graphs without reference to particular
		node states.}\label{tab:clex}
\end{table*}

In reality, we do not know the dependence structure between
the graph nodes. 
However, the practical requirement of truncation of the moment hierarchy
obliges us to assume an independence map for each of the largest motifs. In order to
obtain a consistent decomposition method, we will choose as independence
map the graph obtained by connecting all 
nodes to other nodes in their $d{-}1$ neighbourhood, where $d$ is the
diameter of the motif. The chain rule then leads to a decomposition in
terms of (maximal) cliques of the independence map in the numerator and the node
sets that separate them in the denominator. Table \ref{tab:clex} shows 
example graphs, with their diameter,
their assumed independence maps, their tree decomposition and the resulting 
closure formula. The choice of dependence within a distance $d$ may 
in some cases result in non-chordal
independence maps, such that the decomposition cannot be made. In the 
next sections, we will explain our method in detail and further show
how one can treat motifs with a non-chordal independence map.

\subsection{Definitions and background}\label{sec:closuredefbg}

We will rely on the theory of decomposable Markov networks, following mostly the terminology of \citet[][Ch. 3]{Pearl1988}. 
We generalise the decomposition from a factorisation involving
(1-)cliques to one involving $d$-cliques. he definitions in this section apply to a general
graph ${\cal G}$ with nodes states $X$, but we will apply the decomposition to motifs in Section \ref{sec:motifdecomp}. 

\paragraph{Separation} Given a graph ${\cal G}({\cal V},{\cal E})$ and three disjoint 
subsets of nodes
$\boldsymbol{i},\boldsymbol{j},\boldsymbol{k}\subset {\cal V}$, $\boldsymbol{k}$ \emph{separates}
$\boldsymbol{i}$ and $\boldsymbol{j}$ in ${\cal G}$, written as
$\boldsymbol{i}\indep_{\mkern-4mu{\cal G}}\;\boldsymbol{j}\mid \boldsymbol{k}$, if every
path between $\boldsymbol{i}$ and $\boldsymbol{j}$ has at least one vertex in $\boldsymbol{k}$. 
Here, $\boldsymbol{k}$ is called a separator, or also, a node cut set of $\boldsymbol{i}$ 
and $\boldsymbol{j}$ in ${\cal G}$.

\paragraph{Independence map} An independence map ${\cal M}$ is a graph that represents 
the independence between components of a set of random variables $X$ such that separation
in ${\cal M}$ guarantees conditional independence between corresponding subsets of $X$. 
More precisely, given three disjoint subsets of nodes
$\boldsymbol{i},\boldsymbol{j},\boldsymbol{k}\subset {\cal V}$, $X$ possesses
a spatial Markov property:
\begin{equation}\label{eq:IMap}
\boldsymbol{i}\indep_{\mkern-4mu{\cal M}} \;\boldsymbol{j}\mid \boldsymbol{k}\implies  X_{\boldsymbol{i}}\indep X_{\boldsymbol{j}}\mid X_{\boldsymbol{k}},
\end{equation}
where the $\indep$ notation on the right refers to independence of the random variables: $P(X_{\boldsymbol{i}}{=}x_{\boldsymbol{i}},X_{\boldsymbol{j}}{=}x_{\boldsymbol{j}}{\mid} X_{\boldsymbol{k}}{=}x_{\boldsymbol{k}})=P(X_{\boldsymbol{i}}{=}x_{\boldsymbol{i}}{\mid} X_{\boldsymbol{k}}{=}x_{\boldsymbol{k}})P(X_{\boldsymbol{j}}{=}x_{\boldsymbol{j}}{\mid} X_{\boldsymbol{k}}{=}x_{\boldsymbol{k}})$.
The pair $(X,{\cal M})$ defines what is known as a \emph{Markov network}. 

\paragraph{Independence beyond distance $d$} Let ${\cal G}=({\cal V},{\cal E})$ be a
graph where the nodes in ${\cal V}$ have (random) states $X$. We define ${\cal G}^d$ as the graph in 
which all nodes of ${\cal G}$ are neighbours if they are at most a shortest distance $d$ 
away from each other, i.e.
\begin{equation}\label{eq:Gk}
{\cal G}^d{=}({\cal V},{\cal E}^d),\;\; \mathrm{with}\;\;{\cal E}^d{=}\{(i,j){\in}{\cal V}: \d_{\cal G}(i,j){\leq} d\}\mbox{.}
\end{equation}
We then say that $({\cal G},X)$ has \emph{independence beyond distance $d$} if ${\cal G}^d$ is the independence map of $X$, or for all distinct subsets of nodes $\boldsymbol{i},\boldsymbol{j},\boldsymbol{k}\subset {\cal V}$ holds
\begin{equation}\label{eq:indepk}
\boldsymbol{i}\indep_{\mkern-4mu{\cal G}^d} \;\boldsymbol{j}\mid \boldsymbol{k}\implies  X_{\boldsymbol{i}}\indep X_{\boldsymbol{j}}\mid X_{\boldsymbol{k}}.
\end{equation}
This means that states of two non-neighbouring sets in ${\cal G}^d$, which by definition
\eqref{eq:Gk} are further than $d$ steps apart in ${\cal G}$, are independent of each
other given the state of their separator $\boldsymbol{k}$. 

\paragraph{Maximal $d$-cliques and $d$-clique graph}  A maximal clique is a complete subgraph not contained in a larger complete subgraph \cite{Harary1973}. As a generalisation, \emph{maximal 
$d$-cliques} are maximal subgraphs with distance between any two nodes 
not greater than $d$ \cite{Mathematica}. Correspondingly, maximal cliques in ${\cal G}^d$ are maximal $d$-cliques in ${\cal G}$. The graph ${\cal C}^d$ is the
\emph{$d$-clique graph} of ${\cal G}$ if
each node in ${\cal C}^d$ corresponds to a $d$-clique in ${\cal G}$
with links between nodes in ${\cal C}^d$ occurring when the corresponding $d$-cliques
overlap. Hence, while nodes in ${\cal C}^d$ correspond to maximal $d$-cliques in
${\cal G}$, links in ${\cal C}^d$ correspond to intersections between overlapping
maximal $d$-cliques in ${\cal G}$.

\paragraph{Junction graph of $d$-cliques} A junction graph of
${\cal G}$'s maximal $d$-cliques, denoted further as ${\cal J}^d({\cal G})$, is a subgraph of the $d$-clique graph ${\cal C}^d$ obtained by
removing redundant links from ${\cal C}^d$. Denoting $d$-cliques 
corresponding to nodes $i,j$ in ${\cal C}^d$ as $c_i,c_j\subset{\cal V}$, a link between
$i$ and $j$ in ${\cal C}^d$ is \emph{redundant} when there is an alternative path between $i$ 
and $j$ in ${\cal C}^d$ passing by a series of other nodes in ${\cal C}^d$ of which the corresponding 
$d$-cliques all contain $c_i\cap c_j$. The junction graph ${\cal J}^d({\cal G})$ is then obtained
by iteratively removing redundant links from the ${\cal C}^d$ until there are no further redundant links. While the $d$-clique
graph is unique, there may be several junction graphs ${\cal J}^d({\cal G})$ of $d$-cliques for one graph ${\cal G}$.
Note that for chordal graphs (defined below), the junction graph equals what is known as a junction
tree, which can also be obtained via the junction tree algorithm \cite{barber2012},
applied to the $d$-clique graph.

\paragraph{$d$-chordality} A graph ${\cal G}$ is chordal when for every cycle of length greater than $3$, there exists a link in ${\cal G}$ between two non-consecutive nodes of the cycle (thus, giving a short-cut, also called \emph{chord} to the cycle). As a generalisation, we will call a graph ${\cal G}$ $d$-chordal 
if ${\cal G}^d$ is chordal. If a graph ${\cal G}$ is $d$-chordal, then 
${\cal J}^d({\cal G})$ is a tree, or equivalently, if ${\cal G}^d$ is chordal,
then ${\cal J}({\cal G}^d)$ is a tree. Non-chordal graphs
can always be converted to a chordal graph via \emph{triangulation},
i.e. adding chords to every chordless cycle of length greater than $3$.
We will write below $\text{tr}({\cal G}^d)$ as a \emph{minimal triangulation} of
a non-chordal ${\cal G}^d$, obtained by adding the smallest number of links
that leads to chordality, unless stated otherwise. 

\paragraph{Decomposability at distance $d$} If there is independence
beyond distance $d$ \eqref{eq:indepk} and ${\cal G}$ is $d$-chordal,
then the joint probability of the network nodes of ${\cal G}$ being in a
given state $P(X=\boldsymbol{x})$ can be factorised over the
$d$-cliques of ${\cal G}$.  We will call this property of the graph
${\cal G}$ and its node states $X$ \emph{decomposability at distance
$d$}. We call ${\cal J}$ the set of
$d$-cliques in ${\cal G}$, ordering its elements consistent with the
resulting junction tree structure for ${\cal G}^d$ (such that
${\cal J}_1$ is the chosen root and parent nodes have lower index
than their leaves), and call $\pa({\cal J}_i)$ the parent node of
${\cal J}_i$. In such cases, the factorisation is possible because
the tree structure between $d$-cliques allows application of the chain
rule of conditional probability.  
\begingroup \allowdisplaybreaks
\begin{align}\label{eq:dec:chord}
P(X{=}\boldsymbol{x})&=\prod_{i=1}^{|{\cal J}|}P(X_{{\cal J}_i}{=}\boldsymbol{x}_{{\cal J}_i}\mid X_{\pa({\cal J}_i)}{=}\boldsymbol{x}_{\pa({\cal J}_i)}),\nonumber\\
&=\prod_{i=1}^{|{\cal J}|}P(X_{{\cal J}_i}{=}\boldsymbol{x}_{{\cal J}_i}\mid X_{\pa({\cal J}_i)\cap {\cal J}_i}{=}\boldsymbol{x}_{\pa({\cal J}_i)\cap {\cal J}_i}),\nonumber\\
&=\frac{\prod_{i=1}^{|{\cal J}|}P(X_{{\cal J}_i}{=}\boldsymbol{x}_{{\cal J}_i})}{\prod_{i=2}^{|{\cal J}|}P(X_{\pa({\cal J}_i)\cap {\cal J}_i}{=}\boldsymbol{x}_{\pa({\cal J}_i)\cap {\cal J}_i})}.
\end{align}
\endgroup
The steps in \eqref{eq:dec:chord} are explained as follows. As $d$-chordality
makes ${\cal J}^d({\cal G})$ a tree and any $d$-clique
separates its neighbours, one can recursively use conditional independence of children given parents (line 1). We use the convention that $\pa({\cal J}_1)=\emptyset$, such that the first factor is $P(X_{{\cal J}_1}{=}\boldsymbol{x}_{{\cal J}_1})$. In line 2 we exploit that any two $d$-cliques in 
${\cal G}$ are also separated by their intersection to condition instead on intersections.

\paragraph{Non-$d$-chordal graphs} There are two alternative ways to decompose non-$d$-chordal ${\cal G}$: (\emph{i}) perform the decomposition \eqref{eq:dec:chord} 
on the (more conservative) independence map after triangulation, $\text{tr}({\cal G}^d)$. 
In this case, the factors in \eqref{eq:dec:chord} may still contain subgraphs of diameter $d{+}1$, but of smaller size than ${\cal G}$, such that one may have to apply \eqref{eq:dec:chord} recursively to achieve smaller diameter for all factors. Furthermore, the resulting decomposition will depend on the choice of triangulation. 
Alternatively, (\emph{ii}), one can
start from the non-tree ${\cal J}^d({\cal G})$ and use the ad-hoc formula (without prior triangulation)
\begin{equation}\label{eq:dec:nchord}
P(X{=}\boldsymbol{x})\approx\zeta\frac{\prod_{i}^{|{\cal J}|}P(X_{{\cal J}_i}{=}\boldsymbol{x}_{{\cal J}_i})}{\prod_{i,j\ne i}^{|{\cal J}|}P(X_{{\cal J}_i\cap {\cal J}_j}{=}\boldsymbol{x}_{{\cal J}_i\cap {\cal J}_j})}.
\end{equation}
Because the fraction in \eqref{eq:dec:nchord} does not result from
application of the chain rule as in \eqref{eq:dec:chord}, it is not a
product of conditional probabilities and hence it does not guarantee
the property that each of the node states in $\boldsymbol{x}$ has to
appear one more time in the numerator than in the denominator, which
in turn leads to inconsistency between closure formulas that assume
different $d$. The factor $\zeta$ in \eqref{eq:dec:nchord}
corrects for this inconsistency -- see Section \ref{sec:motifdecomp} for more
detail.  After applying (\emph{i}) to non-chordal graphs, the
nodes in the resulting $d$-clique tree are not all maximal
$d$-cliques of ${\cal G}$ any more \footnote{More precisely, they
consist of the nodes of the maximal cliques of
$\text{tr}({\cal G}^d)$ and the links of ${\cal G}$.}. 
 When applying
(\emph{ii}) to non-chordal graphs, the subgraphs in
${\cal J}^d({\cal G})$ are still maximal $d$-cliques of ${\cal G}$ but
the clique graph is not a tree, thus, violating the assumptions
  behind the decomposition \eqref{eq:dec:chord}.
  
\paragraph{Maximum motif diameter} The decomposition explained in this
section implies that if independence beyond distance $d$ is valid
in the whole network $({\cal G},X)$, then we know that we 
only need to consider motifs up to diameter $d$, justifying truncation
of the moment hierarchy. Note, however, that truncation is usually
done at a give size, not at a given diameter. We expand on this 
in point \ref{pt:exactness}.

\subsection{Motif decomposition}\label{sec:motifdecomp}
\paragraph{Decomposition of motifs at the individual level} We
  apply the decomposition \eqref{eq:dec:chord} to motifs with
connectivity $\mathsf{a}$ and chordal independence map
$\mathsf{a}^{d}$ embedded in the network. For now, we ignore that this
may in some cases break the independence assumption, but see section
\ref{pt:exactness} for more detail on this issue. Considering a set
  $\boldsymbol{i}$ of nodes that have connectivity $\mathsf{a}$ in
our network ${\cal G}$ and taking
$P(X_{\boldsymbol{i}}{=}\boldsymbol{x})$ as the probability that these
are in states with labels $\boldsymbol{x}$, we can write the
decomposition  \eqref{eq:dec:chord} for of $P(X_{\boldsymbol{i}}{=}\boldsymbol{x})=\langle\left[\boldsymbol{x}_{\boldsymbol{i}}^\mathsf{a}\right]\rangle$
to obtain 
\begin{align}\label{eq:dec:moti:chord}
\langle\left[\boldsymbol{x}_{\boldsymbol{i}}^\mathsf{a}\right]\rangle=\frac{\prod_{j=1}^{|{\cal J}|}\langle[\boldsymbol{x}_{\boldsymbol{i}_{{\cal J}_j}}]\rangle}{\prod_{j=2}^{|{\cal J}|}\langle[\boldsymbol{x}_{\boldsymbol{i}_{\pa({\cal J}_j)\cap {\cal J}_j}}]\rangle},
\end{align}
where now ${\cal J}$ is the set of $d$-cliques in $\mathsf{a}$. Choosing 
$d=\text{diam}(\mathsf{a})-1$ ensures that the decomposition results in component motifs with
diameter decreased by one compared to the decomposed motif. For motifs with non-chordal $\mathsf{a}^{d}$ one
can, as noted above, either triangulate $\mathsf{a}^{d}$ first
or use the \emph{ad-hoc} approximation \eqref{eq:dec:nchord}. 
We relied on the package \texttt{Chordal Graph} \cite{Bulatov2011} for 
triangulation of non-chordal $\mathsf{a}^{d}$. The 
ad-hoc formula \eqref{eq:dec:nchord} applied to the motif at $\boldsymbol{i}$  is
\begin{equation}\label{eq:dec:moti:nchord}
\langle[\boldsymbol{x}_{\boldsymbol{i}}^\mathsf{a}]\rangle\approx\frac{\prod_{j}^{|{\cal J}|}\langle[\boldsymbol{x}_{\boldsymbol{i}_{{\cal J}_j}}]\rangle}{\prod_{j,k\ne j}^{|{\cal J}|}\langle[\boldsymbol{x}_{\boldsymbol{i}_{{\cal J}_j\cap {\cal J}_k}}]\rangle}\prod_{j\in\boldsymbol{i}}\langle[ x_j]\rangle^{\textstyle\gamma_j}.
\end{equation}
The consistency correction (written as $\zeta$ in \eqref{eq:dec:nchord}) here equals $\prod_{j\in\boldsymbol{i}}\langle\left[x_j\right]\rangle^{\gamma_j}$ 
\bwr{\jsq{if one uses the notation with links ${\cal L}_\ell$ of the junction graph, is then
  \begin{align*}
    \gamma_\nu=1+\sum_{\ell=1}^{|{\cal L}|}\mathds{1}_{\nu\in{\cal L}_\ell}-\sum_{j=1}^{|{\cal J}|}\mathds{1}_{\nu\in{\cal J}_j},
  \end{align*} where $\mathds{1}_e$ is the logical indicator function, equalling $1$ if the logical expression $e$ is true, $0$ otherwise?}}
and ensures that the ad-hoc extension of closures to motifs with non-chordal $\mathsf{a}^{d}$
does not result in inconsistency with MF$1$ \cite[condition 1 of][Ch. 21]{Dieckmann2000} 
under independence between node states: when
all motifs of order greater than one are replaced by products of order
one motifs, i.e.
$\langle\left[\boldsymbol{x}_{\boldsymbol{i}_\cdot}^\mathsf{a}\right]\rangle\to\prod_{j\in\boldsymbol{i}_\cdot}\langle\left[x_j\right]\rangle$, the right hand side of (\ref{eq:dec:moti:nchord}) should reduce to MF$1$. Therefore, for each $j$, $\gamma_j$ is chosen such that this is fulfilled. These ad-hoc steps usually result in violation of the
conservation relations of Section \ref{sec:alg:iv}
\cite{Dieckmann2000}. In
the approximations used in Section~\ref{sec:SIS}, the bias introduced due to this violation is small. For mitigation of this problem, see
\cite{Dieckmann2000,Peyrard2008}. 

\paragraph{Decomposition of motifs at the population level}  If we take the following \emph{spatial homogeneity} assumption for all motifs $\boldsymbol{x}^\mathsf{a}$ of sizes $m$ up to our maximal considered size $k$
\begin{equation}\label{eq:spatialhom}
\forall \boldsymbol{i},\boldsymbol{i'}\in {\cal I}(\mathsf{a}): \langle\left[\boldsymbol{x}_{\boldsymbol{i}}^{\mathsf{a}}\right]\rangle=\langle\left[\boldsymbol{x}_{\boldsymbol{i'}}^{\mathsf{a}}\right]\rangle=\frac{\langle\left[\boldsymbol{x}^{\mathsf{a}}\right]\rangle}{\left[\mathsf{a}\right]}=\langle\llbracket\boldsymbol{x}^{\mathsf{a}}\rrbracket\rangle,
\end{equation}
where ${\cal I}(\mathsf{a}):=\{\boldsymbol{j}\in S(m,N):\mathsf{A}_{\boldsymbol{j}}=\mathsf{a}\}$,
then (\ref{eq:dec:moti:chord},~\ref{eq:dec:moti:nchord}) are independent of $\boldsymbol{i}$, such that
we can write \eqref{eq:dec:moti:chord} as 
\begin{align}\label{eq:dec:mot:chord}
\langle\llbracket\boldsymbol{x}^\mathsf{a}\rrbracket\rangle\approx\frac{\prod_{j=1}^{|{\cal J}|}\langle\llbracket\boldsymbol{x}_{{\cal J}_j}\rrbracket\rangle}{\prod_{j=2}^{|{\cal J}|}\langle\llbracket\boldsymbol{x}_{\pa({\cal J}_j)\cap {\cal J}_j}\rrbracket\rangle},
\end{align}
and \eqref{eq:dec:moti:nchord} as
\begin{equation}\label{eq:dec:mot:nchord}
\langle\llbracket\boldsymbol{x}^\mathsf{a}\rrbracket\rangle\approx\frac{\prod_{j}^{|{\cal J}|}\langle\llbracket\boldsymbol{x}_{{\cal J}_j}\rrbracket\rangle}{\prod_{j,k}^{|{\cal J}|}\langle\llbracket\boldsymbol{x}_{{\cal J}_j\cap {\cal J}_k}\rrbracket\rangle}\prod_{p}^{m}\langle\llbracket x_p\rrbracket\rangle^{\textstyle\gamma_p},
\end{equation}
which may be used to close the population-level equations \eqref{eq:xdotblock}.
In (\ref{eq:dec:mot:chord},~\ref{eq:dec:mot:nchord}), node indexing of a motif is consistent with the node labels $\boldsymbol{x}$. Recall that we use a single consistent indexing for isomorphic motifs. The decomposition (\ref{eq:dec:mot:chord},~\ref{eq:dec:mot:nchord}) is not unique. If the independence assumptions are satisfied each of the alternative ways to decompose motif $\boldsymbol{x}^\mathsf{a}$ should result in the same value. As we do not expect the independence to be perfectly valid, we
take the average of the alternative ways of decomposing $\boldsymbol{x}^\mathsf{a}$ if they exist.

\paragraph{Normalisation} We showed the closure formulas for
normalised motifs, i.e. we first normalised the counts of motifs 
via \eqref{eq:cons:norm}, and then applied the closure. Hence, in this case, the counts
of induced subgraphs in the network enter into the system of equations as normalisation factors in the unclosed system. One can also decide not to normalise (or to do it after applying closure). 
In this latter case, the subgraph counts
enter into the final system of equations when applying closure, as the 
closure formulas for the non-normalised motif counts contain them 
(to see this, substitute each motif count in (\ref{eq:dec:mot:chord},~\ref{eq:dec:mot:nchord})
as $\llbracket\boldsymbol{y}^{\mathsf{b}}\rrbracket\to[\boldsymbol{y}^{\mathsf{b}}]/[\mathsf{b}]$). Hence, structural information specific to the considered network enters the mean field equations either when normalising the motif counts or when applying closure. In simple cases,
such as lattices or random graphs,
the subgraph counts can be found by hand without much effort. In other cases, one can resort to
subgraph counting algorithms -- we used \texttt{IGraph} \cite{Horvat2020} for \texttt{Mathematica} \cite{Mathematica}.

\paragraph{Law of large numbers} We will use the population-level closure to study the 
steady states in a single realisation of a given network. Motif counts are then assumed to be 
the total counts in a single network, instead of their expectations over many
realisations. As the closure formulas apply to expectations, we make the additional assumption that motif counts are close to their expectations. 

\paragraph{Additional bias}\label{pt:exactness} Decomposing 
the whole network ${\cal G}$ at distance $d$
is exact when ${\cal G}$ is $d$-chordal and there is independence beyond
distance $d$ (Section \ref{sec:closuredefbg}). Applying the decomposition to motifs
embedded in ${\cal G}$ instead of to ${\cal G}$ can be done without additional
bias when $\mathsf{a}$
is a distance-hereditary subgraph of ${\cal G}$ (i.e. distances between nodes
in $\mathsf{a}$ are equal to those between corresponding nodes in ${\cal G}$) and
conditional independence relations implied by $\mathsf{a}$ are also valid in ${\cal G}$. 
As a counterexample for the former, take for
${\cal G}$ the six-node graph $\motlatsix$ and for $\mathsf{a}$ its induced subgraph consisting of nodes $\{1,2,4,6,5\}$. Here, $\mathsf{a}$
is not distance-hereditary because $\d_\mathsf{a}(1,5){=}4\neq\d_{\cal G}(1,5){=}2$.
As a counterexample for the latter, take for ${\cal G}$ the square $\motlatfour$
and for $\mathsf{a}$ the 3-node chain $\{1,2,4\}$. In this case, $\mathsf{a}$ is distance hereditary, but, while (when assuming independence beyond distance $d{=}1$) within $\mathsf{a}$ we have the independence relation $\{1\}\indep_{\mkern-5mu\mathsf{a}}\;\{4\}\mid\{2\}$, this is not true in 
${\cal G}$, where $\{1\}\not\indep_{\mkern-5mu{\cal G}}\;\{4\}\mid\{2\}$, because 
the node cut set in ${\cal G}$ for $\{1\}$ and $\{4\}$ is $\{2,3\}$. This occurs
because the decomposed motifs are non-maximal $d$-cliques. Therefore, we expect a bias as a consequence of this in the closed mean-field equation hierarchy
(
\eqref{eq:xdotblock} with \eqref{eq:dec:mot:chord} or \eqref{eq:dec:mot:nchord} at population level under spatial homogeneity). This bias can be avoided when 
expressing the equations in terms of maximal $p$-cliques for $p\in\{0,...,k{+}1\}$ 
and truncating at given diameter instead of at given size. 

\subsection{Examples}\label{sec:clex}

Appendix \ref{sec:clexs} shows $13$ application examples of (\ref{eq:dec:mot:chord},~\ref{eq:dec:mot:nchord}) in table form. 
As the closures can be
written independent of the particular labels, they are shown for subgraphs only, with each node tagged with its index. We have
also dropped the $\langle\cdot\rangle$, assuming that the law of large numbers applies, such that the counts approach their expectations almost surely for increasing network size $N$. The examples can be understood 
by reading the table from left to right. 
Below, we derive
the normalisation factors and the non-normalised closures of examples 1-3 of
Appendix \ref{sec:clexs} for different network types. Note that, unlike in 
Appendix \ref{sec:clexs}, we use letter labels below, for consistency with the main
text and the literature. 

\begin{enumerate}[wide, labelwidth=!]
	\item $\boldsymbol{x}^{\mathsf{a}}=\motabc$:\enspace This diameter-$2$ motif has 
	chordal independence map equal to $\mathsf{a}$ and decomposes with \eqref{eq:dec:mot:chord} as
		\begin{equation}
			\llbracket\motabc\rrbracket\approx\frac{\llbracket\motab\rrbracket\llbracket\motbc\rrbracket}{\llbracket\motb\rrbracket},\label{eq:mclosabcn}
		\end{equation}
	assuming conditional independence beyond distance $d{=}1$.
	Via normalization \eqref{eq:cons:norm} we obtain also the closure for the non-normalized counts:
		\begin{equation}
			\left[\motabc\right]\approx\left[\motabcNL\right]\frac{\left[\motaNL\right]}{\left[\motabNL\right]^2}\frac{\left[\motab\right]\left[\motbc\right]}{\left[\motb\right]}.\label{eq:mclosabc}
		\end{equation}
	The counts of the induced subgraphs of size $2$ and $1$, required for normalisation, are
		\begin{align}
			\left[\motabNL\right]=\kappa N,\quad\left[\motaNL\right]=N,\label{eq:munl123}
		\end{align}
	and total number of triples ($3$-node motifs) in the network is
		\begin{equation}
			\left[\motabcNL\right]+\left[\mottr\right]=\sum_i^N \kappa_i(\kappa_i-1),\label{eq:triples}
		\end{equation}
	with $\kappa_i$ the number of neighbours of node $i$ and $\kappa$ the mean number
	of neighbours over the whole network. For particular network types \eqref{eq:triples} can be simplified. Below are two examples.
	\begin{enumerate}
		\item For a network with fixed degree without triangles (e.g. a square lattice), we have $\forall i:\kappa_i=\kappa$ and $\left[\mottr\right]=0$, 
		such that
			\begin{equation}
				\left[\motabcNL\right]=\kappa(\kappa-1)N.\label{eq:motUNL3hom}
			\end{equation}
		Using \eqref{eq:mclosabc} and \eqref{eq:motUNL3hom}, we obtain
			\begin{equation}
				\left[\motabc\right]\approx\frac{\kappa-1}{\kappa}\frac{\left[\motab\right]\left[\motbc\right]}{\left[\motb\right]}.\label{eq:mclosabchom}
			\end{equation}
		An early use of this closure for networks can be found in \citet{Keeling1997}.
		\item In a large Erd\H{o}s-R{\'e}nyi random network, we have $\kappa_i\sim \text{Pois}(\kappa)$ and $[\mottr]/([\mottr]+[\motabcNL])\approx0$ \footnote{For an ER random network, the expected number of $m$-node cycles and $m$-node chains (counting all ordered $m$-tuples without repetition) are  $L(N)=\frac{N!}{(N-m)!}(\frac{\kappa}{N-1})^{m}$
			and  $C(N)=\frac{N!}{(N-m)!}(1-\frac{\kappa}{N-1})(\frac{\kappa}{N-1})^{m-1}$,
			where $\frac{\kappa}{N-1}$ is the probability of having a link between
			two given nodes. Hence, $\mathrm{lim}_{N\to\infty}L(N)=\kappa^m=O(1)$ and $\mathrm{lim}_{N\to\infty}C(N)=\mathrm{lim}_{N\to\infty}\kappa^{m-1}N=O(N)$. Hence, as $L(N)$ remains finite and $C(N)$ grows with $N$, cycles can be ignored in the limit of large $N$. In the case $m=3$ (and $N$ large), we can hence safely assume that all triples are chains.}. Hence
			\begin{align}
				\left[\motabcNL\right]=&\sum_i^N\kappa_i^2-\sum_i^N\kappa_i,\nonumber\\
				=&N\left[\mathrm{E}(\kappa_i^2)-\mathrm{E}(\kappa_i)\right],\nonumber\\
				=&N\left[\mathrm{Var}(\kappa_i)+\mathrm{E}(\kappa_i)^2-\mathrm{E}(\kappa_i)\right],\nonumber\\
				=&\kappa^2N,\label{eq:motUNL3ER}
			\end{align}
		where replacing the average by the expectation on the second line requires $N\to\infty$ (law of large numbers), on the third line we used $\mathrm{Var}(\cdot):=\mathrm{E}((\cdot)^2)-(\mathrm{E}(\cdot))^2$, and on the fourth line we used that, for $\kappa_i\sim \text{Pois}(\kappa)$, we have $\mathrm{E}(\kappa_i)=\mathrm{Var}(\kappa_i)=\kappa$. We could also have obtained this result directly from the large-$N$ limit of chains \cite{Note1}.
		Using \eqref{eq:mclosabc} and \eqref{eq:motUNL3ER}, we obtain
			\begin{equation}
				\left[\motabc\right]\approx\frac{\left[\motab\right]\left[\motbc\right]}{\left[\motb\right]},\label{eq:mclosabcER}
			\end{equation}
		This closure was, to the best of our knowledge, first used for networks in \citet{Gross2006}.
	\end{enumerate} 
	\item $\boldsymbol{x}^{\mathsf{a}}=\motabctr$:\enspace This diameter-$1$ motif 
	can be decomposed into its three $0$-cliques as
		\begin{equation}
			\llbracket\motabctr\rrbracket\approx\llbracket\mota\rrbracket\llbracket\motb\rrbracket\llbracket\motc\rrbracket,
		\end{equation}
	when assuming independence of nodes ($d{=}0$). With non-normalised counts, this becomes
		\begin{equation}\label{eq:MFtria}
			\llbracket\motabctr\rrbracket\approx\left[\mottr\right]\left[\mota\right]\left[\motb\right]\left[\motc\right]/\left[\motaNL\right]^3.
		\end{equation}
	Alternatively, one can extend the usage of the ad-hoc formula \eqref{eq:dec:mot:nchord} 
	to include non-maximal cliques: using its
	three $1$-cliques in \eqref{eq:dec:mot:nchord}, we obtain
		\begin{equation}
			\llbracket\motabctr\rrbracket\approx\frac{\llbracket\motab\rrbracket\llbracket\motbc\rrbracket\llbracket\motca\rrbracket}{\llbracket\mota\rrbracket\llbracket\motb\rrbracket\llbracket\motc\rrbracket},\label{eq:PAtriaN}
		\end{equation}
	which is known as the Kirkwood closure for triangles \cite{Sharkey2008}. Using \eqref{eq:cons:norm}, we
	obtain for the closure with non-normalised counts 
		\begin{equation}
			\left[\motabctr\right]\approx\left[\mottr\right]\frac{\left[\motaNL\right]^3}{\left[\motabNL\right]^3}\frac{\left[\motab\right]\left[\motbc\right]\left[\motca\right]}{\left[\mota\right]\left[\motb\right]\left[\motc\right]}.\label{eq:PAtria}
		\end{equation}
	The frequency $\left[\mottr\right]$ depends on
	the network type. For instance, if we use the definition of the clustering
	coefficient $\phi:=\left[\mottr\right]/(\left[\mottr\right]+\left[\motabcNL\right])$ \cite{Keeling1999}, we have [via \eqref{eq:munl123} and \eqref{eq:triples}] for a
	network with fixed degree $\left[\mottr\right]=\phi\kappa(\kappa-1)N$, such that 
		\[
		\left[\motabctr\right]\approx\frac{\kappa-1}{\kappa^2}\phi N\frac{\left[\motab\right]\left[\motbc\right]\left[\motca\right]}{\left[\mota\right]\left[\motb\right]\left[\motc\right]},
		\]
	which was first used for networks by \citet{Keeling1999}.
	\item $\boldsymbol{x}^{\mathsf{a}}=\motabcdtre$:\enspace This diameter-$2$ motif has 
	chordal independence map $\mathsf{a}$ and decomposes with \eqref{eq:dec:mot:chord} as
		\begin{equation}\label{eq:clstar1}
			\llbracket\motabcdtre\rrbracket\approx\frac{\llbracket\motab\rrbracket\llbracket\motbc\rrbracket\llbracket\motbd\rrbracket}{\llbracket\motb\rrbracket^2},
		\end{equation}
	when assuming independence beyond distance $d{=}1$.
	Alternatively, extending the ad-hoc formula \eqref{eq:dec:mot:nchord} to the three
	non-maximal $3$-cliques, we obtain
		\begin{equation}
			\llbracket\motabcdtre\rrbracket\approx\frac{\llbracket\motabc\rrbracket\llbracket\motabd\rrbracket\llbracket\motcbd\rrbracket\llbracket\motb\rrbracket}{\llbracket\motab\rrbracket\llbracket\motbc\rrbracket\llbracket\motbd\rrbracket},\label{eq:clstar2}
		\end{equation}
	where a consistency correction $\llbracket\motb\rrbracket$ was required. The non-normalised form of this closure  was first used in \citet{House2009}.
\end{enumerate}

\section{Application to SIS epidemic spreading}\label{sec:SIS}
We apply our method to SIS spreading, which is a
continuous-time discrete-state Markov chain description of epidemic
spreading through a population of susceptibles \cite[e.g.][]{Kiss2017}. As introduced in Section~\ref{sec:ds-ct-mc}, we have $n=2$ species. The $2\times2$ matrix $\mathbf{R}^0$ of
spontaneous conversion rates and the $2\times2\times2$ tensor $\mathbf{R}^1$ of
conversion rates due to nearest-neighbour interaction for SIS
spreading have only two positive entries, $R^0_{2,1}{=}\gamma$, $R^1_{1,2,2}{=}\beta$, corresponding to reaction scheme \eqref{eq:reactionSIS}.
Hence, contagion of susceptibles $S$ occurs over $IS$ links
at rate $\beta$, whereas recovery occurs spontaneously at rate $\gamma$.
In the study of phase transitions and interacting particle systems,
SIS epidemic spreading is known as the contact process \cite{Harris1974},
which is typically studied on $\mathsf{d}$-dimensional lattices.
In this context, it was found to belong to the directed percolation
universality class, of which scaling properties have been widely studied
\cite{Marro1999,Henkel2008,Tome2015}.

We run the simulations with a Gillespie algorithm \cite{Gillespie2007} and stabilise them
via feedback control \cite{SGNWK08,schilder2015experimental,barton2017control}, 
such that steady states can be obtained in a 
more efficient manner than when running regular simulations (see 
Appendix \ref{sec:feedback}).
  
\subsection{Mean-field equations}\label{sec:MFSIS}
We derived the mean-field models up to fifth order for the square lattice
and up to second order for other networks, including cubic/hypercubic lattices, 
random regular networks and Erd\H{o}s-R{\'e}nyi random networks. In the main text, 
we only show a step-by-step 
derivation of the first and second-order mean-field models because they allow
demonstration of our method in the simplest form. Recall that we write in the text $\langle\left[\cdot\right]\rangle$ as $\left[\cdot\right]$, assuming the LLN holds (in the \texttt{Mathematica} file for MF4 in \ref{sec:MFhoeqs}, we use the notation
$\langle\cdot\rangle$ instead). 

To gain insight in the strength of dependence between neighbouring nodes in simulations and higher-order mean-field models, we will observe the correlation between neighbouring node states $a$ and $b$ as in \citet{Keeling1999}, defined by
\begin{equation}\label{eq:corrs0}
C_{ab}=\frac{N}{\kappa}\frac{\left[\motab\right]}{\left[\mota\right]\left[\motb\right]},
\end{equation}
(where $\kappa$ is the mean degree) or when motif counts are normalised, 
\begin{equation}\label{eq:corrs}
C_{ab}=\frac{\llbracket\motab\rrbracket}{\llbracket\mota\rrbracket\llbracket\motb\rrbracket}.
\end{equation}
They are uncentered correlations between species types that are separated
by one link. Values greater than 1 indicate clustering and values
less than 1 avoidance (compared to a uniform random distribution). For
a generalisation of this correlation to arbitrary distances between
end nodes, see Appendix \ref{sec:dcorrs}.

\subsubsection{MF1 \label{sec:First-order}}

The first-order mean field originates from the molecular field approximation
in statistical physics \cite{weiss1907,bragg1934} and is now commonly known 
as the `mean field model' \cite{Marro1999,Henkel2008,Tome2015,porter2016,Kiss2017,newman2018}.
MF1 only considers node states ($k=1$) and neglects
correlations beyond distance 0. It provides a picture of the dynamics when species
are well mixed throughout a large domain. One way to achieve this is when the domain 
is a complete network on which susceptibles and infecteds have contact rate 
$\kappa\beta/N$, and when $N\to\infty$ \cite{Kiss2017}.

There are two
motif types of size 1: $\motS$ and $\motI$. Of these, only $\motI$ is dynamically
relevant. This means that we only need the equation:
\[
\ddt[\motI]=\beta[\motIS]-\gamma[\motI].
\]
To close the system at order 1, $\left[\motIS\right]$
needs to be expressed in terms of $\left[\motI\right]$ and $\left[\motS\right]$.
No correlation beyond distance 0 corresponds to (using \eqref{eq:dec:mot:chord}):
\begin{equation}\label{eq:clos1sis}
\frac{\left[\motIS\right]}{\kappa N}=\frac{[\motI]}{N}\frac{[\motS]}{N}\implies\left[\motIS\right]=\frac{\kappa}{N}[\motI](N-\left[\motI\right]),
\end{equation}
where we have also used the conservation relation $\left[\motS\right]+\left[\motI\right]=N$
to substitute $\left[\motS\right]=N-\left[\motI\right]$. 
The final expression for the first-order mean field is
\[
\ddt[\motI]=\beta\frac{\kappa}{N}[\motI](N-\left[\motI\right])-\gamma[\motI],
\]
which, after normalisation [via \eqref{eq:cons:norm}], yields 
  \begin{align}
\label{eq:mf1:sis}
\ddt\llbracket\motI\rrbracket=\beta\kappa\llbracket\motI\rrbracket(1-\llbracket\motI\rrbracket)-\gamma\llbracket\motI\rrbracket.
  \end{align}
The steady state solutions are then 
\[
\llbracket\motI\rrbracket_{1}^{*}=0,\llbracket\motI\rrbracket_{2}^{*}=1-\frac{\gamma}{\kappa\beta}.
\]
At $\beta/\gamma=\kappa^{-1}$, the solution
$\llbracket\motI\rrbracket_{1}^{*}$ becomes unstable due to a transcritical
bifurcation, also known as the epidemic threshold in epidemiology. 

\subsubsection{MF2\label{sec:Second-order}}
The second-order mean field originates from the Bethe approximation
in statistical physics \cite{bethe1935} and is now commonly known 
as the `pair approximation' \cite{porter2016,Kiss2017,newman2018,Matsuda1992,Keeling1997,Rand1999,Dieckmann2000,Kefi2007a,Gross2006}.
MF2 neglects dependence beyond distance 1 
and is obtained by considering all dynamically relevant motifs up to size 2. Noting
that motifs without infecteds are dynamically irrelevant and omitting zero blocks,
we obtain
\begin{equation*}\label{eq:SISmf2}
	\begin{blockarray}{ccccccccc}
		&& [\motI] & [\motIS] & [\motII] & [\motSSI] & [\motISI] & [\motSSItr] & [\motSIItr]\\
		\begin{block}{cc[c|cc|cccc]}\relax
			\dot{[\motI]} & & -\gamma & \beta & 0 &  &  &  & \\\cline{3-9}
			\dot{[\motIS]} & & & -\beta-\gamma & \gamma & \beta & -\beta & \beta & -\beta\\
			\dot{[\motII]} & & & 2\beta &-2\gamma & 0 & 2\beta & 0 & 2\beta\\	
		\end{block}
	\end{blockarray}\;.
\end{equation*}
Hence, two types of order 3 motifs appear on the right-hand side: chains and triangles.
For the networks we consider in this paper, triangular subgraphs are either not present 
(in case of square, cubic, hypercubic lattices)
or negligible for large $N$ (e.g. for Erd\H{o}s-R{\'e}nyi random networks \cite{Note1}), so 
we only need to consider the system
\begin{equation*}
	\begin{blockarray}{ccccccc}
		&& [\motI] & [\motIS] & [\motII] & [\motSSI] & [\motISI] \\
		\begin{block}{cc[c|cc|cc]}\relax
			\dot{[\motI]} & & -\gamma & \beta & 0 &  &  \\\cline{3-7}
			\dot{[\motIS]} & &  &-\beta-\gamma & \gamma & \beta & -\beta\\
			\dot{[\motII]} & &  & 2\beta &-2\gamma & 0 & 2\beta\\	
		\end{block}
	\end{blockarray}\;.
\end{equation*}
The number of conservation relations used for elimination depends on whether the networks
have a homogeneous degree.

\paragraph{Networks with homogeneous degree} In this case,
the conservation relations are
\begin{align}
\left[\motI\right]+\left[\motS\right] = & N,\label{eq:siscons1}\\
\left[\motIS\right]+\left[\motSS\right] = & \kappa\left[\motS\right],\label{eq:siscons2}\\
\left[\motII\right]+\left[\motIS\right] = & \kappa\left[\motI\right],\label{eq:siscons3}
\end{align}
However, due to the dynamic irrelevance of $[\motS]$ and $[\motSS]$, only \eqref{eq:siscons3}
can be used to eliminate further variables. We use it to eliminate $\left[\motII\right]=\kappa[\motI]-[\motIS]$. Following (\ref{eq:subcons}-\ref{eq:xdotsubcons2}),
this means
\begin{align*}
	\mathbf{x}&=
	\begin{bmatrix}
		[\motI]&
		[\motIS]&
		[\motII]\\
	\end{bmatrix}^T,
	\\
	\tilde{\mathbf{x}}&=
	\begin{bmatrix}
		[\motI]\\
		[\motIS]\\
	\end{bmatrix},
	\;
	\tilde{\mathbf{x}}_3=
	\begin{bmatrix}
		[\motSSI]\\
		[\motISI]\\
	\end{bmatrix},
\end{align*}
and
\begin{align*}
	\mathbf{E}=
	\begin{bmatrix}
		1 & 0\\
		0 & 1\\
		\kappa & -1\\
	\end{bmatrix},
	\;
	\mathbf{c}=
	\begin{bmatrix}
		0\\
		0\\
		0\\
	\end{bmatrix},
	\;
	\tilde{\mathbf{Q}}_{1\cdots2,1\cdots2}=&
	\begin{bmatrix}
		-\gamma & \beta & 0\\
		 0 &-\beta-\gamma & \gamma\\
	\end{bmatrix},
	\\
	\tilde{\mathbf{Q}}_{23}=&
	\begin{bmatrix}
		\beta & -\beta\\
	\end{bmatrix},
\end{align*}
such that we can calculate $\mathbf{Q}',\tilde{\mathbf{c}}$ to obtain
\begin{equation}
	\begin{blockarray}{cccccc}
		&& [\motI] & [\motIS] & [\motSSI] & [\motISI]\\
		\begin{block}{cc[c|c|cc]}\relax
			\dot{[\motI]} & & -\gamma & \beta &  & \\\cline{3-6}
			\dot{[\motIS]} & & \gamma\kappa &-\beta-2\gamma & \beta & -\beta\\
		\end{block}
	\end{blockarray}\;.
\end{equation}
Applying
the closure \eqref{eq:mclosabchom} for degree-homogeneous networks [resulting from \eqref{eq:dec:mot:chord}],
\begin{align}
\left[\motSSI\right] & \approx \kappa(\kappa-1)N\frac{\frac{\left[\motSS\right]}{\kappa N}\frac{\left[\motIS\right]}{\kappa N}}{\frac{\left[\motS\right]}{N}}=\frac{\kappa-1}{\kappa}\frac{[\motSS][\motIS]}{[\motS]},\nonumber\\
\left[\motISI\right] & \approx \kappa(\kappa-1)N\frac{\frac{\left[\motIS\right]^2}{\kappa^2 N^2}}{\frac{\left[\motS\right]}{N}}=\frac{\kappa-1}{\kappa}\frac{[\motIS]^2}{[\motS]},\label{eq:mf2:hom:closure}
\end{align}
we obtain the final nonlinear system
\begin{equation}
	\begin{blockarray}{cccccc}
		&& [\motI] & [\motIS] & \frac{[\motSS][\motIS]}{[\motS]} & \frac{[\motIS]^2}{[\motS]}\\
		\begin{block}{cc[c|c|cc]}\relax
			\dot{[\motI]} & & -\gamma & \beta &  & \\\cline{3-6}
			\dot{[\motIS]} & & \gamma\kappa &-\beta-2\gamma & \beta\frac{\kappa-1}{\kappa} & -\beta\frac{\kappa-1}{\kappa}\\
		\end{block}
	\end{blockarray}\;,
\end{equation}
which, after elimination of $[\motS]$ and $[\motSS]$ via (\ref{eq:siscons1}-\ref{eq:siscons2})
and normalisation via \eqref{eq:cons:norm} becomes
\begin{equation}\label{eq:mf2:sis:1}
	\begin{blockarray}{ccccc}
		& \llbracket\motI\rrbracket & \llbracket\motIS\rrbracket & (1-\frac{\llbracket\motIS\rrbracket}{1-\llbracket\motI\rrbracket})\llbracket\motIS\rrbracket & \frac{\llbracket\motIS\rrbracket^2}{1-\llbracket\motI\rrbracket}\\
		\begin{block}{c[c|c|cc]}\relax
			\dot{\llbracket\motI\rrbracket} & -\gamma & \kappa\beta &  & \\\cline{2-5}
			\dot{\llbracket\motIS\rrbracket} & \gamma &-\beta-2\gamma & \beta(\kappa-1) & -\beta(\kappa-1)\\
		\end{block}
	\end{blockarray},
\end{equation}
From the steady state solutions \eqref{eq:SS2homkap}
we find that the epidemic threshold,
is now located at $\beta/\gamma{=}(\kappa-1)^{-1}$. We also derived non-trivial steady state correlations via (\ref{eq:corrs}) in \eqref{eq:sscorrs2homkap}.

\paragraph{Networks with heterogeneous degree} 
Here, \eqref{eq:siscons2}
and \eqref{eq:siscons3} do not hold, but the total frequency of any given subgraph is still conserved. Hence, the conservation relations are \eqref{eq:siscons1} from order $1$ and 
\begin{equation}\label{eq:alginvSIS2het}
\left[\motSS\right]+\left[\motII\right]+2\left[\motIS\right]=\kappa N.
\end{equation}
This means that we cannot eliminate $[\motII]$ here (unlike in case of networks 
with homogeneous degree), such that we have three instead of two equations.
Also applying the closure for ER networks \eqref{eq:mclosabcER} [resulting from \eqref{eq:dec:mot:chord}], 
  \begin{align}
\left[\motSSI\right]&\approx\frac{\left[\motSS\right]\left[\motIS\right]}{\left[\motS\right]},\qquad\left[\motISI\right]\approx\frac{\left[\motIS\right]^{2}}{\left[\motS\right]},\label{eq:mf2:het:closure}
  \end{align}
we obtain 
\begin{equation}\label{eq:mf2:sis:het}
	\begin{blockarray}{ccccccc}
		&& [\motI] & [\motIS] & [\motII] & \frac{[\motSS][\motIS]}{[\motS]} & \frac{[\motIS]^2}{[\motS]}\\
		\begin{block}{cc[c|cc|cc]}\relax
			\dot{[\motI]} & & -\gamma & \beta & 0 &  & \\\cline{3-7}
			\dot{[\motIS]} & &  &-\beta-\gamma & \gamma & \beta & -\beta\\
			\dot{[\motII]} & &  & 2\beta &-2\gamma & 0 & 2\beta\\	
		\end{block}
	\end{blockarray}\;.
\end{equation}
After substitution of $[\motS]$ and $[\motSS]$
via the conservation relations (\ref{eq:siscons1},\ref{eq:alginvSIS2het}) and normalisation via \eqref{eq:cons:norm}, this becomes
\begin{equation}\label{eq:mf2:sis:het2}
	\begin{blockarray}{cccccc}
		& \llbracket\motI\rrbracket & \llbracket\motIS\rrbracket & \llbracket\motII\rrbracket & \frac{(1-2\llbracket\motIS\rrbracket-\llbracket\motII\rrbracket)\llbracket\motIS\rrbracket}{1-\llbracket\motI\rrbracket} & \frac{\llbracket\motIS\rrbracket^2}{1-\llbracket\motI\rrbracket}\\
		\begin{block}{c[c|cc|cc]}\relax
			\dot{\llbracket\motI\rrbracket} & -\gamma & \kappa\beta & 0 &  & \\\cline{2-6}
			\dot{\llbracket\motIS\rrbracket} & &-\beta-\gamma & \gamma & \kappa\beta & -\kappa\beta\\
			\dot{\llbracket\motII\rrbracket} & & 2\beta &-2\gamma & 0 & 2\kappa\beta\\	
		\end{block}
	\end{blockarray}\;.
\end{equation}
There are three steady states, of which two are in the admissible range \eqref{eq:SS2hetkap}.
Here, the epidemic threshold is located at $\beta/\gamma{=}\kappa^{-1}$, as in the first-order
mean field model. The steady state correlations are given in \eqref{eq:sscorrs2hetkap}.

\subsubsection{MF3-MF5} \label{sec:mf35}
Approximations of higher order than two correspond to cluster variation
approximations in statistical physics \cite{kikuchi1951}.
MF3 results from neglecting dependence 
beyond distance 2 and is commonly known as the `triple approximation' \cite{House2009}. 
Its step-by-step derivation for the square lattice is shown in Appendix \ref{sec:MF3}
and results in four equations (after substitution with conservation relations). 
In Appendix \ref{sec:MFhoeqs}, we show the derivation of the unclosed
MF4 for a general network as \texttt{Mathematica} notebook output. We
also derived the closed MF4 and MF5 for the square lattice. The number of equations after
elimination with conservation relations is respectively 14 and 37. The derivation
of MF4 with closure is shown in Appendix \ref{sec:MFhoeqs} point 2.
The steady states of
MF4 and MF5 are shown in Figure~\ref{fig:Comparison-of-the} of Section \ref{sec:Comparison}. 

\subsection{Comparison to simulations}\label{sec:Comparison}
Here, we compare the steady states of MF1-5 of
SIS epidemic spreading with those of the simulations on a selection of
network types: lattices, regular random networks and Erd\H{o}s-R{\'e}nyi random
networks. 

\begin{figure*}
\centering
\includesvg[inkscapelatex=false,width=2\columnwidth]{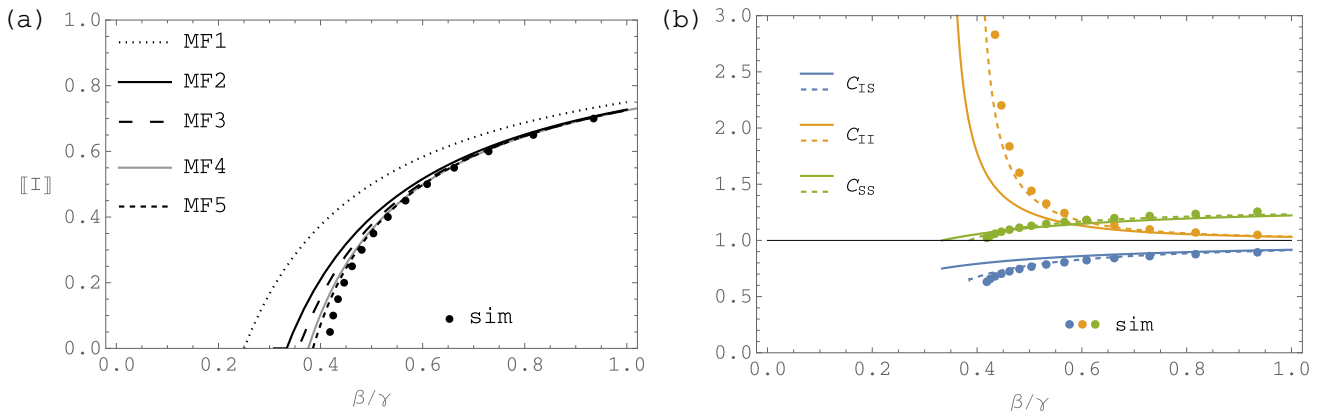}
\caption{Comparison of the mean-field approximations of order $1$ to $5$ (lines) with simulations (markers)
of SIS epidemic spreading on a square lattice as a function of $\beta/\gamma$:
(a) Nontrivial steady states, (b) Correlations at distance $1$ [shown for
MF$2$ (solid), MF$5$ (dashed), and simulations (dots)]. \label{fig:Comparison-of-the}}
\end{figure*}

\paragraph{Square lattice} In Figure \ref{fig:Comparison-of-the}a,
we show, for the square lattice, the steady state fraction of infecteds versus
$\beta/\gamma$ for MF1-MF5 compared to simulations. The steady states of the mean-field
models get closer to those of the simulation with increasing order. All mean-field
models have an increasing bias in their non-trivial (endemic) steady states when approaching the critical value of $\beta/\gamma$ from above. Figure \ref{fig:Comparison-of-the}b compares
the steady state distance-1 correlations between species types from MF2 and MF5
to those in simulations. Species of the same type cluster whereas
different species tend to avoid each other (compared to a random distribution).
The infected-infected correlation diverges when approaching the critical value of $\beta/\gamma$ from above and has a singularity
at the bifurcation. E.g. for MF$2$, via (\ref{eq:sscorrs2homkap}), we have $C_{II}^{*}(\beta/\gamma)\to\infty$ for $\beta/\gamma\, \searrow\, (\kappa-1)^{-1}$
(limit from above in the endemic equilibrium). We recall that the errors in the mean-field models visible in Figure~\ref{fig:Comparison-of-the}a can be due to violation of the statistical dependence assumption beyond distance $k$, non-chordality of 
the independence map, and violation of the spatial homogeneity assumption.
Spatial homogeneity can be violated in two ways: via heterogeneity of structure
and via heterogeneity of dynamics \cite{Sharkey2008}. As the structure
of a lattice is homogeneous and we only study steady states (i.e.
there are no dynamics) that are spatially homogeneous (see Figure \ref{fig:ctrl_t_bif}b), spatial inhomogeneity can not be a source of bias here. 
This leaves statistical dependence and non-chordality
as only sources of bias. In the square lattice there are always cycles with diameter
larger than $\mathsf{d}$ for any chosen $\mathsf{d}$ along which unaccounted for information can spread, unless $\mathsf{d}$ is greater than or equal to the graph diameter of the entire lattice. This means that 
mean-field models that do not consider motifs of size up to the network
diameter minus one cannot be exact. 

As MF1 
assumes no correlation between
the states of neighbouring nodes, the distance between the horizontal line 
through $1$ and the $\boldsymbol{+}$ markers of the correlations in the simulations
is a measure of the bias of MF1 due to neglection of correlations
in MF1 closures. Likewise, the distance between the steady state MF2/MF5
correlations and the simulations is due to neglection of (higher-order
and conditional) dependence in MF2/MF5 closures and higher-level non-chordality.
All models have larger biases closer to the critical point, where higher-order correlations become more important. This is a well-known characteristic of continuous phase transitions, in
which correlations occur on increasingly long ranges when approaching
the phase transition (bifurcation). Figure~\ref{fig:Steady-state-correlation} in Appendix \ref{sec:dcorrs} shows the correlation as a function of distance from a central point
[via \eqref{eq:CabD}], for various values of $\beta/\gamma$. It shows,
as expected from phase transitions theory, that correlations at any
distance are larger closer to the critical point. At the critical
point, theory shows that correlations occur at all distances \cite{Henkel2008}.

\begin{figure*}
\includesvg[inkscapelatex=false,width=2\columnwidth]{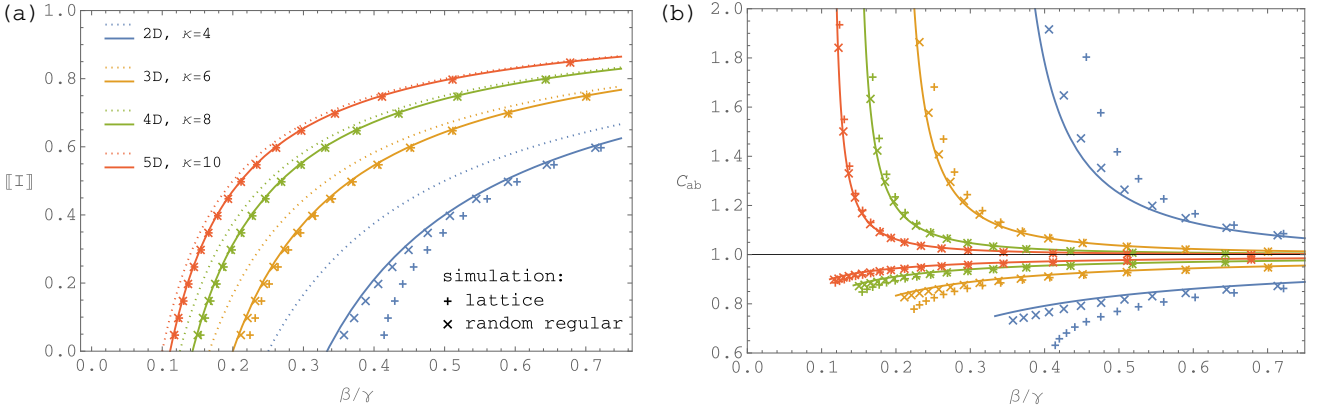}

\caption{Comparison of the mean-field approximations (dotted lines: MF1, solid
lines: MF2) and simulations (markers) of SIS epidemic
spreading on lattices and random regular networks with number of neighbours
4,6,8,10 as a function of $\beta/\gamma$: (a) Nontrivial steady states,
(b) Correlations: $C_{II}^{*}$ (above 1), $C_{SI}^{*}$(below 1).\label{fig:Comparison-of-the-1}}
\end{figure*}

Comparing the square lattice to the 4-neighbour random regular graph, we can see that
there is also a phase transition, but there is substantially less bias than in the square
lattice (Figure \ref{fig:Comparison-of-the-1} blue $\boldsymbol{\times}$ vs $\boldsymbol{+}$). This is because random regular graphs are
locally treelike, and hence, unlike in the square lattice, correlations over longer
distances can be captured well via a decomposition of larger motifs into links (MF2).  
The small remaining bias in the 4-regular random graph we suspect to be 
because the assumed conditional independence is not valid for all states
\cite{Sharkey2015,Sharkey2015a}.

\paragraph{General cubic lattices and random regular networks} We show in Figure \ref{fig:Comparison-of-the-1} how the steady
states and correlations in MF1, MF2 and simulations depend on the number
of neighbours in $\mathsf{d}$-dimensional cubic lattices and random regular networks. When 
$\mathsf{d}$ is
the lattice dimension, the lattice degree is $\kappa=2\mathsf{d}$. The observations
of the square lattice generalise to cubic/hypercubic lattices and
random regular networks (at least up to $\kappa=10$): \emph{i}. there is
a transcritical bifurcation at a particular value of \textbf{$\beta/\gamma$},
where $C_{II}^{*}$ becomes singular, \emph{ii}. MF1 and MF2 capture qualitatively
the steady state fraction of infecteds and correlations are captured
qualitatively by MF2, \emph{iii}. MF2 is less biased than MF1, \emph{vi}. the
bias is larger closer to the bifurcation. According to MF1, the bifurcation
occurs at $\kappa^{-1}=(2\mathsf{d})^{-1}$ and according to MF2 at 
$(2\mathsf{d}-1)^{-1}=(\kappa-1)^{-1}$
(see Sections \ref{sec:First-order} and \ref{sec:Second-order}).
\citet{Liggett2005} proved that for lattices, the critical value
predicted by MF2 is a lower bound. Figure \ref{fig:Comparison-of-the-1}a
shows that MF2 and this lower bound is approached increasingly closely
when the lattice dimension increases. Due to higher clustering of
neighbours, the epidemic threshold in lattices is higher than that
in the corresponding random regular network \cite{Keeling1999},
but this difference decreases with dimension/degree (Figure \ref{fig:Comparison-of-the-1}a).
The steady states of a 5-dimensional hypercubic lattice and of a random
regular network with degree 10 are indistinguishable from each other
and from MF2. This is because random walks in space of dimension 5 or higher
have a finite number of intersections almost surely
\cite{Erdos1960,Heydenreich2017a,Lawler2010}. If 
the path along which an infection travels is seen
as a random walk, having many intersections in $\mathsf{d}\leq4$ means that
one cannot ignore alternative infection paths. In $\mathsf{d}>4$, it is harder
for infections to travel via alternative paths to the same point,
and hence those paths resemble trees more closely. Hence, as explained above 
and in \citet{Sharkey2015a}, MF2 should then be more accurate. 

\begin{figure*}
\includesvg[inkscapelatex=false,width=2\columnwidth]{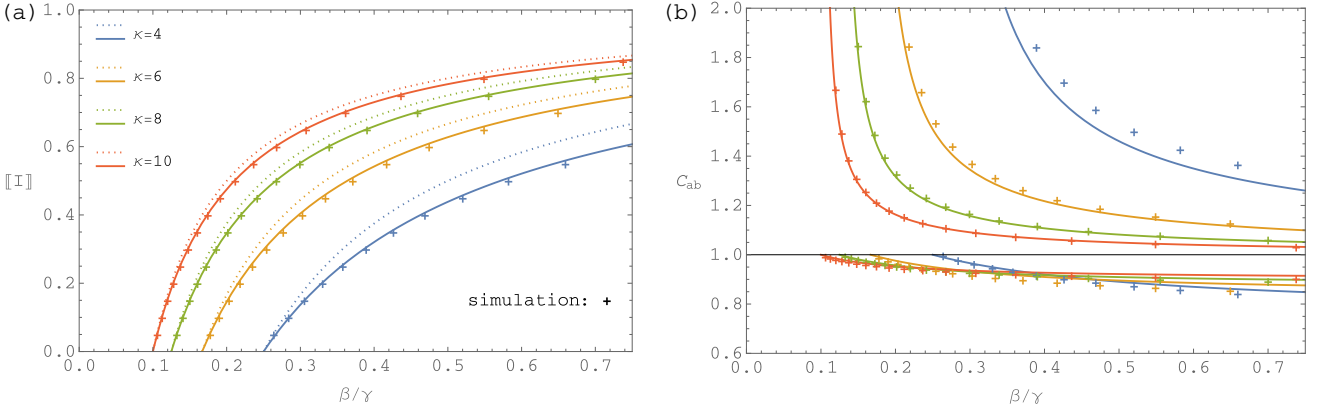}

\caption{Comparison of the mean-field approximations (dotted lines: MF1, solid
lines: MF2) and simulations (markers) of SIS epidemic spreading
on an Erd\H{o}s-R{\'e}nyi random network with number of neighbours 4,6,8,10
as a function of $\beta/\gamma$: (a) Nontrivial steady states, (b)
Correlations: $C_{II}^{*}$ (above 1), $C_{SI}^{*}$(below 1).\label{fig:Comparison-of-the-2}}

\end{figure*}

\paragraph{Erd\H{o}s-R{\'e}nyi random networks} Finally, we show in Figure \ref{fig:Comparison-of-the-2} how the
steady states and correlations in MF1, MF2 and simulations depend on
the number of neighbours in an Erd\H{o}s-R{\'e}nyi random network. Recall
that MF2 on an Erd\H{o}s-R{\'e}nyi random network is different from that of the networks
above because it has one fewer conservation relation. As above, the
behaviour of the steady state solutions is captured qualitatively
by MF1 and MF2 and of the correlations by MF2 alone. Also as above,
networks with a larger degree have lower biases and MF2 is better
than MF1, but now there seems to be a slight increase of bias with $\beta/\gamma$,
at least in the range inspected. Despite spatial heterogeneity
of the degree in Erd\H{o}s-R{\'e}nyi random networks, there is considerably less bias
than in lattices, confirming that the presence of cycles beyond closure distance is the dominant cause for mean-field model biases in the steady states of SIS spreading.

\section{Summary and conclusions}

Previous work found that exact closed individual-level 
moment equations exist for SIR spreading on arbitrary networks,
with the requirement to consider larger-sized motifs, and therefore more equations,
for networks that are decreasingly tree-like \cite{Sharkey2015a}.
While for other dynamics than SIR spreading it may not be possible to
prove exactness for a finite number of closed moment equations,
it is generally found that accuracy 
increases with the order of approximation \cite[e.g.][]{House2009,Kiss2017} (which was confirmed here). Feasibility of automated derivation of exact closed moment equations
for SIR epidemic spreading was shown by \cite{Sharkey2015a}, while an automated procedure
to derive unclosed moment equations for arbitrary dynamics was developed by \cite{Danos2020}.
We developed an automated procedure to both derive and close population-level
moment equations for arbitrary dynamics on networks at any approximation order, allowing us to consider mean-field models of higher orders than typically derived by hand. 
For this purpose, we developed a method to derive closure schemes from predefined
independence assumptions. Our closure formulas rely, besides the requirements of spatial homogeneity and
large network size, on the assumption of conditional independence beyond
distance $k$ and $k$-chordality of the considered network. 
Consistently, our simulations of SIS epidemic spreading showed
that, at given approximation order, the largest biases occurred for networks with many short cycles of any size, such as lattices, in parameter regimes with long-range correlations, such as near continuous phase transitions. Note however that, for lattices, we found the bias of mean-field models to decrease with lattice dimension, which is consistent with results in percolation and phase transitions theory,
where it was shown that the importance of cycles decreases with the lattice dimension
\cite{Erdos1960,Heydenreich2017a,Lawler2010}. We also showed
that the conventional procedure of truncation at a maximum motif size instead 
of at a maximum motif diameter necessitates independence assumptions 
that are inconsistent for different approximated motifs or that may be incompatible with the network (see \ref{pt:exactness}). This suggests that choosing a moment space
that consists of motifs at increasing diameter instead of increasing size
would lead to more accurate mean-field models. 

Whereas our method still needs to be tested more widely, we expect it to lend 
itself well to study dynamics on networks
with density and size of short cycles between that of random networks
and (low-dimensional) lattices, particularly when
the network is structurally homogeneous. In these cases,
derivation by hand may be too tedious while the final
set of moment equations is still more manageable than
the Markov chain simulations.
For networks with considerable
degree heterogeneity 
and/or community structure, we expect
approximate master equation methods 
\cite{marceau2010,gleeson2013,fennell2019,stonge2021,Cui2022} to be more efficient.
Our 
approach focused on static networks with at most
nearest-neighbour interactions, but it can be extended to adaptive 
networks such as those studied in 
\cite{Gross2006,demirel2014,Danos2020}, and to dynamics
with higher-order interactions \cite{battiston2020} -- requiring
reaction rate tensors $\mathbf{R}^p$ with $p\ge2$. Our approximation
scheme of Section \ref{sec:Closure-scheme} can be applied more generally,
to understand precisely which independence assumptions are taken in 
other existing types of mean-field 
approximations than moment closure, or to devise new mean-field
approximations. 
It may also serve to extend the use of message-passing methods 
\cite{karrer2010,wilkinson2014,koher2019} 
for epidemic modelling to graphs with cycles. 

While we were able to derive
closure formulas by assuming statistical independence, 
leading to mean-field models, other assumptions can be used
to obtain closures \cite{kuehn2016},
such as maximum entropy \cite{Rogers2011}, or time scale separation 
between moments at different orders \cite{GrossKevrekidis2008}. 
Sometimes one can find an appropriate moment space and closure by taking account of the characteristic 
features of the process in consideration, leading to a description with greater efficiency compared to what is obtainable by using the size-based moment space and closing via independence assumptions \cite{Wuyts2022}. 
We expect that for finding good moment spaces and closures, equation-free and machine-learning 
methods \cite{Kevrekidis2009,Patsatzis2022} will
play an important role, in particular because the most appropriate low-dimensional descriptions or their closures may not necessarily be available in closed form 
\cite{Rogers2011,GrossKevrekidis2008,Kevrekidis2009,Patsatzis2022}. It is subject of future work to explore if and how these different approaches relate.

\begin{acknowledgments}
This work was supported by the UK Engineering and Physical Sciences Research Council (EPSRC) grants EP/N023544/1 and EP/V04687X/1.
\end{acknowledgments}

\bibliographystyle{apsrev4-1}
\bibliography{./library}

\onecolumngrid
\renewcommand{\thesection}{A-\Roman{section}}
\renewcommand{\thesubsection}{\Alph{subsection}}
\renewcommand{\theparagraph}{\alph{paragraph}}
\renewcommand{\theequation}{A\arabic{equation}}
\renewcommand{\thefigure}{A\arabic{figure}}
\renewcommand{\thetable}{A\arabic{table}}
\setcounter{equation}{0}  
\setcounter{section}{0}  
\setcounter{figure}{0}  
\newpage
\begin{center}
	\fontsize{20}{24}\bfseries\scshape{APPENDIX}
\end{center}

\paragraph*{Notation and references}
\label{notref} References to items in this document are preceded
by `A'. Items not preceded by `A' refer to the main text.

\tableofcontents

\newpage
\section{Derivation of the moment equations}\label{sec:Differential-and}
In this section, we derive
expressions for the expected rate of change and conservation relations 
of motif counts, first shown in (\ref{eq:cons:norm},~\ref{eq:cons2}) of Section
\ref{sec:alg:iv}. Both can be seen as invariance
relations, the former being of a differential type, and derived from the master equation for the Markov chain on the network, and the latter of an
algebraic type, following from the property that the network is fixed. 

\paragraph{Master equation for transitions in a Markov chain} We
recall that the state of the system for our network with $N$ nodes is
given by (using $:=$ for ``defined as'')
$
X:=(X_{1},X_{2},...,X_{N})\in\{1,\ldots,n\}^N,
$
where $X_{i}$ is the label of the species that occupies node $i$. 
Hence, the total number of states is $n^{N}$.
The probabilistic transition from one state to another, following a discrete-state continuous-time Markov chain, defines an evolution equation for the probability of being in each of these states, the so-called
master equation (also known as the Kolmogorov-forward
equation for a Markov jump process \cite{Gardiner2009}). The probability density $P(X,t)$ for a particular
state $X$ at time $t$ changes according to
\begin{eqnarray}
	\ddt P(X,t) & = & \sum_{X'\neq X}\left[w(X'\rightarrow X)P(X',t)-w(X\rightarrow X')P(X,t)\right],\label{eq:ME}
\end{eqnarray}
where $w$ denotes the transition rate between system states. If we
define $\mathbf{W}$ as a $n^{N}\times n^{N}$ transition rate matrix
with non-diagonal entries $w(X'\rightarrow X)$ and diagonal entries
$-\sum_{X'\neq X}w(X\rightarrow X')$, we can rewrite the master equation
as 
\[
\dot{\mathbf{\mathbf{P}}}(t)=\mathbf{W}\mathbf{P}(t),
\]
which describes the evolution of the density for all states as elements of a vector $\mathbf{P}$
(and not just of one particular state as in \eqref{eq:ME}). Because
almost surely at most one node can change state at any one time $t$, only states differing from each other in one node can be directly
transitioned between. Therefore, $\mathbf{W}$ must be sparse, having
only $N(n-1)$ entries in each row or column, as each of the $N$
nodes can convert to any of the $n-1$ other species. This permits
writing \eqref{eq:ME} as
\begin{align}
	\ddt P(X,t)&=\sum_{i=1}^{N}\sum_{k\neq X_{i}}^{n}\left[w_{i}(X_{i\shortto k}\rightarrow X)P(X_{i\shortto k},t)-w_{i}(X\rightarrow X_{i\shortto k})P(X,t)\right]\mbox{,}\label{eq:ME2}
\end{align}
where
\[
X\to X_{i\shortto k}: (X_{1},...,X_{i},...,X_{N})\mapsto(X_{1},...,k,...,X_{N})
\]
is the operator that
replaces the species at the $i$th node by species $k$, and $w_{i}(.)$ 
is the conversion rate at node $i$.
\paragraph{Motifs and their frequencies}
In what follows, we will derive from the master equation the evolution
of the frequency (or total count) of the \emph{motifs} $\boldsymbol{x}^\mathsf{a}$, as defined in Section \ref{sec:ds-ct-mc}. Recall that these motifs, 
are defined by their state label vector $\boldsymbol{x}$ of size $m$ and the 
adjacency structure between motif nodes $\mathsf{a}$. We denote a single \emph{occurrence} at a given $\boldsymbol{i}\in S(m,N)$ as
\begin{align}
	\left[\boldsymbol{x}_{\boldsymbol{i}}^{\mathsf{a}}\right](t) & := \delta_\mathsf{a}(\mathsf{A}_{\boldsymbol{i}})\delta_{\boldsymbol{x}}(X_{\boldsymbol{i}}(t))
	\mbox{,}\label{eq:motif1}
\end{align}
which requires an
exact match of the adjacency $\mathsf{a}$ by
$\mathsf{A}_{\boldsymbol{i}}$ and state label vector $\boldsymbol{x}$ by $X_{\boldsymbol{i}}$. 
For example, if $\mathsf{A}$ contains a
connected triangle between nodes $1$, $2$ and $3$, then motifs with
$\mathsf{a}=\{(1,2),(2,3)\}$ would not occur on
$\boldsymbol{i}=(1,2,3)$. The total count of motifs in the large network is then the number of such exact matches, which is obtained by summing over all indices $\boldsymbol{i}\in S(m,N)$, as shown in \eqref{eq:motifcount:all}.
Since the large network structure is constant in time, we can use the counts
of induced subgraphs $\left[\mathsf{a}\right]$ given by \eqref{eq:cons:norm} as a normalisation factor for the (variable) counts of motifs such that we may consider normalised motif frequencies
\begin{equation}
	\llbracket\boldsymbol{x}^{\mathsf{a}}\rrbracket:=\left[\boldsymbol{x}^{\mathsf{a}}\right]/\left[\mathsf{a}\right]\mbox{.}\label{eq:norm}
\end{equation}
\paragraph{Evolution of expected counts} Total and normalised counts
$\left[\boldsymbol{x}^{\mathsf{a}}\right]$ and
$\llbracket\boldsymbol{x}^{\mathsf{a}}\rrbracket$ 
refer to realisations of states $X$ on the large random network such
that they are random variables. Similarly,
$\left[\boldsymbol{x}_{\boldsymbol{i}}^{\mathsf{a}}\right](X)$ is a
random variable in $\{0,1\}$ for each index vector $\boldsymbol{i}\in S(m,N)$
once we take the randomness of states $X$ into account. Its
expectation is
\begin{equation}\label{eq:motif1exp}
	\langle\left[\boldsymbol{x}_{\boldsymbol{i}}^{\mathsf{a}}\right]\rangle=\sum_{X}\left[\boldsymbol{x}_{\boldsymbol{i}}^{\mathsf{a}}\right](X)P(X,t)\mbox{.}
\end{equation}
The master equation \eqref{eq:ME2} for the density $P$ implies that the expectation satisfies the differential equation
\begin{eqnarray}
	\ddt \langle\left[\boldsymbol{x}_{\boldsymbol{i}}^{\mathsf{a}}\right]\rangle & = & \sum_{X}\left[\boldsymbol{x}_{\boldsymbol{i}}^{\mathsf{a}}\right](X)\ddt P(X,t)\nonumber \\
	& = & \sum_{X}\Bigl\{\delta_\mathsf{a}(\mathsf{A}_{\boldsymbol{i}})\delta_{\boldsymbol{x}}(X_{\boldsymbol{i}})\sum_{i'=1}^{N}\sum_{k\neq X_{i'}}^{n}\left[w_{i'}(X_{i'\shortto k}\rightarrow X)P(X_{i'\shortto k},t)-w_{i'}(X\rightarrow X_{i'\shortto k})P(X,t)\right]\Bigr\},\nonumber \\
	& = & \Bigl\langle\delta_\mathsf{a}(\mathsf{A}_{\boldsymbol{i}})\sum_{i'=1}^{N}\sum_{k\neq X_{i'}}^{n}\Bigl[\delta_{\boldsymbol{x}}\bigl((X_{i'\shortto k})_{\boldsymbol{i}}\bigr)-\delta_{\boldsymbol{x}}\bigl(X_{\boldsymbol{i}}\bigr)\Bigr]w_{i'}(X\rightarrow X_{i'\shortto k})\Bigr\rangle,\nonumber\\
	& = & \Bigl\langle\delta_\mathsf{a}(\mathsf{A}_{\boldsymbol{i}})\sum_{p=1}^{m}\sum_{k\neq X_{i_{p}}}^{n}\Bigl[\delta_{\boldsymbol{x}}\bigl((X_{i_p\shortto k})_{\boldsymbol{i}}\bigr)-\delta_{\boldsymbol{x}}\bigl(X_{\boldsymbol{i}}\bigr)\Bigr]w_{i_{p}}(X\rightarrow X_{i_p\shortto k})\Bigr\rangle,\nonumber \\
	& = & \Bigl\langle\delta_\mathsf{a}(\mathsf{A}_{\boldsymbol{i}})\sum_{p=1}^{m}\sum_{k\neq X_{i_{p}}}^{n}\delta_{\boldsymbol{x}_{p\shortto\emptyset}}(X_{\boldsymbol{i}_{p\shortto\emptyset}})\Bigl[\delta_{x_p}(k)-\delta_{x_p}(X_{i_{p}})\Bigr]w_{i_{p}}(X\rightarrow X_{i_p\shortto k})\Bigr\rangle,\label{eq:ddtmMi}
\end{eqnarray}
where in the second step, we substituted
$\delta_{\boldsymbol{x}}(X_{\boldsymbol{i}})w_{i'}(X_{i'\shortto k}\rightarrow X)P(X_{i'\shortto k},t)$ for
$\delta_{\boldsymbol{x}}\bigl((X_{i'\shortto k})_{\boldsymbol{i}}\bigr)w_{i'}(X\rightarrow X_{i'\shortto k})P(X,t)$
which corresponds to a reordering of the terms in $\sum_{X}$, and
where we have used the notation $\left\langle \cdot\right\rangle $ for
the expectation $\sum_{X}(\cdot)P(X,t)$.  In the third step, we 
used that only changes to the nodes belonging to
$\boldsymbol{i}$ matter. In the last step, we factored out the
common elements in the delta functions corresponding to all but the
$p$th element of $X_{\boldsymbol{i}}$ and $\boldsymbol{x}$, and we
used the subscript notation $(\cdot)_{p\shortto\emptyset}$ to denote
a vector with element $p$ removed. The expression on the
right-hand side in \eqref{eq:ddtmMi} can be understood independent of
the prior algebraic manipulations: the expected change rate of
$\bigl\langle\left[\boldsymbol{x}_{\boldsymbol{i}}^{\mathsf{a}}\right]\bigr\rangle$
equals the sum of expected rates for each of the nodes $i_p$ in
$\boldsymbol{i}$ changing its state \emph{to} $x_p$ minus the rate of
node $i_p$ changing its state \emph{from} $x_p$.

\paragraph{Conversion rates}
Next, we will insert the two types of admissible transitions 
(as discussed in Section \ref{sec:ds-ct-mc}) into
\eqref{eq:ddtmMi}, namely spontaneous conversions with rates given in
matrix $\mathbf{R}^0\in\mathbb{R}^{n\times n}$, and, 
conversions due to interactions with a
single nearest neighbour with rates given in
$\mathbf{R}^1\in\mathbb{R}^{n\times n\times
	n}$. 
The diagonal entries of $\mathbf{R}^0,\mathbf{R}^1$ are zero without loss of generality. 
For these transition types, the rates $w_{i_p}$ in \eqref{eq:ddtmMi} are
\begin{eqnarray}
	w_{i_{p}}(X\rightarrow X_{i_p\shortto k}) & = & \sum_{a,c}^{n}R^1_{akc}\delta_a(X_{i_{p}})\sum_{j}^{N}\mathsf{A}_{i_{p}j}\delta_c(X_{j})+\sum_{a}^{n}R^0_{ak}\delta_a(X_{i_{p}}).\label{eq:wieta}
\end{eqnarray}
With the rates in \eqref{eq:wieta}, and noting that 
\begin{align*}
	\delta_{x_p}(k)-\delta_{x_p}(X_{i_{p}})=
	\begin{cases}
		\phantom{-}1 & \mbox{if $x_{p}=k$,}\\
		-1 & \mbox{if $x_{p}=X_{i_{p}}$,}\\
		\phantom{-}0 & \mbox{if $k\neq x_{p} \mbox{\ and\ } x_{p}\neq X_{i_{p}}$,}  
	\end{cases}
\end{align*}
the sum inside the averaging brackets in \eqref{eq:ddtmMi} has the form
\begin{multline*}
	\Bigl[\sum_{a,c}^{n}R^1_{ax_{p}c}\delta_a(X_{i_{p}})\sum_{j}^{N}\mathsf{A}_{i_{p}j}\delta_c(X_{j})+\sum_{a}^{n}R^0_{ax_{p}}\delta_a(X_{i_{p}})\Bigr]\delta_{\boldsymbol{x}_{p\shortto\emptyset}}(X_{\boldsymbol{i}\setminus i_{p}})\\
	-\sum_{k\neq X_{i_{p}}}^{n}\Bigl[\sum_{c}^{n}R^1_{x_{p}kc}\sum_{j}^{N}\mathsf{A}_{i_{p}j}\delta_c(X_{j})+R^0_{x_{p}k}\Bigr]\delta_{\boldsymbol{x}}\bigl(X_{\boldsymbol{i}}\bigr),
\end{multline*}
where we used for the last two terms that $\delta_{x_p}(X_{i_{p}})\delta_{\boldsymbol{x}_{p\shortto\emptyset}}(X_{\boldsymbol{i}\setminus i_{p}})$ can be combined to $\delta_{\boldsymbol{x}}\bigl(X_{\boldsymbol{i}}\bigr)$.

\begin{figure}
	\centering
	\includegraphics[scale=0.45]{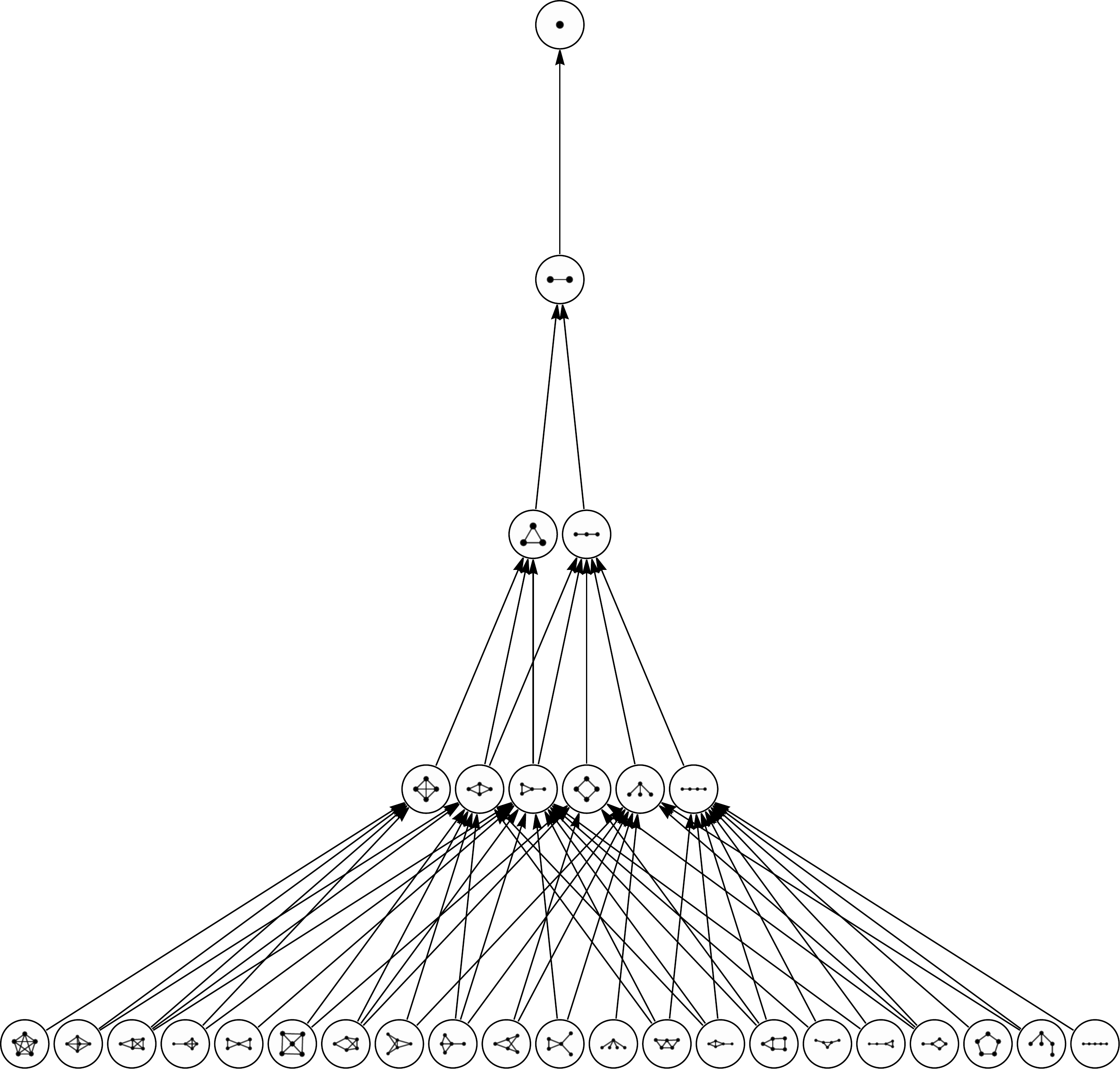}
	\caption{Dependence of moment equations of a given motif on motifs of one order higher due to nearest-neighbour interactions, up to order 4 (ignoring labels). \label{fig:Dependence-of-motif}}
\end{figure}

\paragraph{Differential equations for motif counts} When distributing the products and exploiting the linearity of the
averaging brackets, the differential equation \eqref{eq:ddtmMi} for the expected rate of change at $\boldsymbol{i}$ becomes
\begin{align*}
	\ddt \langle\left[\boldsymbol{x}_{\boldsymbol{i}}^{\mathsf{a}}\right]\rangle  = & \delta_\mathsf{a}(\mathsf{A}_{\boldsymbol{i}})\sum_{p}^{m}\Bigl[\sum_{a,c}^{n}R^1_{ax_{p}c}
	\Bigl\langle\delta_{\boldsymbol{x}_{p\shortto\emptyset}}(X_{\boldsymbol{i}_{p\shortto\emptyset}})\delta_a(X_{i_{p}})\sum_{j}^{N}\mathsf{A}_{i_{p}j}\delta_c(X_{j})\Bigr\rangle+\sum_{a}^{n}R^0_{ax_{p}}\Bigl\langle\delta_{\boldsymbol{x}_{p\shortto\emptyset}}(X_{\boldsymbol{i}_{p\shortto\emptyset}})\delta_a(X_{i_{p}})\Bigr\rangle\\
	& -\sum_{k\neq X_{i_{p}},c}^{n}R^1_{x_{p}kc}\Bigl\langle\delta_{\boldsymbol{x}}\bigl(X_{\boldsymbol{i}}\bigr)\sum_{j}^{N}\mathsf{A}_{i_{p}j}\delta_c(X_{j})\Bigr\rangle-\sum_{k\neq X_{i_{p}}}^{n}R^0_{x_{p}k}\Bigl\langle\delta_{\boldsymbol{x}}\bigl(X_{\boldsymbol{i}}\bigr)\Bigr\rangle\Bigl].
\end{align*}
After replacing index label $a$ by $k$, using the definition of $[\boldsymbol{x}^\mathsf{a}_{\boldsymbol{i}}]$ in \eqref{eq:motif1}
and using $(\boldsymbol{x}_{\boldsymbol{i}}^{\mathsf{a}})_{p\shortto k}$
to indicate $\boldsymbol{x}_{\boldsymbol{i}}^{\mathsf{a}}$ with its
$p$th element replaced by species $k$, this becomes
\begin{eqnarray*}
	\ddt \langle\left[\boldsymbol{x}_{\boldsymbol{i}}^{\mathsf{a}}\right]\rangle & = & \sum_{p}^{m}\sum_{k}^{n}\left[\sum_{c}^{n}R^1_{kx_{p}c}\Bigl\langle\left[(\boldsymbol{x}_{\boldsymbol{i}}^{\mathsf{a}})_{p\shortto k}\right]\sum_{j}^{N}\mathsf{A}_{i_{p}j}\delta_c(X_{j})\Bigr\rangle+R^0_{kx_{p}}\bigl\langle\left[(\boldsymbol{x}_{\boldsymbol{i}}^{\mathsf{a}})_{p\shortto k}\right]\bigr\rangle\right.\\
	&  & \left.-\sum_{c}^{n}R^1_{x_{p}kc}\Bigl\langle\left[\boldsymbol{x}_{\boldsymbol{i}}^{\mathsf{a}}\right]\sum_{j}^{N}\mathsf{A}_{i_{p}j}\delta_c(X_{j})\Bigr\rangle-R^0_{x_{p}k}\Bigl\langle\left[\boldsymbol{x}_{\boldsymbol{i}}^{\mathsf{a}}\right]\Bigr\rangle\right].
\end{eqnarray*}
This shows that $\bigl\langle\left[\boldsymbol{x}_{\boldsymbol{i}}^{\mathsf{a}}\right]\bigr\rangle$
can increase by a conversion from motifs that differ from $\boldsymbol{x}_{\boldsymbol{i}}^{\mathsf{a}}$
in only one node (first two terms), or decrease by having any of the
nodes in $\boldsymbol{x}_{\boldsymbol{i}}^{\mathsf{a}}$ convert to another species
(last two terms), with both increase and decrease possible via interaction
with neighbours and via spontaneous conversion. The factors of form
$\left[\cdot\right]\sum_{j}^{N}\mathsf{A}_{i_{p}j}\delta_c(X_{j})$ count the number of 
$c$-connections of motif $[\cdot]$ (located at $\boldsymbol{i}$) at node $i_p$. 
The neighbouring node $j$ with species $c$ can be part of the motif $[\cdot]$, or it can be outside of $[\cdot]$, in which case it gives rise to higher-order motifs. Therefore, by splitting
the neighbourhood sums as follows,
\begin{eqnarray*}
	\sum_{j}^{N}\mathsf{A}_{i_{p}j}\delta_c(X_{j}) & = & \sum_{j\in\boldsymbol{i}}^{N}\mathsf{A}_{i_{p}j}\delta_c(X_{j})+\sum_{j\notin\boldsymbol{i}}^{N}\mathsf{A}_{i_{p}j}\delta_c(X_{j}),
\end{eqnarray*}
their products with $\left[\boldsymbol{x}_{\boldsymbol{i}}^{\mathsf{a}}\right]$ and $\left[(\boldsymbol{x}_{\boldsymbol{i}}^{\mathsf{a}})_{p\shortto k}\right]$ count contributions of $c-$neighbours from within
versus from outside the motif separately: 
\begin{eqnarray}
	\left[(\boldsymbol{x}_{\boldsymbol{i}}^{\mathsf{a}})_{p\shortto k}\right]\sum_{j}^{N}\mathsf{A}_{i_{p}j}\delta_c(X_{j}) & = & \boldsymbol{\delta}_c(\boldsymbol{x})\mathsf{a}\boldsymbol{e}_p\left[(\boldsymbol{x}_{\boldsymbol{i}}^{\mathsf{a}})_{p\shortto k}\right]+\sum_{\boldsymbol{y}^\mathsf{b}\in{\cal N}^c_p\left((\boldsymbol{x}^{\mathsf{a})_{p\shortto k}}\right)}\sum_{j\notin\boldsymbol{i}}^{N}\left[\boldsymbol{y}_{\boldsymbol{i}, j}^{\mathsf{b}}\right],\\
	\left[\boldsymbol{x}_{\boldsymbol{i}}^{\mathsf{a}}\right]\sum_{j}^{N}\mathsf{A}_{i_{p}j}\delta_c(X_{j}) & = & \boldsymbol{\delta}_c(\boldsymbol{x})\mathsf{a}\boldsymbol{e}_p\left[\boldsymbol{x}_{\boldsymbol{i}}^{\mathsf{a}}\right]+\sum_{\boldsymbol{y}^\mathsf{b}\in{\cal N}^c_p\left({\boldsymbol{x}^{\mathsf{a}}}\right)}\sum_{j\notin\boldsymbol{i}}^{N}\left[\boldsymbol{y}_{\boldsymbol{i}, j}^{\mathsf{b}}\right]\mbox{.}\label{eq:split}
\end{eqnarray}
On the right-hand side, $\boldsymbol{e}_p$ is an $m$-dimensional vector with a $1$ at position $p$ and zeroes elsewhere, and $\boldsymbol{\delta}_c(\boldsymbol{x})$ a vector Kronecker delta function that returns a vector of the size of $\boldsymbol{x}$ with ones where the elements equal $c$ and zeroes elsewhere. 
Thus, the term $\boldsymbol{\delta}_c(\boldsymbol{x})\mathsf{a}\boldsymbol{e}_p$ counts the number of $c$-connections at position $p$ in the motif $\boldsymbol{x}^\mathsf{a}$.
We use the notation
$\boldsymbol{i}, j$ for
the vector $\boldsymbol{i}$ with an extra
node index $j$ appended
at position $m+1$. We used ${\cal N}^c_p(\boldsymbol{x}^\mathsf{a})$ to denote the set of all $(m+1)$th order connected motifs that can be obtained by linking a new $c$-node to the $p$th position in motif $\boldsymbol{x}^\mathsf{a}$, i.e.,
\begin{equation}\label{eq:nbhdef}
	{\cal N}_{p}^c(\boldsymbol{x}^\mathsf{a}):=\bigcup_{\ell=1}^{m+1}
	\left\{\boldsymbol{y}^\mathsf{b}:|\boldsymbol{y}|=m+1, y_\ell=c, \boldsymbol{y}^\mathsf{b}_{\ell\shortto\emptyset}=\boldsymbol{x}^\mathsf{a},(\ell,p)\in\mathsf{b} \right\},
\end{equation}
where $\boldsymbol{y}^\mathsf{b}_{\ell\shortto\emptyset}$ denotes the $m$th order connected motif obtained by deleting the $\ell$th node of $\boldsymbol{y}^\mathsf{b}$. The sum over elements of ${\cal N}^c_p(\cdot)$ in \eqref{eq:split} is taken because any of the other motif nodes can also link to the new node. The types of
higher-order motifs appearing in the differential equation depend
on the considered motif. In Figure~\ref{fig:Dependence-of-motif},
we show this dependence structure (ignoring the labels). 

With the
substitutions from above, we obtain the 
\emph{individual-level moment equations}:
\begin{eqnarray}\label{eq:ddtxai}
	\ddt \langle\left[\boldsymbol{x}_{\boldsymbol{i}}^{\mathsf{a}}\right]\rangle & = & \sum_{p}^{m}\sum_{k}^{n}\biggl\{\Bigl(\sum_{c}^{n}\boldsymbol{\delta}_c(\boldsymbol{x})\mathsf{a}\boldsymbol{e}_p R^1_{kx_{p}c}+R^0_{kx_{p}}\Bigr)\bigl\langle\left[(\boldsymbol{x}_{\boldsymbol{i}}^{\mathsf{a}})_{p\shortto k}\right]\bigr\rangle+\sum_{c}^{n}\sum_{\boldsymbol{y}^\mathsf{b}\in{\cal N}^c_p\left(\boldsymbol{x}^{\mathsf{a}}_{p\shortto k}\right)}\,\sum_{j\notin\boldsymbol{i}}^{N}R^1_{kx_{p}c}\bigl\langle\left[\boldsymbol{y}_{\boldsymbol{i}, j}^{\mathsf{b}}\right]\bigr\rangle\nonumber \\
	&  & -\Bigl(\sum_{c}^{n}\boldsymbol{\delta}_c(\boldsymbol{x})\mathsf{a}\boldsymbol{e}_p R^1_{x_{p}kc}+R^0_{x_{p}k}\Bigr)\bigl\langle\left[\boldsymbol{x}_{\boldsymbol{i}}^{\mathsf{a}}\right]\bigr\rangle-\sum_{c}^{n}\sum_{\boldsymbol{y}^\mathsf{b}\in{\cal N}^c_p\left({\boldsymbol{x}^{\mathsf{a}}}\right)}\sum_{j\notin\boldsymbol{i}}^{N}R^1_{x_{p}kc}\bigl\langle\left[\boldsymbol{y}_{\boldsymbol{i}, j}^{\mathsf{b}}\right]\bigr\rangle\biggr\}.
\end{eqnarray}
The \emph{population-level moment equations} are then obtained by
taking the sum $\sum_{\boldsymbol{i}\in S(m,N)}$ [see definitions of $[\boldsymbol{x}^\mathsf{a}]$ and $[\boldsymbol{x}^\mathsf{a}_{\boldsymbol{i}}]$ in \eqref{eq:motifcount:all} and \eqref{eq:motif1}] of \eqref{eq:ddtxai} over all index sets $\boldsymbol{i}$:
\begin{equation}\label{eq:ddtxa}
	\begin{aligned}
		\ddt \langle\left[\boldsymbol{x}^{\mathsf{a}}\right]\rangle =&
		\sum_p^m\sum_k^n\sum_c^n\Biggl\{\left(\frac{R^0_{kx_{p}}}{n}+\boldsymbol{\delta}_c(\boldsymbol{x})\mathsf{a}\boldsymbol{e}_p R^1_{kx_{p}c}\right)\bigl\langle\left[(\boldsymbol{x}^{\mathsf{a}})_{p\shortto k}\right]\bigr\rangle-
		\left(\frac{R^0_{x_{p}k}}{n}+\boldsymbol{\delta}_c(\boldsymbol{x})\mathsf{a}\boldsymbol{e}_p R^1_{x_pkc}\right)\langle\left[\boldsymbol{x}^{\mathsf{a}}\right]\bigr\rangle
		\\
		&+\sum_{\boldsymbol{y}^\mathsf{b}\in{\cal N}^c_p\left(\boldsymbol{x}^{\mathsf{a}}_{p\shortto k}\right)}
		R^1_{kx_{p}c}\langle[\boldsymbol{y}^\mathsf{b}]\rangle-
		\sum_{\boldsymbol{y}^\mathsf{b}\in{\cal N}^c_p(\boldsymbol{x}^\mathsf{a})}
		R^1_{x_{p}kc}\langle[\boldsymbol{y}^\mathsf{b}]\rangle\Biggr\},
	\end{aligned}
\end{equation}
where we collected the sums at the front.
In closing this section, we make the following remarks. \emph{i}. The evolution of expected motif counts of size $m$ is a function of expected motif counts of size $m$ and $m+1$. \emph{ii}. Equation \eqref{eq:ddtxa} leads to the same equation for motifs in the network that are \emph{isomorphic}, where we call two motifs $\boldsymbol{x}^\mathsf{a}$ and $\boldsymbol{y}^\mathsf{b}$ isomorphic if there exists a permutation $\pi$ of their node indices that maps them onto each other, i.e. $\forall\boldsymbol{x}^{\mathsf{a}},\boldsymbol{y}^{\mathsf{b}}:\boldsymbol{x}^{\mathsf{a}}\simeq\boldsymbol{y}^{\mathsf{b}}\iff\exists \pi:\pi\boldsymbol{x}=\boldsymbol{y},\,\pi\mathsf{a}\pi^T=\mathsf{b}$. Hence, different motifs that are isomorphic belong to the same equivalence class. Isomorphic motifs have equal total counts such that we may consider only one representative from each equivalence class. In our implementation, we therefore make sure that we count all isomorphic motifs under a single representative node indexing. \emph{iii}. The summation in \eqref{eq:motifcount:all} over index tuples $\boldsymbol{i}\in S(m,N)$ that we used to go from \eqref{eq:ddtxai} to \eqref{eq:ddtxa} leads to multiple counting of motifs that possess \emph{automorphisms} other than the identity transformation. A motif has an automorphism if it can be mapped onto itself by a permutation of its node indices. An automorphism is therefore an isomorphism with itself, i.e. $\mathrm{Aut}(\boldsymbol{x}^{\mathsf{a}}):=\{ \pi:\pi\boldsymbol{x}=\boldsymbol{x},\,\pi\mathsf{a}\pi^T=\mathsf{a}\}$. The multiplicity with which a particular motif $\boldsymbol{x}^{\mathsf{a}}$ is counted (via \eqref{eq:motifcount:all}) is then equal to $|\mathrm{Aut}(\boldsymbol{x}^{\mathsf{a}})|$.

\paragraph{Conservation relations\label{sec:Conservation-equations}}

In fixed networks, the frequencies of induced subgraphs of a given
type (e.g. nodes, links, polygons, chains) remain fixed. The sum over all
possible motif label orderings then yields these fixed frequencies, or using
the notation from above,
\begin{equation}\label{eq:alginv1}
	\sum_{\boldsymbol{x}}\left[\boldsymbol{x}^{\mathsf{a}}\right]=\left[\mathsf{a}\right],\qquad\sum_{\boldsymbol{x}}\llbracket\boldsymbol{x}^{\mathsf{a}}\rrbracket=1,
\end{equation}
where the right is the normalised form of the left [via
\eqref{eq:norm}]. For networks with homogeneous degree and motifs and every stub $p$ of motif $\boldsymbol{x}^\mathsf{a}$ (also called leaf), i.e. a node with degree $1$, there is the additional conservation relation
\begin{equation}\label{eq:alginv2}
	\sum_{k}\left[\boldsymbol{x}^{\mathsf{a}}_{p\shortto k}\right]=\frac{\left[\mathsf{a}\right]}{\left[\mathsf{a}\setminus p\right]}\left[\boldsymbol{x}^\mathsf{a}_{p\shortto\emptyset}\right],\qquad\sum_{k}\llbracket\boldsymbol{x}^{\mathsf{a}}_{p\shortto k}\rrbracket=\llbracket\boldsymbol{x}^\mathsf{a}_{p\shortto\emptyset}\rrbracket,
\end{equation}
where the right is again the normalised form 
of the left (via \eqref{eq:norm}).
Under the conditions mentioned above, $\left[\mathsf{a}\right]/\left[\mathsf{a}\setminus p\right]$
is the number of out-motif connections at the node that connects to
the stub. As \eqref{eq:alginv1} and \eqref{eq:alginv2} can be derived for each motif up to the chosen
truncation order, they form an additional system of equations that can be used to
reduce the dimensionality of the mean-field model via elimination. To avoid multiplicity of 
equations while deriving relations \eqref{eq:alginv2}, we set up one equation per set of stubs that
lead to isomorphic variants of the considered subgraph when their indices are permuted.

\paragraph{Number of equations} A simple lower bound on the number of equations (before elimination) 
can be found by considering only chains. At each order
there is only one chain graph. Order 1 contributes $n$ equations, where $n$ is the number of  possible node states.  Each chain
of order $m>1$ has $\left(\binom{n}{m}\right)=\binom{m+n-1}{n}=\frac{(m+n-1)!}{(m-1)!n!}$
ways of labelling its $m$ nodes with $n$ species. The total number
of equations from chains is then $n_{\mathrm{ch}}=2+\sum_{m=2}^{k}\left(\binom{n}{m}\right).$
For networks in which the number of short cycles goes to zero with
$N\rightarrow\infty$, as in Erd\H{o}s-R{\'e}nyi random networks, the total
number of equations is equal to this lower bound plus the number of equations due to non-chain trees. When cycles need
to be taken into account however, there will be additional
connected motif types at each order; see Table~\ref{tab:Number-of-equations} column $n_{g}$ \cite[Table 4.2.1]{Harary1973}.
To obtain the number of equations contributed by each of these, one
has to consider all labelling orderings, knowing that some orderings
lead to isomorphic motifs and hence do not add to the total. We have listed in Table \ref{tab:Number-of-equations} the total number
of equations $n_{\mathrm{eq}}$ resulting from our enumeration algorithm
for dynamics with $n=2$, such as SIS epidemic spreading.
For a particular network, these are still reduced by 
the number of motifs not occurring in the considered network and by
the number of variables via the conservation relations. Column 
$n_4$ in Table~\ref{tab:Number-of-equations} shows the number of equations for the
square lattice and column $n_{4c}$ shows the remaining number after eliminating
variables by using conservation relations.

\begin{table}
	\centering
	\begin{tabular}{cccccc}
		\hline 
		$k$ & $n_{\mathrm{ch}}$ & $n_{g}$ & $n_{\mathrm{eq}}$ & $n_{\mathrm{4}}$ & $n_{\mathrm{4c}}$\tabularnewline
		\hline 
		$1$ & 2 & 1 & 2 & 2 & 1\tabularnewline
		$2$ & 5 & 2 & 5 & 5 & 2\tabularnewline
		$3$ & 11 & 4 & 15 & 11 & 4\tabularnewline
		$4$ & 21 & 10 & 65 & 35 & 14\tabularnewline
		$5$ & 36 & 31 & 419 & 113 & 38\tabularnewline
		\hline 
	\end{tabular}
	\caption{Cumulative number of equations ($n_{\mathrm{ch}}$, $n_{\mathrm{eq}}$, $n_{\mathrm{4c}}$) or motif
		types ($n_{g}$) as a function of order $k$. $n_{\mathrm{ch}}$: number of chain motifs (if the number of species $n=2$), $n_{g}$: number of subgraph types, $n_{\mathrm{eq}}$: total number of equations ignoring conservation relations (if $n=2$), $n_{\mathrm{4}}$: number of equations for the square lattice (if $n=2$), $n_{\mathrm{4c}}$: number of equations for the square lattice after elimination via	conservation relations (if $n=2$). \label{tab:Number-of-equations}}
\end{table}

\section{Moment equations for SIS spreading up to order 3}\label{sec:MFSIS3}
The moment equations up to third order can be written (via \eqref{eq:xdotblock}) as
\begin{equation}\label{eq:xdotblock3sis}
	\begin{blockarray}{ccccc}
		& \mathbf{x}_1 & \mathbf{x}_2 & \mathbf{x}_3 & \mathbf{x}_4\\ 
		\begin{block}{c[cccc]}
			\dot{\mathbf{x}}_1 & \mathbf{Q}_1 & \mathbf{Q}_{12} & \mathbf{0} & \mathbf{0}\\
			\dot{\mathbf{x}}_2 & \mathbf{0} & \mathbf{Q}_2 & \mathbf{Q}_{23} & \mathbf{0}\\
			\dot{\mathbf{x}}_3 & \mathbf{0} & \mathbf{0} & \mathbf{Q}_3 & \mathbf{Q}_{34}\\
		\end{block}
	\end{blockarray}
\end{equation}
Then, the coefficients for motifs up to order three are (omitting zero blocks)
\begin{equation*}\label{eq:SIS3}
	{\tiny
		\makeatletter\setlength\BA@colsep{.2pt}\makeatother
		\begin{blockarray}{cccccccccccc}
			& [\motI] & [\motIS] & [\motII] & [\motSIS] & [\motSSI] & [\motSII] & [\motISI] & [\motIII] & [\motSSItr] & [\motSIItr] & [\motIIItr]\\
			\begin{block}{cc|cc|cccccccc|}
				\dot{[\motI]} & -\gamma  & \beta  & 0 & & & & & & & & \\\cline{2-12}
				\dot{[\motIS]} & & -\beta -\gamma  & \gamma  & 0 & \beta  & 0 & -\beta  & 0 & \beta  & -\beta  & 0 \\
				\dot{[\motII]} & & 2 \beta  & -2 \gamma  & 0 & 0 & 0 & 2 \beta  & 0 & 0 & 2 \beta  & 0 \\\cline{2-12}
				\dot{[\motSIS]} & & & & -2 \beta -\gamma  & 0 & 2 \gamma  & 0 & 0 & 0 & 0 & 0 \\
				\dot{[\motSSI]} & & & & 0 & -\beta -\gamma  & \gamma  & \gamma  & 0 & 0 & 0 & 0 \\
				\dot{[\motSII]} & & & & \beta  & \beta  & -\beta -2 \gamma  & 0 & \gamma  & 0 & 0 & 0 \\
				\dot{[\motISI]} & & & & 0 & 0 & 0 & -2 \beta -2 \gamma  & \gamma  & 0 & 0 & 0 \\
				\dot{[\motIII]} & & & & 0 & 0 & 2 \beta  & 2 \beta  & -3 \gamma  & 0 & 0 & 0 \\
				\dot{[\motSSItr]} & & & & 0 & 0 & 0 & 0 & 0 & -2 \beta -\gamma  & 2 \gamma  & 0 \\
				\dot{[\motSIItr]} & & & & 0 & 0 & 0 & 0 & 0 & 2 \beta  & -2 \beta -2 \gamma  & \gamma  \\
				\dot{[\motIIItr]} & & & & 0 & 0 & 0 & 0 & 0 & 0 & 6 \beta  & -3 \gamma  \\
			\end{block}
		\end{blockarray},
	}
\end{equation*}
while those for fourth-order motifs in $\mathbf{Q}_{34}$ are
\begin{equation*}
	{\tiny
		\makeatletter\setlength\BA@colsep{.2pt}\makeatother
		\begin{blockarray}{ccccccccccccccccccccccccccccc}
			& [\motSSSI] & [\motSSSItre] & [\motSSSIsqo] & [\motSSSIstr] & [\motSISSsqi] & [\motSISI] & [\motSISIsqo] & [\motSISIstr] & [\motSISIsqi] & [\motISSI] & [\motSSIItre] & [\motSSIIsqo] & [\motISSIstr] & [\motSSIIstr] & [\motIISSsqi] & [\motIISI] & [\motSIIIsqo] & [\motIISIstr] & [\motIISIsqi] & [\motISIItre] & [\motISIIstr] & [\motIIISsqi] & [\motISSSstr] & [\motISSSsqi] & [\motSSSIsqii] & [\motISISsqi] & [\motSSIIsqii] & [\motSIIIsqii]\\\cline{2-29}
			\begin{block}{c|cccccccccccccccccccccccccccc}
				\dot{[\motSSS]} & -2 \beta  & -\beta  & -2 \beta  & -4 \beta  & -3 \beta  & 0 & 0 & 0 & 0 & 0 & 0 & 0 & 0 & 0 & 0 & 0 & 0 & 0 & 0 & 0 & 0 & 0 & 0 & 0 & 0 & 0 & 0 & 0 \\
				\dot{[\motSIS]} & 0 & \beta  & 0 & 2 \beta  & \beta  & -2 \beta  & -2 \beta  & -2 \beta  & -2 \beta  & 0 & 0 & 0 & 0 & 0 & 0 & 0 & 0 & 0 & 0 & 0 & 0 & 0 & 0 & 0 & 0 & 0 & 0 & 0 \\
				\dot{[\motSSI]} & \beta  & 0 & \beta  & \beta  & \beta  & 0 & 0 & 0 & 0 & -\beta  & -\beta  & -\beta  & -2 \beta  & -\beta  & -2 \beta  & 0 & 0 & 0 & 0 & 0 & 0 & 0 & 0 & 0 & 0 & 0 & 0 & 0 \\
				\dot{[\motSII]} & 0 & 0 & 0 & 0 & 0 & \beta  & \beta  & \beta  & \beta  & 0 & \beta  & 0 & \beta  & \beta  & \beta  & -\beta  & -\beta  & -\beta  & -\beta  & 0 & 0 & 0 & 0 & 0 & 0 & 0 & 0 & 0 \\
				\dot{[\motISI]} & 0 & 0 & 0 & 0 & 0 & 0 & 0 & 0 & 0 & 2 \beta  & 0 & 2 \beta  & 2 \beta  & 0 & 2 \beta  & 0 & 0 & 0 & 0 & -\beta  & -2 \beta  & -\beta  & 0 & 0 & 0 & 0 & 0 & 0 \\
				\dot{[\motIII]} & 0 & 0 & 0 & 0 & 0 & 0 & 0 & 0 & 0 & 0 & 0 & 0 & 0 & 0 & 0 & 2 \beta  & 2 \beta  & 2 \beta  & 2 \beta  & \beta  & 2 \beta  & \beta  & 0 & 0 & 0 & 0 & 0 & 0 \\
				\dot{[\motSSStr]} & 0 & 0 & 0 & 0 & 0 & 0 & 0 & 0 & 0 & 0 & 0 & 0 & 0 & 0 & 0 & 0 & 0 & 0 & 0 & 0 & 0 & 0 & -3 \beta  & -6 \beta  & -3 \beta  & 0 & 0 & 0 \\
				\dot{[\motSSItr]} & 0 & 0 & 0 & 0 & 0 & 0 & 0 & 0 & 0 & 0 & 0 & 0 & -2 \beta  & 0 & -2 \beta  & 0 & 0 & 0 & 0 & 0 & 0 & 0 & \beta  & 2 \beta  & \beta  & -2 \beta  & -2 \beta  & 0 \\
				\dot{[\motSIItr]} & 0 & 0 & 0 & 0 & 0 & 0 & 0 & 0 & 0 & 0 & 0 & 0 & 2 \beta  & 0 & 2 \beta  & 0 & 0 & 0 & 0 & 0 & -\beta  & -2 \beta  & 0 & 0 & 0 & 2 \beta  & 2 \beta  & -\beta  \\
				\dot{[\motIIItr]} & 0 & 0 & 0 & 0 & 0 & 0 & 0 & 0 & 0 & 0 & 0 & 0 & 0 & 0 & 0 & 0 & 0 & 0 & 0 & 0 & 3 \beta  & 6 \beta  & 0 & 0 & 0 & 0 & 0 & 3 \beta  \\
			\end{block}
		\end{blockarray}
	}.
\end{equation*} 	
When written out, this corresponds to the equations:
\begingroup
\allowdisplaybreaks
{\footnotesize 
	\begin{align} 
		\ddt \left[\motI\right] = & -\gamma\left[\motI\right]+\beta\left[\motIS\right],\label{eq:mf2}\\
		\ddt \left[\motIS\right] = & \gamma\left[\motII\right]-(\beta+\gamma)\left[\motIS\right]+\beta\left[\motSSI\right]+\beta\left[\motSSItr\right]-\beta\left[\motISI\right]-\beta\left[\motSIItr\right],\label{eq:mf4}\\
		\ddt \left[\motII\right]  = & 2\beta\left[\motIS\right]-2\gamma\left[\motII\right]+2\beta\left[\motISI\right]+2\beta\left[\motSIItr\right],\label{eq:mf5}\\
		\ddt \left[\motSIS\right]  = & 2\gamma\left[\motSII\right]-(2\beta+\gamma)\left[\motSIS\right]-2\beta\left[\motSISI\right]+\beta\left[\motSSSItre\right]-2\beta\left[\motSISIsqo\right]-2\beta\left[\motSISIstr\right]+2\beta\left[\motSSSIstr\right]-2\beta\left[\motSISIsqi\right]+\beta\left[\motSISSsqi\right],\label{eq:mf7}\\
		\ddt \left[\motSSI\right]  = & \gamma\left[\motSII\right]+\gamma\left[\motISI\right]-(\beta+\gamma)\left[\motSSI\right]-\beta\left[\motISSI\right]+\beta\left[\motSSSI\right]-\beta\left[\motSSIItre\right]-2\beta\left[\motISSIstr\right]-\beta\left[\motSSIIsqo\right]+\beta\left[\motSSSIsqo\right]\label{eq:mf8}\\
		& +\beta\left[\motSSSIstr\right]-\beta\left[\motSSIIstr\right]-2\beta\left[\motIISIsqi\right]+\beta\left[\motSISSsqi\right],\nonumber \\
		\ddt \left[\motSII\right]  = & \gamma\left[\motIII\right]-(\beta+2\gamma)\left[\motSII\right]+\beta\left[\motSIS\right]+\beta\left[\motSSI\right]-\beta\left[\motIISI\right]+\beta\left[\motSISI\right]+\beta\left[\motSSIItre\right]+\beta\left[\motISSIstr\right]-\beta\left[\motSIIIsqo\right]+\beta\left[\motSISIsqo\right]\label{eq:mf9}\\
		& -\beta\left[\motIISIstr\right]+\beta\left[\motSISIstr\right]
		+\beta\left[\motSSIIstr\right]+\beta\left[\motIISSsqi\right]-\beta\left[\motIISIsqi\right]+\beta\left[\motSISIsqi\right],\nonumber \\
		\ddt \left[\motISI\right]  = & \gamma\left[\motIII\right]-2(\beta+\gamma)\left[\motISI\right]+2\beta\left[\motISSI\right]-\beta\left[\motISIItre\right]-2\beta\left[\motISIIstr\right]+2\beta\left[\motISSIstr\right]+2\beta\left[\motSSIIsqo\right]-\beta\left[\motIIISsqi\right]+2\beta\left[\motIISSsqi\right],\label{eq:mf10}\\
		\ddt \left[\motIII\right]  = & -3\gamma\left[\motIII\right]+2\beta\left[\motSII\right]+2\beta\left[\motISI\right]+2\beta\left[\motIISI\right]+\beta\left[\motISIItre\right]+2\beta\left[\motISIIstr\right] +2\beta\left[\motSIIIsqo\right]+2\beta\left[\motIISIstr\right]\label{eq:mf11}\\
		& +\beta\left[\motIIISsqi\right]+2\beta\left[\motIISIsqi\right],\nonumber\\
		\ddt \left[\motSSItr\right]  = & 2\gamma\left[\motSIItr\right]-(2\beta+\gamma)\left[\motSSItr\right]-2\beta\left[\motISSIstr\right]+\beta\left[\motISSSstr\right]-2\beta\left[\motISISsqi\right]+2\beta\left[\motISSSsqi\right]
		-2\beta\left[\motIISSsqi\right]-2\beta\left[\motSSIIsqii\right]+\beta\left[\motSSSIsqii\right],\label{eq:mf13}\\
		\ddt \left[\motSIItr\right]  = & \gamma\left[\motIIItr\right]-2(\beta+\gamma)\left[\motSIItr\right]+2\beta\left[\motSSItr\right]-\beta\left[\motISIIstr\right]+2\beta\left[\motISSIstr\right]-2\beta\left[\motIIISsqi\right]+2\beta\left[\motISISsqi\right]+2\beta\left[\motIIISsqi\right],\label{eq:mf14}\\
		& -\beta\left[\motSIIIsqii\right]+2\beta\left[\motSSIIsqii\right],\nonumber\\
		\ddt \left[\motIIItr\right]  = & -3\gamma\left[\motIIItr\right]+6\beta\left[\motSIItr\right]+3\beta\left[\motISIIstr\right]+6\beta\left[\motIIISsqi\right]+3\beta\left[\motSIIIsqii\right].\label{eq:mf15}
	\end{align}
}
\endgroup
\newpage

\section{Obtaining steady states from simulations}\label{sec:feedback} 

\paragraph{Feedback control}
As we aim to compare the steady states
of mean-field models to those of the simulation, we need a way to obtain the steady
states of the simulations, even if they are unstable or marginally stable. 
Treating the simulation like an ideal physical experiment, the general approach to finding equilibria regardless of stability is to introduce a stabilizing feedback loop of the form
\begin{equation}
	\label{eq:control}
	r(t)=r_0+g([\boldsymbol{x}^\mathsf{a}](t)-[\boldsymbol{x}^\mathsf{a}]_\mathrm{ref})\mbox{,}
\end{equation}
as was done in
\cite{SGNWK08,schilder2015experimental,barton2017control} for
continuation of unstable vibrations in mechanical experiments.  In
\eqref{eq:control}, $r$ is one of the conversion rates in $\mathbf{R}^0$ or $\mathbf{R}^1$. 
The feedback control \footnote{\eqref{eq:control} is
	the simplest possible form of feedback control as it changes only
	one input ($r$) depending on one output ($[\boldsymbol{x}^\mathsf{a}]$) 
	and is static (no further processing
	of $[\boldsymbol{x}^\mathsf{a}]$ before feeding it back). Linear
	control theory ensures that unstable equilibria of ODEs can be
	stabilized with single-input-single-output dynamic feedback control
	using any single input and any single output satisfying some
	genericity conditions (linear controllability and
	observability).} makes
this rate time dependent by coupling it to the motif frequency
$[\boldsymbol{x}^\mathsf{a}](t)$ of a chosen motif
$\boldsymbol{x}^\mathsf{a}$ through the relation
\eqref{eq:control}. The factor $g$ is called the
feedback control gain and is problem specific.
When performing bifurcation analysis, it is convenient if the rate
$r$ used as the control input is also the bifurcation parameter
varied for the bifurcation diagram. In this case, whenever the
simulation with feedback control \eqref{eq:control} settles to an
equilibrium
$(r_\mathrm{c}^*,[\boldsymbol{x}^\mathsf{a}]_\mathrm{c}^*)$ (in
the limit of large $N$) the point
$(r_\mathrm{c}^*,[\boldsymbol{x}^\mathsf{a}]_\mathrm{c}^*)$ will
be on the equilibrium branch of the simulation without feedback control \cite{BS13,SOW14,renson2017experimental}.
For SIS spreading, we choose
\begin{equation}
	\label{eq:sis:control}
	r=\beta\mbox{,}\quad\boldsymbol{x}^\mathsf{a}=\motI\mbox{,}\quad g=\infty\mbox{.}
\end{equation}
This limit for feedback control results in what is called the conserved contact process, as proposed by \citet{Tome2001}. While SIS spreading does not have any unstable steady states,
points that are marginally stable, as near the
continuous phase transition (epidemic threshold) are stabilised as well by the control. This stabilisation suppresses fluctuations,
even close to the bifurcation, which results in faster convergence of the mean
and absence of absorption for any positive $[\motI]$. 

\begin{figure}
	\centering
	\includegraphics[width=1\columnwidth]{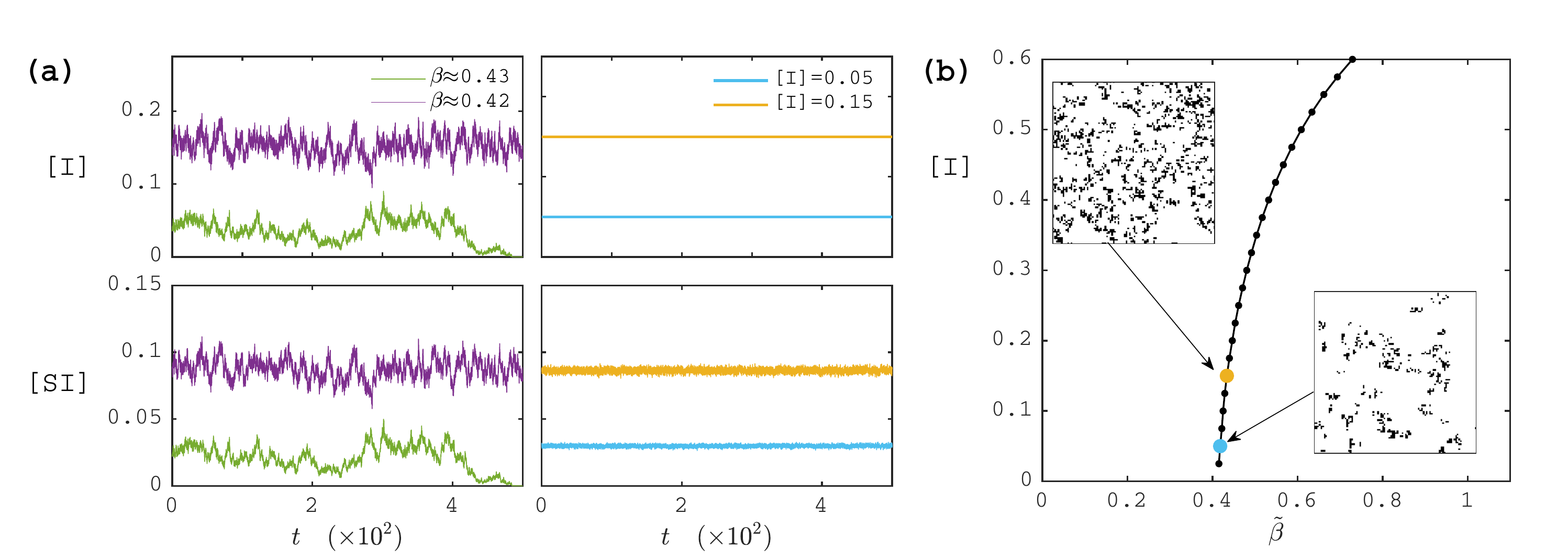}
	\caption{Feedback control of SIS epidemic spreading on a square lattice:
		(a) example time profiles of infected fraction \texttt{[I]} and
		fraction of susceptible-infected links
		\texttt{[SI]} in conventional (left) and controlled (right) simulation,
		(b) bifurcation diagram with control;
		insets: example spatial patterns;
		grid size 100$\times$100 (white=susceptible, black=infected), gain
		$g\gg1$, recovery rate $\gamma=1$. Near the critical point, 
		the conventional simulation shows large fluctuations and short 
		extinction times (absorption event occurs at the end
		of the time series for $\beta\approx0.43$ in green) whereas the
		controlled simulation has small fluctuations and no extinction.
		\label{fig:ctrl_t_bif}}
\end{figure}

\paragraph{Simulation algorithm} In the limit of infinite gain $g$, 
each recovery event forces a simultaneous
infection event, such that the number of infected nodes $[\motI]$ stays constant.
In a simulation based on the Gillespie algorithm \cite{Gillespie2007} this is
done in the following steps:
\begin{enumerate}[topsep=0pt,itemsep=-1ex,partopsep=1ex,parsep=1ex]
\item start with a number of randomly distributed infected
nodes, 
\item recover an infected node selected uniformly at
random, 
\item infect a randomly selected susceptible node,
with selection probability proportional to its number of infected
neighbours, 
\item advance time with
$\Delta t=-\log(\xi)/(\gamma\left[\motI\right])$, (where $\xi$ is a
uniform random variable on the interval $[0,1]$),
\item go to step $2$.
\end{enumerate}
This loop runs until we observe that $\left[\motIS\right](t)$ is stationary (call this time $t_\mathrm{e}$). Then for some additional time $T$ we observe the fluctuations of $\left[\motIS\right](t)$ around its mean. \citet{Tome2001} derived the effective
infection rate $\tilde{\beta}$ for each chosen count 
$\left[\motI\right]$ of infected nodes by noting that, because for every 
infection event there is a recovery event,
\[
\tilde{\beta}\langle\left[\motIS\right]\rangle=\gamma\left[\motI\right],
\]
where $\langle\cdot\rangle$ is an average over many independent realisations, such that
\[
\tilde{\beta}=\gamma\frac{\left[\motI\right]}{\langle\left[\motIS\right]\rangle},
\]
which they found to lead to the same nontrivial steady states 
$\left[\motI\right]^*(\tilde{\beta}/\gamma)$ as $\left[\motI\right]^*(\beta/\gamma)$
in the model without control.
As the model with control is ergodic (no absorbing states exist) we need to run only a single realisation and compute the effective infection rate as
\begin{equation}
\tilde{\beta}=\gamma\frac{\left[\motI\right]}{\int_{t_\mathrm{e}}^{t_\mathrm{e}+T}{\left[\motIS\right](t)\de t}/T}.\label{eq:beteff_tav}
\end{equation}
In this manner, the error bars 
of estimates of $\tilde\beta$ can be made arbitrarily small by increasing $T$.

\section{Correlations at given distance}\label{sec:dcorrs}

\begin{figure}[H]
\centering
\includesvg[inkscapelatex=false,width=.45\columnwidth]{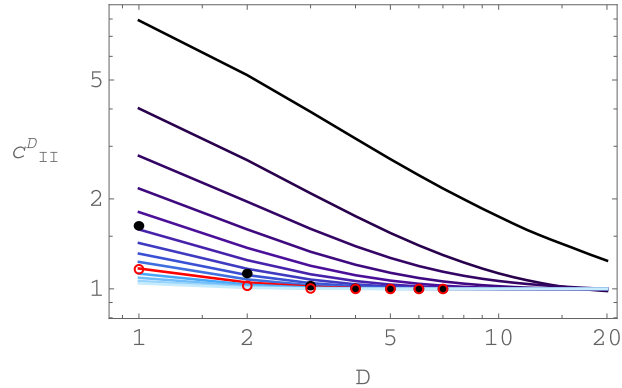}
\caption{Steady state correlation function $C_{II}^{D}$ versus distance $D$
	\eqref{eq:CabD} for simulations of SIS spreading on the square
	lattice (lines). Colour
	scale indicates values of $\beta/\gamma$ corresponding to the markers
	in Figure \ref{fig:Comparison-of-the} and sorted such that higher correlations
	occur for lower values of $\beta/\gamma$ (closer to the critical
	point). Markers show approximated correlations when assuming MF2 conditions
	hold (using \eqref{eq:CabDmots}) and factorising the numerator in
	pair fractions). Filled black markers are MF2 correlations for $\beta/\gamma=0.4137$
	(which for simulations is the top black line). Open red markers are
	MF2 correlations for $\beta/\gamma=0.5602$ (which for simulations
	is the red line). \label{fig:Steady-state-correlation}}
\end{figure}

In Figure \ref{fig:Steady-state-correlation} we show for SIS epidemic spreading the correlation at given distance between two nodes (only between infected nodes shown).
Its general definition is
\begin{equation}
C_{ab}^{D}=\frac{N^2}{\sum_{ij}\mathds{1}_{\d(i,j)=D}}\frac{\sum_{ij}[\mota_{i}\motb_{j}]_{\d(i,j)=D}}{[\mota][\motb]},\label{eq:CabD}
\end{equation}
or with normalised motifs:
\begin{equation}
C_{ab}^{D}=\frac{\sum_{ij}\llbracket\mota_{i}\motb_{j}\rrbracket_{\d(i,j)=D}}{\llbracket\mota\rrbracket\llbracket\motb\rrbracket},\label{eq:CabDn}
\end{equation}
where the distance $\d(i,j)$ is the length of the shortest path between $i$ and $j$. 
Note that the correlations between
neighbouring nodes (\ref{eq:corrs0}) is a special case of this, i.e.
$C_{ab}=C_{ab}^{1}$. An alternative way to write \eqref{eq:CabDn} is 
\begin{equation}
C_{ab}^{D}=\frac{\sum_{c_{1},...,c_{D-1}}\llbracket a c_{1}...c_{D-1}b\rrbracket}{\llbracket\mota\rrbracket\llbracket\motb\rrbracket},\label{eq:CabDmots}
\end{equation}
where $\llbracket a c_{1}...c_{D-1}b\rrbracket$ is a chain motif of size $D+1$
with indicated states. This definition
was used to derive $C_{ab}^{D}$ as approximated by mean-field
models in Figure \ref{fig:Steady-state-correlation}, by applying the closure formula to the chain in the numerator.

\section{Steady states of MF1-MF2}\label{sec:ss12}
Here, we show the expressions for the steady states and steady state correlations of MF1 and MF2.
\subsection{MF1}
The MF1 model equals SIS epidemic spreading under well-mixed conditions. Its steady states are the trivial and the endemic state,
\[
\llbracket\motI\rrbracket_{1}^{*}=0,\llbracket\motI\rrbracket_{2}^{*}=1-\frac{\gamma}{\kappa\beta}.
\]

\subsection{MF2}
For the MF2 model the equations and, hence, their steady states, depend on the type of network.

\subsubsection{Degree-homogeneous networks}
The steady states of the dynamic equations for MF2 on degree-homogeneous networks \eqref{eq:mf2:sis:1} are 
\begin{equation}
(\llbracket\motI\rrbracket^{*},\llbracket\motIS\rrbracket^{*})_{1}=(0,0),\qquad(\llbracket\motI\rrbracket^{*},\llbracket\motIS\rrbracket^{*})_{2}=(1-\frac{\kappa-1}{\kappa\beta(\kappa-1)/\gamma-1},\frac{\beta(\kappa-1)-1}{\beta(\beta\kappa(\kappa-1)-1)}).\label{eq:SS2homkap}
\end{equation}
The non-trivial correlations, obtained by substitution of the above into \eqref{eq:corrs}, are
\begin{equation}
C_{II}^{*}=\frac{(\gamma-\beta\kappa)(\gamma-\beta(\kappa-1)\kappa)}{\beta\kappa^{2}(\beta(\kappa-1)-\gamma)},\quad C_{SI}^{*}=1-\frac{\gamma}{\beta\kappa(\kappa-1)},\quad C_{SS}^{*}=\frac{\beta(\kappa-1)\kappa-\gamma}{\beta(\kappa-1)^{2}}.\label{eq:sscorrs2homkap}
\end{equation}

\subsubsection{Degree-heterogeneous networks}
The steady states of the dynamic equations for MF2 on degree-heterogeneous networks \eqref{eq:mf2:sis:het2} are 
\begin{eqnarray}
(\llbracket\motI\rrbracket^{*},\llbracket\motII\rrbracket^{*},\llbracket\motIS\rrbracket^{*})_{1} & = & (0,0,0),\nonumber \\
(\llbracket\motI\rrbracket^{*},\llbracket\motII\rrbracket^{*},\llbracket\motIS\rrbracket^{*})_{2} & = & \frac{\gamma  \left(\sqrt{\beta  (\kappa -1)^2+4 \gamma }-\sqrt{\beta } (\kappa +3)\right)}{2 \beta ^{3/2} \kappa }+1,\nonumber\\
&  & \frac{\gamma\sqrt{\beta}(\kappa+1)-\gamma\sqrt{4\gamma+\beta(\kappa-1)^{2}}}{2\kappa\beta^{3/2}}).\label{eq:SS2hetkap}
\end{eqnarray}
The non-trivial correlations, via substitution of the above into \eqref{eq:corrs}, are
\begin{align}\label{eq:sscorrs2hetkap}
C_{II}^{*}&=\frac{4\beta^{3/2}\kappa-2\sqrt{\beta}\gamma(\kappa+3)+2\gamma\sqrt{\beta(\kappa-1)^{2}+4\gamma}}{\sqrt{\beta}\kappa\left(\sqrt{\beta(\kappa-1)^{2}+4\gamma}-\sqrt{\beta}(\kappa+1)\right)^{2}},\\ C_{SI}^{*}&=C_{SS}^{*}=\frac{2\gamma}{\kappa\sqrt{\beta^{2}(\kappa-1)^{2}+4\beta\gamma}-\beta(\kappa-1)\kappa}.
\end{align}

\section{MF$3$ for SIS spreading}\label{sec:MF3}
We only apply MF3 to the square lattice, such that we
can ignore all motifs that contain triangles and set
$\left[\mottr\right]=\left[\motstr\right]=\left[\motsqi\right]=\left[\motsqii\right]=0$
in the remaining equations, such that we obtain 
\begin{equation*}
{\footnotesize
	\makeatletter\setlength\BA@colsep{.1pt}\makeatother
	\begin{blockarray}{cccccccccccccccccccc}
		& & [\motI] & [\motIS] & [\motII] & [\motSIS] & [\motSSI] & [\motSII] & [\motISI] & [\motIII]& [\motSSSI] & [\motSSSItre] & [\motSSSIsqo] & [\motSISI] & [\motSISIsqo] & [\motISSI] & [\motSSIItre] & [\motSSIIsqo] & [\motIISI] & [\motSIIIsqo] & [\motISIItre]\\
		\begin{block}{ccc|cc|ccccc|cccccccccc}
			\dot{[\motI]} & & -\gamma  & \beta  & 0 & & & & & & & & & & & & & & & & \\\cline{2-21}
			\dot{[\motIS]} & & & -\beta -\gamma  & \gamma  & 0 & \beta  & 0 & -\beta  & 0 & & & & & & & & & & & \\
			\dot{[\motII]} & & & 2 \beta  & -2 \gamma  & 0 & 0 & 0 & 2 \beta & 0 & & & & & & & & & & & \\\cline{2-21}
			\dot{[\motSIS]} & & & & & -2 \beta -\gamma  & 0 & 2 \gamma  & 0 & 0 & 0 & \beta  & 0 & -2 \beta  & -2 \beta & 0 & 0 & 0 & 0 & 0 & 0\\
			\dot{[\motSSI]} & & & & & 0 & -\beta -\gamma  & \gamma  & \gamma  & 0 & \beta  & 0 & \beta  & 0 & 0 & -\beta  & -\beta  & -\beta & 0 & 0 & 0\\
			\dot{[\motSII]} & & & & & \beta  & \beta  & -\beta -2 \gamma  & 0 & \gamma & 0 & 0 & 0 & \beta  & \beta & 0 & \beta  & 0 & -\beta & -\beta  & 0\\
			\dot{[\motISI]} & & & & & 0 & 0 & 0 & -2 \beta -2 \gamma  & \gamma & 0 & 0 & 0 & 0 & 0 & 2 \beta  & 0 & 2 \beta  & 0 & 0 & -\beta\\
			\dot{[\motIII]} & & & & & 0 & 0 & 2 \beta  & 2 \beta  & -3 \gamma & 0 & 0 & 0 & 0 & 0 & 0 & 0 & 0 & 2 \beta  & 2 \beta& \beta\\
		\end{block}
	\end{blockarray}.
}
\end{equation*}	
The 7 conservation relations are
\begingroup
\allowdisplaybreaks
\begin{eqnarray*}
[\motI]+[\motS] & = & N,\\\relax
[\motIS]+[\motSS] & = & \kappa[\motS],\\\relax
[\motII]+[\motIS] & = & \kappa[\motI],\\\relax
[\motSSI]+[\motSSS] & = & (\kappa-1)[\motSS],\\\relax
[\motSII]+[\motSIS] & = & (\kappa-1)[\motIS],\\\relax
[\motISI]+[\motSSI] & = & (\kappa-1)[\motIS],\\\relax
[\motIII]+[\motSII] & = & (\kappa-1)[\motII].
\end{eqnarray*}
\endgroup
We eliminate $[\motII]$, $[\motSIS]$,
$[\motSSI]$, $[\motSII]$.  Following (\ref{eq:subcons}-\ref{eq:xdotsubcons2}),
this means
\begin{equation*}
\mathbf{x}=
\begin{bmatrix}
	[\motI] \\ [\motIS] \\ [\motII] \\ [\motSIS] \\ [\motSSI] \\ [\motSII] \\ [\motISI] \\ [\motIII]\\
\end{bmatrix},
\;\;
\tilde{\mathbf{x}}=
\begin{bmatrix}
	[\motI] \\ [\motIS] \\ [\motISI] \\ [\motIII]\\
\end{bmatrix},
\;\;
\tilde{\mathbf{x}}_4=
\begin{bmatrix}
	[\motISSI] \\ [\motSSIIsqo] \\ [\motIISI] \\ [\motSIIIsqo] \\ [\motISIItre]\\
\end{bmatrix},
\end{equation*}
and
\begin{align*}
\tilde{\mathbf{Q}}_{1\cdots3,1\cdots3}&=
\setcounter{MaxMatrixCols}{20}
\left(\begin{smallmatrix}
	-\gamma  & \beta  & 0 &0 &0 &0 &0 &0 \\
	0& -\beta -\gamma  & \gamma  & 0 & \beta  & 0 & -\beta  & 0\\
	0& 2 \beta  & -2 \gamma  & 0 & 0 & 0 & 2 \beta & 0\\
	0&0 &0 & -2 \beta -\gamma  & 0 & 2 \gamma  & 0 & 0\\
	0&0 &0 & 0 & -\beta -\gamma  & \gamma  & \gamma  & 0\\
	0&0 &0 & \beta  & \beta  & -\beta -2 \gamma  & 0 & \gamma\\
	0&0 &0 & 0 & 0 & 0 & -2 \beta -2 \gamma  & \gamma\\
	0&0 &0 & 0 & 0 & 2 \beta  & 2 \beta  & -3 \gamma\\
\end{smallmatrix}\right),
&
\mathbf{E}&=
\left(\begin{smallmatrix}
	1 & 0 & 0 & 0 \\
	0 & 1 & 0 & 0 \\
	\kappa  & -1 & 0 & 0 \\
	-\kappa(\kappa-1) & 2 (\kappa -1) & 0 & 1 \\
	0 & \kappa -1 & -1 & 0 \\
	\kappa(\kappa -1)  & 1-\kappa  & 0 & -1 \\
	0 & 0 & 1 & 0 \\
	0 & 0 & 0 & 1 \\
\end{smallmatrix}\right),\\
\tilde{\mathbf{Q}}_{34}&=
\left(\begin{smallmatrix}
	2 \beta  & 2 \beta  & 0 & 0 & -\beta  \\
	0 & 0 & 2 \beta  & 2 \beta  & \beta  \\
\end{smallmatrix}\right),
&
\mathbf{c}&=
\left(\begin{smallmatrix}
	0&0&0&0&0&0&0&0\\
\end{smallmatrix}\right)^T,
\end{align*}
which results in the four remaining equations ($\tilde{\mathbf{c}}=\mathbf{0}$ so not shown)
\begin{equation}\label{eq:sqlatopen}
{\footnotesize
	\makeatletter\setlength\BA@colsep{1pt}\makeatother
	\begin{blockarray}{ccccccccccc}
		& & [\motI] & [\motIS] & [\motISI] & [\motIII] & [\motISSI] & [\motSSIIsqo] & [\motIISI] & [\motSIIIsqo] & [\motISIItre]\\
		\begin{block}{ccc|c|cc|ccccc}
			\dot{[\motI]} & & -\gamma & \beta  & & & & & & & \\\cline{2-11}
			\dot{[\motIS]} & & 4\gamma & 2\beta -2\gamma & -2\beta  & 0 & & & & & \\\cline{2-11}
			\dot{[\motISI]} & & & & -2 \beta -2 \gamma  & \gamma & 2 \beta  & 2 \beta  & 0 & 0 & -\beta\\
			\dot{[\motIII]} & &24\beta & -6\beta & 2 \beta  & -3 \gamma-2\beta & 0 & 0 & 2 \beta  & 2 \beta& \beta\\
		\end{block}
	\end{blockarray}.
}
\end{equation}
where we have used $\kappa=4$ for the square lattice. In normalised form this is
\begin{equation}
{\footnotesize
	\makeatletter\setlength\BA@colsep{1pt}\makeatother
	\begin{blockarray}{ccccccccccc}
		& & \llbracket\motI\rrbracket & \llbracket\motIS\rrbracket & \llbracket\motISI\rrbracket & \llbracket\motIII\rrbracket & \llbracket\motISSI\rrbracket & \llbracket\motSSIIsqo\rrbracket & \llbracket\motIISI\rrbracket & \llbracket\motSIIIsqo\rrbracket & \llbracket\motISIItre\rrbracket\\
		\begin{block}{ccc|c|cc|ccccc}
			\dot{\llbracket\motI\rrbracket} & & -\gamma & 4\beta  & & & & & & & \\\cline{2-11}
			\dot{\llbracket\motIS\rrbracket} & & \gamma & 2\beta -2\gamma & -6\beta  & 0 & & & & & \\\cline{2-11}
			\dot{\llbracket\motISI\rrbracket} & & & & -2 \beta -2 \gamma  & \gamma & \frac{14}{3}\beta  & \frac{4}{3}\beta  & 0 & 0 & -2\beta\\
			\dot{\llbracket\motIII\rrbracket} & &2\beta & -2\beta & 2 \beta  & -3 \gamma-2\beta & 0 & 0 & \frac{14}{3}\beta  & \frac{4}{3}\beta& 2\beta\\
		\end{block}
	\end{blockarray}.
}
\end{equation}
To close the system of equations, we apply \eqref{eq:dec:mot:chord} to the chains and star, and \eqref{eq:dec:mot:nchord} to 
the cycles (using the extension to non-maximal 2-cliques) and obtain the normalised closures 
(see also Section \ref{sec:clex}, examples 3-5)
\begingroup
\allowdisplaybreaks
\begin{align}\label{eq:sis:clos:SSSItre}
\llbracket\motISSI\rrbracket & \approx  \frac{\llbracket\motSSI\rrbracket^2}{\llbracket\motSS\rrbracket}=\frac{(\llbracket\motIS\rrbracket - \llbracket\motISI\rrbracket)^2}{(1 - \llbracket\motI\rrbracket - \llbracket\motIS\rrbracket)},\nonumber\\
\llbracket\motIISI\rrbracket & \approx  \frac{\llbracket\motSII\rrbracket\llbracket\motISI\rrbracket}{\llbracket\motIS\rrbracket}=\frac{(\llbracket\motI\rrbracket - \llbracket\motIII\rrbracket - \llbracket\motIS\rrbracket) \llbracket\motISI\rrbracket}{\llbracket\motIS\rrbracket},\nonumber\\
\llbracket\motISIItre\rrbracket & \approx  \frac{\llbracket\motIS\rrbracket^3}{\llbracket\motS\rrbracket^2}=\frac{\llbracket\motIS\rrbracket^3}{(1 -\llbracket\motI\rrbracket)^2},\\
\llbracket\motSSIIsqo\rrbracket & \approx  \frac{\llbracket\motSII\rrbracket^2\llbracket\motSSI\rrbracket^{2}}{\llbracket\motIS\rrbracket^2\llbracket\motSS\rrbracket\llbracket\motII\rrbracket}=\frac{(\llbracket\motI\rrbracket - \llbracket\motIII\rrbracket - \llbracket\motIS\rrbracket)^2 (\llbracket\motIS\rrbracket - \llbracket\motISI\rrbracket)^2}{(1 - \llbracket\motI\rrbracket - \llbracket\motIS\rrbracket) (\llbracket\motI\rrbracket - \llbracket\motIS\rrbracket) \llbracket\motIS\rrbracket^2},\nonumber\\
\llbracket\motSIIIsqo\rrbracket & \approx  \frac{\llbracket\motSII\rrbracket^{2}\llbracket\motIII\rrbracket\llbracket\motISI\rrbracket}{\llbracket\motIS\rrbracket^2\llbracket\motII\rrbracket^2}=\frac{\llbracket\motIII\rrbracket (\llbracket\motI\rrbracket - \llbracket\motIII\rrbracket - \llbracket\motIS\rrbracket)^2 \llbracket\motISI\rrbracket}{(\llbracket\motI\rrbracket - \llbracket\motIS\rrbracket)^2 \llbracket\motIS\rrbracket^2},\nonumber
\end{align}
\endgroup
where we also used the conservation relations to substitute any previously eliminated motif.
The steady state solutions of the final system
are roots of a ninth-order polynomial, of which two are admissible (see Figure~\ref{fig:Comparison-of-the} for numerical results).

\section{MF4 for SIS spreading in \texttt{Mathematica}\label{sec:MFhoeqs}}

We uploaded two application examples of our algorithm at fourth order to \href{https://notebookarchive.org/moment-closure-and-moment-equations-for-sis-spreading--2022-06-6hd82es/}{\texttt{Mathematica}'s
notebook archive} (best viewed locally).
In the 
first section of the file, the procedure as explained in the main text is followed for MF4 on the square lattice.
The second section in the file shows the general coefficient matrix $\mathbf{Q}$ for the (unclosed) moment equations up to fourth order, with columns labelled by $\mathbf{x}_1,\dots,\mathbf{x}_5$ (stored in the 
variable \texttt{mots}). 
The algorithm can be instructed to write the resulting set of equations as \texttt{MATLAB}
\cite{MathWorks2021} functions in a format compatible
with continuation software such as \texttt{COCO} \cite{COCO} and
\texttt{MatCont} \cite{matcont}. 

\section{Closure examples}\label{sec:clexs}

Here we show 13 examples of subgraph decompositions based on the
method explained in Section \ref{sec:Closure-scheme} (\ref{eq:dec:mot:chord},~\ref{eq:dec:mot:nchord}). They are shown in table form 
on the next page
but we also discussed three commonly used ones (1-3 in the table) in 
Section \ref{sec:clex}.

As the closures can be
written independent of the particular labels, they are shown in the table for subgraphs only, with each node tagged with its index. We have
also dropped the $\langle\cdot\rangle$, assuming that the law of large numbers applies, such that the counts approach their expectations almost surely for increasing network size $N$. The examples can be understood 
by reading the table from left to right. Comments are added in the last column. By column, it shows: 1. the 
example number, 2. the considered subgraph, 3. its diameter,
4. its independence map assuming independence beyond distance $\text{diam}(\mathsf{a}-1)$, 5. whether the independence map is chordal,
6. the derived junction graph of $d$-cliques where $d=\text{diam}(\mathsf{a}-1)$, 7. the resulting
closure formula, 8. a triangulation of the independence map if the independence map
is non-chordal, 9. the junction graph
of $d$-cliques based on the triangulation, 10. the closure formula based on the 
triangulation, 11. an ad-hoc extension of the method to non-maximal cliques for some
subgraphs, 12. the resulting closure formula based on this extension, 13. comments. 
Our \texttt{Mathematica} script \cite{Wuyts2022mfmcl} generates these closures automatically. 

\newpage
\eject \pdfpagewidth=70cm \pdfpageheight=60cm
\includegraphics{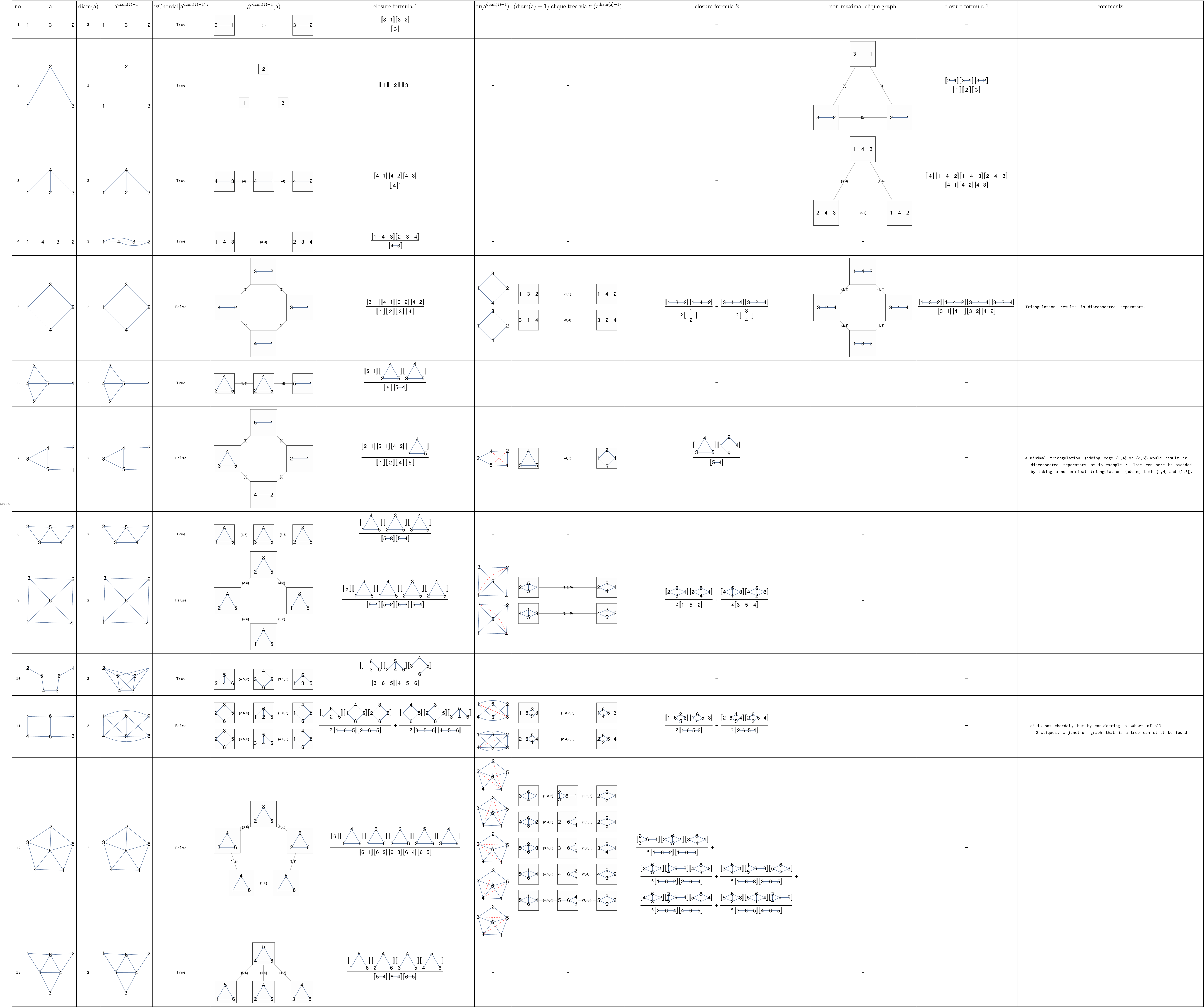}

\end{document}